\documentclass[pdflatex,sn-mathphys-num]{sn-jnl}

\usepackage{graphicx}%
\usepackage{multirow}%
\usepackage{amsmath,amssymb,amsfonts}%
\usepackage{amsthm}%
\usepackage{mathrsfs}%
\usepackage[title]{appendix}%
\usepackage{xcolor}%
\usepackage{textcomp}%
\usepackage{manyfoot}%
\usepackage{booktabs}%
\usepackage{algorithm}%
\usepackage{algorithmicx}%
\usepackage{algpseudocode}%
\usepackage{listings}%

\usepackage[utf8]{inputenc}
\usepackage{bbm}
\usepackage{hyperref}
\usepackage{cleveref}
\usepackage{mathtools}
\usepackage[dvipsnames]{xcolor}
\usepackage{caption}
\usepackage{placeins}
\usepackage{diagbox}
\usepackage{tabularx}
\usepackage{cancel}
\usepackage[symbol]{footmisc}
\usepackage[T1]{fontenc}

\usepackage{scalerel}
\usepackage{stackengine,wasysym}

\usepackage{tikz}

\numberwithin{equation}{section}%

\newcommand{\E}{\mathbbm{E}}
\newcommand{\R}{\mathbbm{R}}

\newcommand{\bcdot}{\ensuremath{\boldsymbol{\cdot}}}

\newcommand{\tr}{{\scriptscriptstyle \rm T}}

\theoremstyle{thmstyleone}%
\theoremstyle{thmstyletwo}%
\newtheorem{remark}{Remark}%

\theoremstyle{thmstylethree}%

\begin{document}

\title[Stochastic interpretation of the Water Waves equations]{Stochastic interpretation of the Water Waves equations}


\author*[1]{\fnm{Antoine} \sur{Moneyron}}\email{antoine.moneyron@inria.fr}

\affil[1]{Univ Rennes, INRIA, IRMAR UMR 6625, F35000, Rennes, France}


\abstract{The aim of this article is to derive a stochastic interpretation of the Water Waves equations, based on the location uncertainty principle (LU). These equations can be used to infer simpler models from scaling limits, such as the Saint-Venant, Boussinesq and Serre-Green-Naghdi ones. Once adapted to the stochastic setting, these models have been successfully implemented numerically, and we show several experiments to assess their properties -- for example the breaking of the wave symmetry, and the influence of time correlated noises.}

\keywords{water wave equations, stochastic modelling, non-linear stochastic partial differential equations, numerical simulations}



\maketitle

\section{Introduction}\label{sec:intro}

Understanding the ocean surface waves is crucial to decipher the dynamical behaviour of the ocean, due to energy transfers between the waves and the current. Yet, these phenomena follow different dynamics since they operate different time and spatial scales. In addition, they are modelled through dramatically different assumptions: wave dynamics are usually derive from potential -- or irrotational -- flows, the current is represented by velocity with non-zero vorticity. Consequently, the wave-current interaction is very complex to model numerically. For this purpose, designing simplified stochastic models of surface waves would be advantageous to account for both dynamics. This can be achieved by adding as a specific noise term in the stochastic representation of ocean circulation \cite{LDLM_2023, TML_2023}.

Several approaches have been explored to adapt the classical deterministic derivation to the stochastic setting in order to derive such stochastic wave models. For instance, deep-ocean stochastic long waves have been shown to emerge from a linearised stochastic shallow water system \cite{MLLTC_2023}, as well as from a stochastic Hamiltonian formulation \cite{Dinvay-Memin-PRSA-22}. Water waves travelling from deep-water regions into shallower areas — where the water depth is much smaller than their wavelength — undergo significant transformations. Rapid changes in height, velocity and direction lead to significant modifications of the free surface profiles. Starting out as nearly perfect sinusoidal shapes, these profiles gradually evolve into an asymmetric form. Mathematically, coastal waves are described by a range of dynamical models, each corresponding to different approximations of the irrotational Euler equations, ultimately averaged over the water column.

More specifically, we rely on the location uncertainty principle developed in \cite{Mémin_2014}, which is based on adding a noise term to the Lagrangian formulation of the displacement. This approach is in fact a central aspect of stochastic ocean modelling and has been widely discussed: an SDE-based stochastic generalisation of the deterministic Lagrangian expression of the flow was proposed in \cite{Restrepo_2007,RRmWB_2011}, from which the authors derived a corresponding Eulerian expression. A similar idea also appears in turbulence modelling: McWilliam and Berloff, for example, proposed a stochastic parametrisation built on Langevin models of turbulence in \cite{BmW_2002}, developed in the wake of Kraichnan's foundational work \cite{Kraichnan59}.

In this work, we aim to derive a location uncertainty interpretation of nonlinear coastal wave models, building on the modelling principle introduced in \cite{Mémin_2014}. These models are naturally accompanied by nonlinear variability arising from numerical or modelling uncertainties. More fundamentally, this stochastic formalism offers a way to account for the effect of the unresolved vertical turbulent component. In this work, we focus on the family of models associated with the Serre-Green-Naghdi equations \cite{Serre_1953, GN_1976}, which capture non-hydrostatic phenomena related to vertical acceleration. In particular, we aim to assess the behaviour of the corresponding numerical simulations. Our main findings show that the LU interpretations of the Saint-Venant, Boussinesq and Serre-Green-Naghdi wave models break the symmetry of the wave and introduce variability relative to their deterministic counterparts. Nevertheless, they still satisfy the same conservation principles, subject to mild conditions on the noise structure. The LU framework thus offers a way to explore the influence of "conservative randomness" on the classical dynamics of wave models. We further provide graphical illustrations of this influence based on numerical simulations.

The paper is organised as follows: after an overview of the location uncertainty (LU) principle, we derive a stochastic representation of the Serre-Green-Naghdi shallow water wave model, following the same physical approximations used in the deterministic setting \cite{Barthelemy-04}. We then show how simpler models can be obtained through further approximations, and discuss the conservation of standard quantities such as mass, momentum and energy. The final section presents and evaluates numerical results for the different models.


\section{The LU formalism} \label{section:prelim-LU}

The LU formulation rests on the following separation of time scales in the flow:
\begin{equation}
    dX_t = u(X_t,t) dt + \sigma(X_t,t) dW_t.
\end{equation}
This decomposition, understood in the It\={o} sense, defines $X$ as the Lagrangian displacement over a three-dimensional domain $S \subset \mathbb{R}^3$. The velocity term $u(X_t, t)$ represents the large-scale Eulerian velocity (correlated in both space and time), while $\sigma(X_t, t)dW$ captures a rapidly oscillating, unresolved velocity component (uncorrelated in time but correlated in space). This second term will be treated as noise in the following.

Let us specify this noise term more precisely: consider a Wiener process $W$ defined on the space of square integrable functions $\mathcal{W} := L^2(\mathcal{S},\R^3)$. Then, there exists a Hilbert orthonormal basis $(e_i)_{i \in \mathbbm{N}}$ of $\mathcal{W}$ and a sequence of independent standard Brownian motions $(\hat{\beta}^i)_{i \in \mathbbm{N}}$, defined on a filtered probability space $(\Omega, \mathcal{F}, (\mathcal{F}_t)_t, \mathbbm{P})$, such that,
\begin{equation*}
    W=\sum_{i \in \mathbbm{N}} \hat{\beta}^i e_i.
\end{equation*}
The tuple $(\Omega, \mathcal{F}, (\mathcal{F}_t)_t, \mathbbm{P}, W)$ thus forms a stochastic basis. Notice that the sum $\sum_{i \in \mathbbm{N}} \hat{\beta}^i e_i$ does \emph{not} converge in $\mathcal{W}$. Hence, the previous identity only makes sense in a space $\mathcal{U}$ including $\mathcal{W}$, chosen so that the embedding $\mathcal{W} \hookrightarrow \mathcal{U}$ is Hilbert-Schmidt. For instance, $\mathcal{U}$ can be the dual space of any reproducing kernel Hilbert subspace
of $\mathcal{W}$ for the inner product $(\bcdot, \bcdot)_{\mathcal{W}}$, e.g. $H^{-s}(\mathcal{S}, \R^3)$ with $s > \frac{3}{2}$. The noise is then defined using a deterministic time dependent correlation operator $\sigma_t$: let $\hat{\sigma} : [0,T] \rightarrow L^2(\mathcal{S}^2,\R^3)$ be a bounded symmetric kernel and define
\begin{equation*}
    (\sigma_t f)(x) = \int_{\mathcal{S}} \hat{\sigma}(x,y,t) f(y) dy, \quad \forall f \in \mathcal{W}.
\end{equation*}
From this definition, $\sigma_t$ is a Hilbert-Schmidt operator mapping $\mathcal{W}$ into itself, hence the noise can be rewritten as
\begin{equation*}
    \sigma_t W_t = \sum_{i \in \mathbbm{N}} \hat{\beta}_t^i \sigma_t e_i,
\end{equation*}
this series converging in $\mathcal{W}$ almost surely, and in $L^p(\Omega, \mathcal{W})$ for all $p \in \mathbbm{N}$. Here $\mathcal{W}$ refers to the space carrying the Wiener process $W_t$, while the notation $L^2(\mathcal{S},\R^3)$ is kept for denoting the space of tridimensional velocities. Furthermore, there exists a Hilbert basis $(\phi_n)_n$ consisting of eigenfunctions of the noise operator $\sigma_t$ -- these eigenfunctions by their corresponding eigenvalues. Through a change of basis, there also exists a sequence of standard Brownian motions $(\beta_t^k)_k$, defined on the same filtered space, such that
\begin{equation*}
    \sigma_t W_t = \sum_{i \in \mathbbm{N}} \beta_t^k \phi_k.
\end{equation*}
We can further associate a covariance tensor to the random field $\sigma W_t$: given two spatial points $x,y \in \mathcal{S}$, and two time points $t,s \in \R^+$, we define $Q$ formally by 
\begin{equation*}
    Q(x,y,t,s) =  \E[(\sigma_t dW_t)(x) (\sigma_s dW_s)(y)] = \int_{\mathcal{S}} \hat{\sigma}(x,z,t) \hat{\sigma}(z,y,s) dz \delta(t-s).
\end{equation*}
The diagonal part of this covariance tensor is called the variance tensor, and denoted by
\begin{equation}
    a(x,t) = \int_{\mathcal{S}} \hat{\sigma}(x,y,t) \hat{\sigma}(y,x,t) dy = \sum_{k=0}^\infty \phi_k(x,t) \phi_k(x,t)^\tr \in \R^{3 \times 3}.
\end{equation}
The variance tensor $a$ is furthermore assumed to fulfil $a \in L^1([0,T], L^2(\mathcal{S}, \R^{ 3\times 3}))$, which can be guaranteed by the regularity assumptions on $\sigma$. Consequently, the integral $\int_0^t \sigma_s dW_s $ is a $\mathcal{W}$-valued Gaussian process with zero mean and bounded variance: $\E\big[\|\int_0^t \sigma_s dW_s\|_{L^2}^2\big]<\infty$. In addition, the quadratic variation of $\int_0^t (\sigma_t d W_s)(x)$ is given by the bounded-variation process $\int_0^t a(x,s) ds.$

As for classical Navier-Stokes equations, the LU Navier-Stokes equations can be derived using a stochastic version of the Reynolds Transport Theorem (SRTT) \cite{Mémin_2014}. Given a random scalar quantity $q$ transported by the flow within a volume $\mathcal{V}(t)$. Then, for incompressible unresolved flows -- that is $\nabla \bcdot \sigma_t=0$ -- the SRTT reads
\begin{gather}
    d\Big( \int_{\mathcal{V}(t)} q(x,t) dx \Big) = \int_{\mathcal{V}(t)} \big(\mathbbm{D}_t q + q \nabla \bcdot (u-u_s) dt\big) dx,\\
    \mathbbm{D}_t q = d_t q + (u - u_s) \bcdot \nabla q \: dt + \sigma dW_t \bcdot \nabla q - \frac{1}{2} \nabla \bcdot (a \nabla q) dt, \label{transport-operator}
\end{gather}
where an extra drift $u_s = \frac{1}{2} \nabla \bcdot a$ is involved, referred to as the It\={o}-Stokes drift in \cite{BCCLM_2020}. In this expression, $d_t q(x,t) = q(x,t+dt) - q(x,t)$ denotes the forward time increment at a fixed spatial point $x$, and $\mathbbm{D}_t q$ the stochastic transport operator as introduced in \cite{Mémin_2014, RMC_2017_Pt1}, which plays the role of the material derivative. The It\={o}-Stokes drift is directly tied to the divergence of the variance tensor $a$, and represents the effects of the noise inhomogeneity on the large-scale dynamics. Incorporating such advection terms is common practice in ocean large-scale modelling, since they can be regarded as practical correction terms \cite{CL_1976, GOP_2004, MSM_1997}. As shown in \cite{BCCLM_2020}, the LU framework naturally produces similar features, generalizing the effects of the small-scale inhomogeneity on the large-scale flow.

The stochastic transport operator contains physically interpretable terms describing the large-scale representation of flows. The last term on the right-hand side of \eqref{transport-operator} corresponds to an inhomogeneous diffusion term, accounting for small-scale mixing. This stochastic diffusion is fully determined by the noise variance and can be seen as a matrix generalization of the Boussinesq eddy-viscosity assumption. The third term on the right-hand side represents the transport of the large-scale quantity $q$ by the unresolved velocity. Remarkably, the energy carried by this backscattering term is exactly balanced by the stochastic diffusion term \cite{RMC_2017_Pt1}. This equilibrium can be interpreted as an instance of the fluctuation-dissipation theorem.

By interpreting the stochastic diffusion term $- \frac{1}{2} \nabla \bcdot (a \nabla q) dt$ as the It\={o}-Stratonovitch correction of the back-scattering It\={o} noise term $\sigma dW_t \bcdot \nabla q$ in the context of ``pure'' transports -- i.e. $\mathbbm{D}_t q = 0$ -- we may define a Stratonovitch transport operator
$$\mathbbm{D}^{\circ}_t q := d_t q + (u - u_s) \bcdot \nabla q \: dt + \sigma dW_t \circ \nabla q.$$
Here, $\sigma dW_t \circ \nabla q$ denotes a Stratonovitch noise term. Assuming this term is well-defined, the two transport operators coincide whenever the forcing terms have bounded variation, meaning that $(\mathbbm{D}_t q = dF)$ and $(\mathbbm{D}^{\circ}_t q = dF)$ are equivalent, whenever $dF = f dt$ for some process $f$.

In the rest of the article, our attention will focus on the Stratonovitch operator $\mathbbm{D}_t^\circ$, rather than the It\={o} one $\mathbbm{D}_t$. This modelling choice arises from the fact that the Stratonovitch product preserves the classical chain rule, while it is not the case for the It\={o} one. Moreover, one may observe the following: although the two transport operators are equivalent for a bounded variation forcing, the presence of a martingale component of force would dramatically affect the structure of the equation. As a matter of fact, we introduce below a martingale component of pressure, which precisely is a martingale forcing contribution. Hence, new terms representing the interaction between the transport noise and this martingale pressure would arise when rewriting the Stratonovitch form of the equation into the It\={o} one. Rewriting the problem accordingly, these terms can be interpreted as a transport of martingale pressure by the noise, and can be neglected under suitable scaling assumptions on the noise \cite{DMM2025}. When these assumptions are met, the two formulations (It\={o} and Stratonovitch) become equivalent.

\section{Derivation of the LU Water Waves equations}

\noindent
In this section, we derive a stochastic interpretation of the Water Waves equations. After recovering a Bernoulli-type equation on the internal potential in the stochastic setting, we rewrite it as an expression on the surface potential. The new relation involves the Dirichlet-to-Neumann operator, as is standard in the mathematical Water Waves literature \cite{lannes2013water}. Eventually, we propose an ``LU type of reasoning'' to show how our approach contains and goes beyond the perturbation of the deterministic Bernoulli equations.

To derive a location uncertainty interpretation of the Water Waves equations, we state first the LU Euler equations in Stratonovitch form, which are
\begin{gather}
    \mathbbm{D}_t^\circ u = d_t u + (u^* \bcdot \nabla) u dt + (\sigma dW_t \circ \nabla) u = - \frac{1}{\rho_0} \nabla (p dt + dp_t^\sigma) - g e_z dt, \label{eq-3d-velocity}\\
    \nabla \bcdot u = \nabla \bcdot u_s = \nabla \bcdot \sigma dW_t = 0, \label{eq-3d-div-free-cond}
\end{gather}
where $u=(v,w)^\tr$ denotes the tridimensional velocity, $v$ and $w$ being respectively the horizontal and vertical components. We also write $u^* = u-u_s = (v^*,w^*)^\tr$, where $u_s = \frac{1}{2} \nabla \bcdot a$, as defined in Section \ref{section:prelim-LU}. The notation $\sigma dW_t = (\sigma^H dW_t, \sigma^z dW_t)^\tr$ stands for the noise and is decomposed in horizontal and vertical contributions similarly. Moreover, $p$ and $p^\sigma$ represent the bounded variation and martingale pressure contributions, respectively. The total increment of (semi-martingale) pressure $dp^{tot} = p dt + dp_t^\sigma$ can in turn be interpreted as a Lagrange multiplier associated to the divergence-free condition $\nabla \bcdot u = 0$, where $u$ is itself a semi-martingale. 

\subsection{Boundary conditions on the top and the bottom}

We develop below the boundary conditions associated to the equations \eqref{eq-3d-velocity} and \eqref{eq-3d-div-free-cond}. Let $h_0>0$ represent the typical depth. At the bottom of the ocean $z=-h_0$, the non-penetration condition reads $w=0$. At the free surface $z= h(x,t)$, we denote by $\eta$ the surface elevation variation, so that,
\[
    z= h(x,t) = h_0\bigl(1 + \eta (x,t) \bigr),
\]
which leads to 
\begin{equation}
    \tilde z = \frac{1}{h_0}h(x,t) =\bigl(1+ \eta (x,t)\bigr).
\end{equation}
\begin{figure}
    \captionsetup[figure]{labelformat=empty}
    \centering
    \begin{tikzpicture}[>=latex, thick]

      \def\W{9}      
      \def\Hmid{3}   
    
      \def\ySurf{4.2}
      \draw[thick] (0,\ySurf-0.3)
        .. controls (0.8,\ySurf+0.8) and (1.5,\ySurf+0.9) .. (2.5,\ySurf+0.5)
        .. controls (3.2,\ySurf+0.2) and (3.8,\ySurf-0.5) .. (4.5,\ySurf-0.4)
        .. controls (5.2,\ySurf-0.3) and (6.0,\ySurf+0.6) .. (7.0,\ySurf+0.3)
        .. controls (7.8,\ySurf+0.1) and (8.4,\ySurf-0.2) .. (\W,\ySurf+0.1);

      \def\yBot{0}
      \draw[thick] (0, \yBot) -- (\W, \yBot);
    
      \node[right] at (\W+0.3, \ySurf + 0.3) {$z=0$};
    
      \draw[dashed] (0, \ySurf) -- (\W+0.3, \ySurf);
    
      \draw[->] (\W+1.1, 1) -- (\W+1.1, 2.5) node[right] {$z$};
    
      \node[above right] at (2.8, \ySurf+0.45) {Surface};
    
    
      \def\xDelta{1.35}
      \def\ySurfPic{4.9}  
      \draw[<->] (\xDelta, \ySurf) -- (\xDelta, \ySurfPic)
        node[midway, right] {$\eta$};
    
    
    
      \def\xH{4.5}
      \draw[<->] (\xH, \yBot) -- (\xH, \ySurf)
        node[midway, right] {$h_0$};
      \draw[<->] (\xH+1.5, \yBot) -- (\xH+1.5, \ySurf+0.3)
        node[midway, right] {$h$};
    

    \end{tikzpicture}
    \caption[Illustration of the water depth and surface deformation.]{Illustration of the water depth $h$ and surface deformation $\eta$. The bottom topography is chosen flat for simplicity.} \label{fig:summary}
\end{figure}
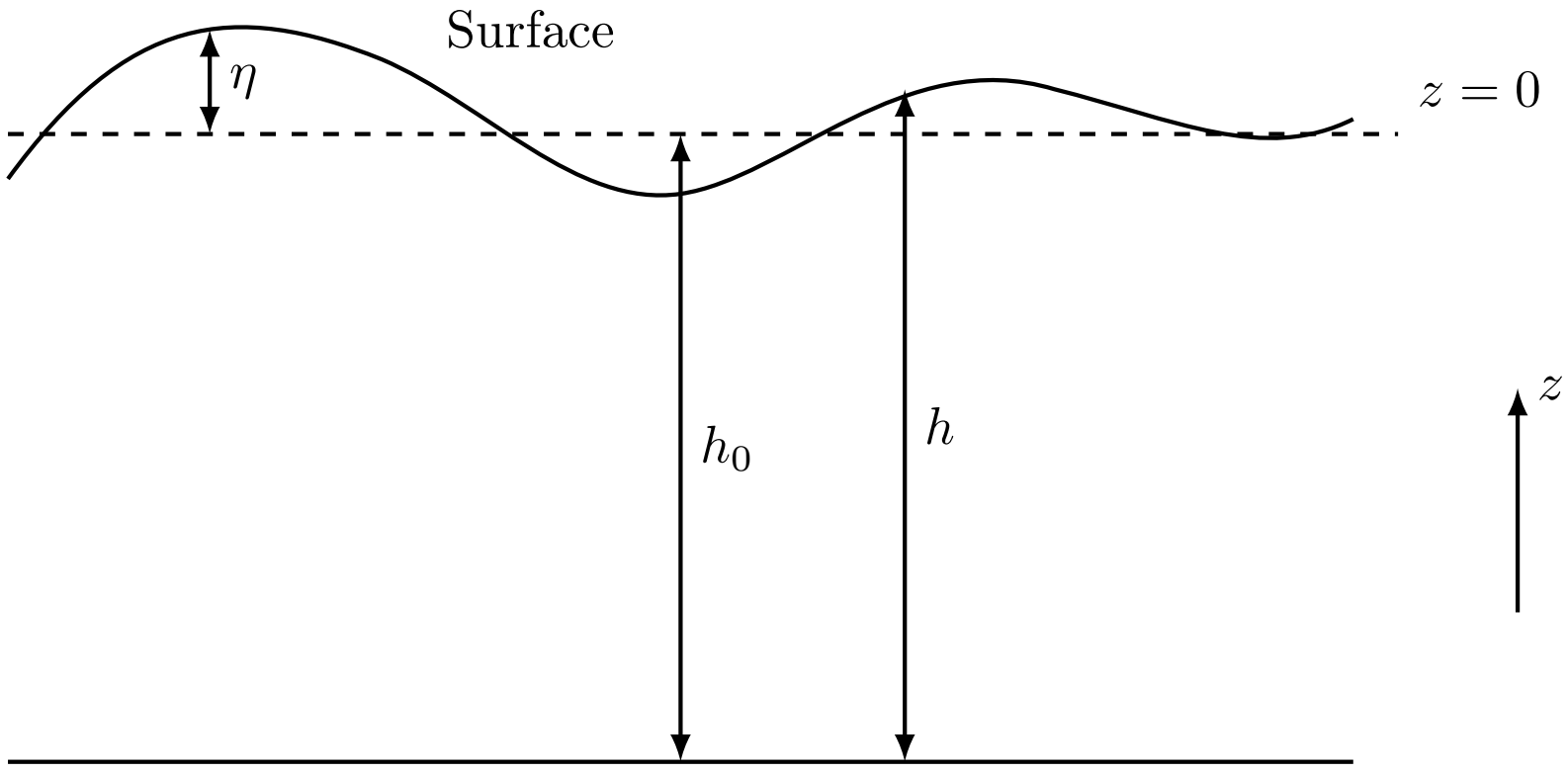

\noindent
Hence, the water surface is given by $\{\Sigma(X_t,t) = 0\}$, where $\Sigma(X_t,t) = \eta(X_t,t) - z$ and $X_t$ denotes the Lagragian displacement. Notice that $\eta(X_t,t)$ depends only on the horizontal coordinates of $X_t$. Next, the evolution of the surface elevation is given by
\begin{equation}
    \mathbbm{D}_t^\circ \Sigma(X_t,t) = 0. \label{eq:evolution-Sigma}
\end{equation}
Importantly, this modelling framework implicitly assumes that the system remains in the regime of non-breaking waves. The previous relation \eqref{eq:evolution-Sigma} can be rewritten as,
\begin{equation}
    \mathbbm{D}_t^\circ \eta = (w - w_s) dt  + \sigma^z dW_t. \label{eq:evolution-eta}
\end{equation}
Thus the boundary conditions at the top and the bottom of the ocean reads,
\begin{equation}
    \left\{
    \begin{array}{c}
         w\vert_{z=-h_0} = 0,\\
         \quad d_t \eta + (v^* dt + \sigma^H dW_t)\vert_{z=\eta} \circ \nabla_H \eta = (w^* dt + \sigma^z dW_t)\vert_{z=\eta}. 
    \end{array}
    \right.
    \label{eq-kinematic}
\end{equation}

\subsection{Deriving the stochastic Bernoulli equation}

In addition, as is standard in the (deterministic) derivation of the Water Waves system, we assume that the large-scale velocity $u$ is irrotational inside the fluid \cite{lannes2013water}, in order to infer a stochastic version of the Bernoulli equation. Notice that we make a weaker assumption for the noise below, as we will consider noises with non-zero vorticity. The large scale velocity thus reads as follows,
\begin{equation}
    u = \nabla \phi,
\end{equation}
where $\phi$ is the large-scale velocity potential. In particular, we infer that $\Delta \phi = 0$ from the divergence-free assumption. Moreover, by the Helmholtz decomposition of $\sigma dW_t$ and $u_s$, for a regular enough noise, we can write that
\begin{equation}
    \sigma dW_t = \nabla \sigma^\phi dW_t + \sigma^\omega dW_t, \quad u_s = \nabla \phi^S + \tilde{u}_s,
\end{equation}
where
\begin{equation}
    \sigma^\phi dW_t := \sum_k \phi_k^\sigma d\beta^k, \quad \sigma^\omega dW_t := \sigma dW_t - \sum_k \nabla \phi_k^\sigma d\beta_k, \quad \tilde{u}_s = u_s - \nabla \phi^S,
\end{equation}
where $\phi_k^\sigma,\phi^S$ denote potential scalar fields, $\tilde{u}_s$ and $\sigma^\omega dW_t$ being divergence-free. For modelling purposes, we assume that the vertical components of $\sigma^{\omega}dW_t$ and $u_s$ are negligible, which, from the divergence-free conditions on the noise and the It\={o}-Stokes drift, allows us to write the following,
\begin{equation}
    \sigma^{\omega} dW_t = \nabla_H^\perp \sigma^\chi dW_t = \sum_k \begin{pmatrix}
         \nabla_H^\perp \chi_k^{\sigma} \\ 0
    \end{pmatrix} d\beta^k, \quad \tilde{u}_s = \begin{pmatrix}
        \nabla_H^\perp \chi^{S} \\ 0
    \end{pmatrix},
\end{equation}
where $\chi_k^\sigma, \chi^S$ denote scalar fields, which can interpreted as $(x,y,z)$-dependent noise stream functions. In practise, this condition can be interpreted as the rotational components being purely horizontal, which makes sense in a context of non-breaking waves. 
Furthermore, we assume that the stream functions are harmonic in the following sense,
\begin{equation}
    \forall k>0, \Delta_3 \chi_k^{\sigma} = 0, \quad \Delta_3 \chi^{S} = 0.
\end{equation}
Consequently, from equation \eqref{eq-kinematic}, we find
\begin{equation}
    d_t \eta + [\nabla_H \phi \bcdot \nabla_H \eta dt + (\sigma^H dW_t - v_s dt)\vert_{z=\eta} \circ \nabla_H \eta - \partial_z \phi dt - (\sigma^z dW_t - w_s dt)]\vert_{z=\eta} = 0.
\end{equation}
Regarding the momentum equation \eqref{eq-3d-velocity}, we assume that
\begin{equation}
    p dt = dp_t^\sigma = 0, \text{ at } z = \eta. \label{eq:bc-pressure}
\end{equation}
Thus, equation \eqref{eq-3d-velocity} reads in terms of $\nabla \phi$ as follows,
\begin{multline}
    d_t \nabla \phi + \frac{1}{2} \nabla |\nabla \phi |^2 dt - (u_s \bcdot \nabla) \nabla \phi dt + (\sigma dW_t \circ \nabla) \nabla \phi + g \nabla z dt\\
    + \frac{1}{\rho_0} \nabla(p dt + dp_t^\sigma)= 0 \label{eq-before-proj}.
\end{multline}
Remind that the Stratonovitch product preserves the classical chain rule, which will be useful for the derivations below. At this point, we aim to infer an equation on the potential $\phi$. To achieve this, we must write the term
\begin{equation}
    (\sigma dW_t \circ \nabla) \nabla \phi - (u_s \bcdot \nabla) \nabla \phi dt
\end{equation}
as a gradient. We give hereafter the reasoning for the Stratonovitch noise component, as the It\={o}-Stokes one is treated similarly. Using Einstein's convention with $j \in \{x,y,z\}$, we infer that for all $i \in \{x,y,z\}$, 
\begin{align}
    [(\sigma dW_t \circ \nabla) \nabla \phi]_i = \sigma^j dW_t \circ \partial_j \partial_i \phi &= \partial_i (\sigma^j dW_t \circ \partial_j \phi) - (\partial_i \sigma^j dW_t) \circ (\partial_j \phi). \label{eq:bad-term-WW}
\end{align}
Then, the last term on the RHS of \eqref{eq:bad-term-WW} reads,
\begin{align}
    (\partial_i \sigma^j dW_t) \circ (\partial_j \phi) &=  (\partial_j \phi) \circ (\partial_{ij} \sigma^\phi dW_t + \partial_i \sigma^{\omega,j} dW_t) \nonumber\\
    &= (u_j \circ \partial_j) (\partial_{i} \sigma^\phi dW_t) +  (u_j \circ \partial_i) \sigma^{\omega,j} dW_t \nonumber\\
    &= (u_j \circ \partial_j) (\sigma^i dW_t) + (u_j \circ \partial_i \sigma^{\omega,j} dW_t - u_j \circ \partial_j \sigma^{\omega,i} dW_t) \nonumber\\
    &= [(u \circ \nabla) \sigma dW_t]_i + 2[\mathcal{A}( \nabla \sigma^{\omega} dW_t) \circ u]_i,
\end{align}
where, for any matrix $M$, $\mathcal{A} M = \frac{1}{2}(M-M^\tr)$ denotes its antisymmetric component. Now, we write that,
\begin{align*}
    2 [\mathcal{A}( \nabla \sigma^{\omega} dW_t) \circ u]_i &= \sum_{j \in \{x,y,z\}} \sum_k [\partial_{i^*} \partial_{j} \chi_k^\sigma - \partial_{j^*} \partial_{i} \chi_k^\sigma] \circ \partial_j \phi d\beta_t^k,\\
    &= \sum_{j \in \{x,y\}} \sum_k [\partial_{i^*} \partial_{j} \chi_k^\sigma - \partial_{j^*} \partial_{i} \chi_k^\sigma] \circ \partial_j \phi d\beta_t^k + \sum_k \partial_{i^*} \partial_z \chi_k^\sigma \circ \partial_z \phi d\beta_t^k\\
    &= \sum_{j \in \{x,y\}} \sum_k [\partial_{i^*} \partial_{j} \chi_k^\sigma - \partial_{j^*} \partial_{i} \chi_k^\sigma] \circ \partial_j \phi d\beta_t^k + \sum_k \partial_{i^*} \partial_z \chi_k^\sigma \circ \partial_z \phi d\beta_t^k.
\end{align*}
where we define formally $\partial_{l^*}$ by,
\begin{align*}
    \partial_{l^*} = \begin{cases}
        -\partial_y &\text{if $l=x$,}\\
        \partial_x &\text{if $l=y$},\\
        0 &\text{if $l=z$.}
    \end{cases}
\end{align*}
Hence, we infer that,
\begin{align*}
    2 [\mathcal{A}( \nabla \sigma^{\omega} dW_t) \circ u]_{i=x} &= - \sum_k \Delta_H \chi_k^\sigma \circ \partial_y \phi d\beta_t^k - \sum_k \partial_{y} \partial_z \chi_k^\sigma \circ \partial_z \phi d\beta_t^k,\\
    2 [\mathcal{A}( \nabla \sigma^{\omega} dW_t) \circ u]_{i=y} &= \sum_k \Delta_H \chi_k^\sigma \circ \partial_x \phi d\beta_t^k + \sum_k \partial_{x} \partial_z \chi_k^\sigma \circ \partial_z \phi d\beta_t^k,\\
    2 [\mathcal{A}( \nabla \sigma^{\omega} dW_t) \circ u]_{i=z} &= - \sum_{j \in \{x,y\}} \sum_k  \partial_{j^*} \partial_{z} \chi_k^\sigma \circ \partial_j \phi d\beta_t^k,
\end{align*}
which also reads more compactly as,
\begin{equation}
    2\mathcal{A}( \nabla \sigma^{\omega} dW_t) \circ u = \sum_k \Delta_H \chi_k^\sigma d\beta_t^k \circ \begin{pmatrix}
        \nabla_H^\perp \phi \\ 0
    \end{pmatrix} + \begin{pmatrix}
        \partial_z \phi \circ \partial_z \nabla_H^\perp \chi_{k}^\sigma d\beta_t^k \\
        - \nabla_H \phi \circ \partial_z \nabla_H^\perp \chi_{k}^\sigma d\beta_t^k
    \end{pmatrix}.
\end{equation}
Therefore, taking the divergence,
\begin{align}
    2\nabla \bcdot [\mathcal{A}( \nabla \sigma^{\omega} dW_t) \circ u] &= \nabla_H \Delta_H \chi_k^\sigma d\beta^k \circ \nabla_H^\perp \phi - \nabla_H \phi \circ \partial_{zz} \nabla_H^\perp \chi_{k}^\sigma d\beta_t^k \nonumber\\
    &= - \nabla_H \phi \circ \Big(\nabla_H^\perp  \Delta_H \chi_k^\sigma + \nabla_H^\perp \partial_{zz} \chi_{k}^\sigma\Big) d\beta_t^k = 0,
\end{align}
using the condition $\Delta_3 \chi_k^{\sigma} = 0$, and the fact that $\nabla_H^\perp f \bcdot \nabla_H g = - \nabla_H^\perp g \bcdot \nabla_H f$ for any differentiable functions $f,g$. Hence, we deduce that $\mathcal{A}(\nabla \chi_k^\sigma d\beta_k) \circ u$ is divergence-free under these assumptions. Consequently, the term $(\partial_i \sigma^j dW_t) \circ (\partial_j \phi)$ consists in two contributions: one that is divergence-free and one that is neglected in the LU framework, through an hypothesis of ``inertial'' noise balance. This corresponds to the Kraichnan-Tennekes hypothesis of turbulence transported by the large-scale random flow, which ensures that the acceleration of the smaller-scale dynamics is balanced by pressure terms associated to the same scale. \cite{Mémin_2014} Consequently, we make the following approximation,
\begin{align}
    (\partial_i \sigma^j dW_t) \circ (\partial_j \phi) &= [(u \circ \nabla) \sigma dW_t]_i + 2[\mathcal{A}( \nabla \sigma^{\omega} dW_t) \circ u]_i \\
    &\approx 2[\mathcal{A}( \nabla \sigma^{\omega} dW_t) \circ u]_i,\nonumber
\end{align}
which, in particular, implies that $(\partial_i \sigma^j dW_t) \circ (\partial_j \phi)$ is (quasi) divergence-free. 

\noindent
Now, denote by $\mathcal{P}$ the Leray projector onto the space of divergence-free functions fulfilling $w = 0$ at $z=-h_0$ and $u \bcdot n = 0$ at $z=\eta$. Then, we denote by $\mathcal{Q} := \mathcal{I} - \mathcal{P}$, so that, by applying $\mathcal{Q}$ to equation \eqref{eq-before-proj}, and using the aforementioned assumption, the final system reads,
\begin{gather}
    d_t \eta + \Big[\nabla_H \phi \bcdot \nabla_H \eta dt + (\sigma^H dW_t - v_s dt) \circ \nabla_H \eta - \partial_z \phi dt - (\sigma^z dW_t - w_s dt)\Big]\Big\vert_{z=\eta} = 0, \\
    d_t \phi + \frac{1}{2} |\nabla \phi |^2 dt + g z dt + (\sigma dW_t - u_s dt) \circ \nabla \phi = - \frac{1}{\rho_0}(pdt + dp_t^\sigma). \label{eq:stochastic-Bernoulli}
\end{gather}
This system provides an intuitive extension of the deterministic Bernoulli equation, that is coupled with a free-surface kinematic boundary condition. As previously mentioned, the noise can be of non-zero vorticity through the parameters $\chi_k^\sigma$.

\subsection{Scaling the equations}

\noindent
In order to scale the stochastic Water Waves equations we aim to derive, we introduce the following parameters. First, we define $H$ and $L$ as the typical vertical and horizontal lengths of the system, respectively. In addition, $A$ denotes the typical height of the wave i.e. the order of magnitude of $\eta$, and $L_\sigma$ represents the typical length associated to $\sigma dW_t$. Then, we define several non-dimensional parameters, as is standard in the Water Waves literature \cite{lannes2013water},
\begin{itemize}
    \item[$\bullet$] $\epsilon := A/H$ the nonlinear parameter,
    \item[$\bullet$] $\mu := (H/L)^2$ the squared aspect ratio,
    \item[$\bullet$] $\nu := \frac{\tanh(2\pi \sqrt{\mu})}{2\pi \sqrt{\mu}}$, the scaling parameter. This allows to distinguish the deep-water regimes ($\nu \sim 2 \pi \sqrt{\mu}$ from the shallow-water ones ($\nu \sim 1$).
\end{itemize}
Notice that $\mu$ corresponds to $\beta^2$ as defined in \cite{DMM_2024}. We adopt this convention to be coherent with the Water Wave literature. We also define stochastic a few stochastic parameters as follows,
\begin{itemize}
    \item[$\bullet$] $r_\sigma := L_\sigma/L$, the ratio between the noise length and the resolved length,
    \item[$\bullet$] $\epsilon_\sigma := \nu r_\sigma$ the parameter of the non-linearity modelled by the noise,
    \item[$\bullet$] $\epsilon_S := \nu r_\sigma^2 = \epsilon_\sigma^2/\nu$, which can be interpreted as a non-linearity parameter associated to the It\={o}-Stokes drift.
\end{itemize}
In the definitions above, the equality $L_\sigma/L = \epsilon_\sigma/\nu$ is set \emph{a priori} for convenience, and thus can be thought of as a link between the typical length scale of the noise $L_\sigma$ and the typical water height variation due to the noise. As a matter of fact, if we this typical variation as $A_\sigma = \epsilon_\sigma H$, we infer from $\frac{L_\sigma}{L} = \frac{\epsilon_\sigma}{\nu}$,
\begin{equation}
    \frac{L_\sigma}{L} = \frac{A_\sigma}{H \nu} \Rightarrow A_\sigma = \sqrt{\mu} \nu L_\sigma = \frac{1}{2\pi} \tanh(2\pi \sqrt{\mu}) L_\sigma. 
\end{equation}
Implicitly, we assume that the bounded variation and martingale contributions follow similar regimes terms of the parameters $\mu$ and $\nu$. Moreover, we introduce the scaling parameters $t_0$ on the time variable and $\phi_0$ on the resolved potential field,
\begin{equation}
    t_0 := \frac{L}{\sqrt{gH\nu}}\quad \phi_0 := \frac{\epsilon L^2}{\nu t_0} = \frac{A L}{H} \sqrt{\frac{gH}{\nu}}.
\end{equation}
They are defined as in \cite{lannes2013water}, where they are derived from linear solutions analysis. However, as our noise is already scaled by a ratio of lengths $r_\sigma$, the $t_0$ dependence would be redundant. For this reason, only the scaling in terms of $\phi_0$ appears in the equations below, through the factor $\frac{\epsilon}{\nu}$ associated to the quadratic term $|\nabla \phi|^2$ in the Bernoulli equation. 

\noindent
Using the scaling parameters above, the conditions on $\phi$, $\phi_k^\sigma$, $\phi^S$, $\chi_k^{\sigma}$ and $\chi^{S}$ become
\begin{gather}
    (\mu \Delta_H + \partial_{zz}) \phi = (\mu \Delta_H + \partial_{zz}) \phi_k^\sigma= (\mu \Delta_H + \partial_{zz}) \phi^S = 0,\\
    (\mu \Delta_H + \partial_{zz})\chi_k^{\sigma} = (\mu \Delta_H + \partial_{zz}) \chi^{S} = 0.
\end{gather}

\noindent
Then, denoting by $\nabla^\mu = (\sqrt{\mu} \nabla_H, \partial_z)^\tr$, we write the non-dimensional system,
\begin{align}
    d_t \eta &+ \frac{1}{\mu \nu} \Big[ \epsilon \mu \nabla_H (\phi - \frac{\epsilon_S}{\epsilon} \phi^S) \bcdot \nabla_H \eta - \partial_z (\phi - \frac{\epsilon_S}{\epsilon} \phi^S) \Big]\Big\vert_{z=\eta} dt \\
    + \frac{\epsilon_\sigma/\epsilon}{\mu \nu} &\Big[ \epsilon \mu \nabla_H \sigma^{\phi} dW_t \circ \nabla_H \eta - \partial_z \sigma^{\phi} dW_t \Big]\Big\vert_{z=\eta} 
    + \Big( \frac{\epsilon_\sigma}{\nu} \nabla_H^\perp \sigma^{\chi} dW_t - \frac{\epsilon_\sigma^2}{\nu^2} \nabla_H^\perp \chi^{S} dt \Big)\Big\vert_{z=\eta} \circ \nabla_H \eta = 0, \nonumber\\
    d_t \phi &+ \frac{\epsilon}{2\mu\nu} |\nabla^\mu \phi|^2 dt + z dt + \frac{r_\sigma}{\mu} (\nabla^\mu \sigma^\phi dW_t + \sqrt{\mu} \nabla_H^\perp \sigma^{\chi} dW_t) \circ \nabla^\mu \phi\\
    &\qquad \qquad \qquad \qquad \qquad - \frac{r_\sigma^2}{\mu} (\nabla^\mu \phi^S + \sqrt{\mu} \nabla_H^\perp \chi^S) dt \bcdot \nabla^\mu \phi = - \frac{C}{\rho_0}(pdt + dp_t^\sigma),\nonumber
\end{align}
where we denote by $C$ a constant depending on the parameters. As we shall see below, the pressure terms will be cancelled out due to the boundary condition at the surface \eqref{eq:bc-pressure}. Moreover, we write $\nabla_H^\perp \sigma^{\chi} dW_t$ and $\nabla_H^\perp \chi^{S}$ in place of
$\begin{pmatrix}
        \nabla_H^\perp \sigma^{\chi} dW_t \\ 0
\end{pmatrix}$ and $\begin{pmatrix}
        \nabla_H^\perp \chi^S \\ 0
\end{pmatrix}$, respectively, by abuse of notation. For convenience, we further denote by $h = 1 + \epsilon \eta$. Notice that we have rewritten the problem in terms of $\sigma^\phi dW_t$, $\phi^S dt$, $\sigma^\chi dW_t$ and $\chi^S dt$. With such expressions, we are in position to introduce the Dirichlet-to-Neumann operator, as is standard in the water waves literature. 

\subsection{The Dirichlet-to-Neumann operator}

\noindent
In this section, we introduce the Dirichlet-to-Neumann operator (see \cite{lannes2013water}). Let $\Phi$ denote the trace of the large-scale velocity potential $\phi$ at the surface, that is $\Phi(x,t) = \phi(x,\eta(x,t),t)$. Thus, $\phi$ is a solution of the elliptic problem
\begin{equation}
    (\mu \Delta_H + \partial_{zz}) \phi = 0, \quad \partial_z \phi = 0 \text{ at $z=-h_0$}, \quad \phi = \Phi \text{ at $z=\eta$.}
\end{equation}
In this formalism, we introduce the Dirichlet-to-Neumann operator
\begin{equation}
    G(\eta)\Phi := (\partial_z \phi)\vert_{z=\eta} - \epsilon \mu (\nabla_H \eta) \bcdot (\nabla_H \phi)\vert_{z=\eta}.
\end{equation}
This operator aims to transform the Dirichlet boundary condition on the velocity to a Neumann one on the potential, hence its name. As a result, the free-surface evolution reads,
\begin{equation}
    d_t \eta = \frac{1}{\mu \nu} G(\eta)(\Phi) dt + dF^\gamma(\eta) + dF^\omega(\eta).
\end{equation}
where
\begin{gather}
    dF^\gamma(\eta) = \frac{\epsilon_\sigma/\epsilon}{\mu \nu} \sum_k G(\eta)(\Phi_k^\sigma) \circ d\beta^k - \frac{\epsilon_S/\epsilon}{\mu \nu} G(\eta)(\Phi^S) dt\\
    dF^\omega(\eta) = (- r_\sigma \nabla_H^\perp \sigma^{\chi} dW_t + r_\sigma^2 \nabla_H^\perp \chi^{S} dt)\vert_{z=\eta} \circ \nabla_H \eta.
\end{gather}
Several identities relating the partial derivatives to the Dirichlet-to-Neumann operator will be necessary in the following. Using the (deterministic) chain rule, we find that
\begin{equation}
    \nabla_H \Phi = (\nabla_H \phi)\vert_{z=\eta} + \epsilon (\partial_z \phi)\vert_{z=\eta} \nabla_H \eta.
\end{equation}
Taking the inner product of the equation above with $\epsilon \mu \nabla_H \eta$ yields,
\begin{equation}
    \epsilon \mu (\nabla_H \eta \bcdot \nabla_H \Phi) = \epsilon \mu \Big( \nabla_H \eta \bcdot (\nabla_H \phi)\vert_{z=\eta}\Big) + \epsilon^2 \mu (\partial_z \phi)\vert_{z=\eta} |\nabla_H \eta|^2.
\end{equation}
Consequently,
\begin{equation}
    G(\eta)\Phi = (\partial_z \phi)\vert_{z=\eta}  (1 + \epsilon^2 \mu |\nabla_H \eta|^2) - \epsilon \mu (\nabla_H \eta \bcdot \nabla_H \Phi),
\end{equation}
that is
\begin{equation}
    (\partial_z \phi)\vert_{z=\eta} = \frac{G(\eta)\Phi + \epsilon \mu (\nabla_H \eta \bcdot \nabla_H \Phi)}{1 + \epsilon^2 \mu |\nabla_H \eta|^2} =: \mu \mathcal{Z}(\eta)\Phi.
\end{equation}
In addition, we infer
\begin{equation}
    (\nabla_H \phi)\vert_{z=\eta} = \nabla_H \Phi -  \epsilon \mu (\nabla_H \eta) \mathcal{Z}(\eta)\Phi.
\end{equation}
Reasoning similarly with the potential quantities $\phi_k^\sigma, \phi^S$ and their traces at the surface $\Phi_k^\sigma, \Phi^S$, we obtain
\begin{gather}
    (\nabla \phi)\vert_{z=\eta} = \begin{pmatrix}
        \nabla_H \Phi - \epsilon \mu \nabla_H \eta \mathcal{Z}(\eta)\Phi \\
        \mathcal{Z}(\eta)\Phi
    \end{pmatrix}\\
    (\nabla \phi_k^\sigma)\vert_{z=\eta} = \begin{pmatrix}
        \nabla_H \Phi_k^\sigma -  \epsilon \mu \nabla_H \eta \mathcal{Z}(\eta)\Phi_k^\sigma  \\
        \mathcal{Z}(\eta)\Phi_k^\sigma
    \end{pmatrix}\\
    (\nabla \phi^S)\vert_{z=\eta} = \begin{pmatrix}
        \nabla_H \Phi^S - \epsilon \mu \nabla_H \eta \mathcal{Z}(\eta)\Phi^S \\
        \mathcal{Z}(\eta)\Phi^S
    \end{pmatrix}.
\end{gather}
Applying It\={o}'s lemma to the potential at the surface $\Phi$ yields,
\begin{equation}
    d_t \Phi = (d_t \phi)\vert_{z=\eta} + \epsilon (\partial_z \phi)\vert_{z=\eta} \circ d_t \eta.
\end{equation}
Hence, we find
\begin{multline}
    d_t \Phi + \frac{\epsilon}{2\mu\nu} \Bigg[\mu\Big[\nabla_H \Phi -  \epsilon \mu \nabla_H \eta \mathcal{Z}(\eta)\Phi\Big]^2 + \Big[ \mu \mathcal{Z}(\eta)\Phi\Big]^2 \Bigg] dt\\
    + \frac{r_\sigma}{\mu} \sum_k \begin{pmatrix}
        \sqrt{\mu} \nabla_H \Phi_k^\sigma -  \epsilon \mu \sqrt{\mu} \nabla_H \eta \mathcal{Z}(\eta)\Phi_k^\sigma \\
        \mu \mathcal{Z}(\eta)\Phi_k^\sigma
    \end{pmatrix} d\beta^k \circ \begin{pmatrix}
        \sqrt{\mu} \nabla_H \Phi -  \epsilon \mu \sqrt{\mu} \nabla_H \eta \mathcal{Z}(\eta)\Phi \\
        \mu \mathcal{Z}(\eta)\Phi
    \end{pmatrix}\\
    - \frac{r_\sigma^2}{\mu} \begin{pmatrix}
        \sqrt{\mu} \nabla_H \Phi^S -  \epsilon \mu \sqrt{\mu} \nabla_H \eta \mathcal{Z}(\eta)\Phi^S \\
        \mu \mathcal{Z}(\eta)\Phi^S
    \end{pmatrix} \bcdot \begin{pmatrix}
        \sqrt{\mu} \nabla_H \Phi -  \epsilon \mu \sqrt{\mu} \nabla_H \eta \mathcal{Z}(\eta)\Phi \\
        \mu \mathcal{Z}(\eta)\Phi
    \end{pmatrix} dt\\
    +\Big(r_\sigma \nabla_H^\perp \sigma^{\chi} dW_t - r_\sigma^2 \nabla_H^\perp \chi^S dt\Big)\Big\vert_{z=\eta} \circ \Big(
    \nabla_H \Phi -  \epsilon \mu \nabla_H \eta \mathcal{Z}(\eta)\Phi\Big)\\
    - \epsilon \Big[\frac{1}{\mu \nu} G(\eta)\Phi dt + dF^\gamma(\eta) + dF^\omega(\eta)\Big] \circ \mathcal{Z}(\eta)\Phi + \eta = 0,
\end{multline}
where the pressure terms cancel due to the boundary condition at the surface \eqref{eq:bc-pressure}. By a classical computation, the last equation can be rewritten as
\begin{equation}
    d_t \Phi + \frac{\epsilon}{2 \nu} |\nabla \Phi|^2 + \eta - \frac{\epsilon \mu}{\nu} \frac{[\frac{1}{\mu} G(\eta)\Phi + \epsilon \nabla_H \Phi \bcdot \nabla_H \eta]^2}{2( 1 + \epsilon^2 \mu |\nabla \eta|^2)} + dF^\sigma(\eta,\Phi) = 0, 
\end{equation}
where we denote by
\begin{multline}
    dF^\sigma(\eta, \Phi) := \frac{r_\sigma}{\mu} \sum_k \begin{pmatrix}
        \sqrt{\mu} \nabla_H \Phi_k^\sigma -  \epsilon \mu \sqrt{\mu} \nabla_H \eta \mathcal{Z}(\eta)\Phi_k^\sigma \\
        \mu \mathcal{Z}(\eta)\Phi_k^\sigma
    \end{pmatrix} d\beta^k \circ \begin{pmatrix}
        \sqrt{\mu} \nabla_H \Phi -  \epsilon \mu \sqrt{\mu} \nabla_H \eta \mathcal{Z}(\eta)\Phi \\
        \mu \mathcal{Z}(\eta)\Phi
    \end{pmatrix}\\
    - \frac{r_\sigma^2}{\mu} \begin{pmatrix}
        \sqrt{\mu} \nabla_H \Phi^S -  \epsilon \mu \sqrt{\mu} \nabla_H \eta \mathcal{Z}(\eta)\Phi^S \\
        \mu \mathcal{Z}(\eta)\Phi^S
    \end{pmatrix} \bcdot \begin{pmatrix}
        \sqrt{\mu} \nabla_H \Phi -  \epsilon \mu \sqrt{\mu} \nabla_H \eta \mathcal{Z}(\eta)\Phi \\
        \mu \mathcal{Z}(\eta)\Phi
    \end{pmatrix} dt\\
    +\Big(r_\sigma \nabla_H^\perp \sigma^{\chi} dW_t - r_\sigma^2 \nabla_H^\perp \chi^S dt\Big)\Big\vert_{z=\eta} \circ \Big(
    \nabla_H \Phi - \epsilon \mu \nabla_H \eta \mathcal{Z}(\eta)\Phi\Big)\\
    - \epsilon \Big[dF^\gamma(\eta) + dF^\omega(\eta)\Big] \circ \mathcal{Z}(\eta)\Phi.
\end{multline}
The term $dF^\sigma(\eta, \Phi)$ can be rewritten as follows,
\begin{multline}
    dF^\sigma(\eta, \Phi) :=  r_\sigma \sum_k \Big(\nabla_H \Phi_k^\sigma - \epsilon \mu \nabla_H \eta \mathcal{Z}(\eta)\Phi_k^\sigma \Big) d\beta^k \circ \Big(\nabla_H \Phi - \epsilon \mu \nabla_H \eta \mathcal{Z}(\eta)\Phi\Big) \\
    + \sum_k \frac{\epsilon_\sigma}{\nu} \epsilon \mu (\nabla_H \eta \bcdot \nabla_H \Phi_k^\sigma) d\beta^k \circ \mathcal{Z}(\eta)\Phi\\ 
    - r_\sigma^2 \Big(\nabla_H \Phi^S - \epsilon \mu \nabla_H \eta \mathcal{Z}(\eta)\Phi^S\Big) \bcdot \Big(\nabla_H \Phi -  \epsilon \mu \nabla_H \eta \mathcal{Z}(\eta)\Phi\Big) dt \\
    - \frac{\epsilon_\sigma^2}{\nu^2} \epsilon \mu (\nabla_H \eta \bcdot \nabla_H \Phi^S) \mathcal{Z}(\eta)\Phi + \Big(r_\sigma \nabla_H^\perp \sigma^{\chi} dW_t - r_\sigma^2 \nabla_H^\perp \chi^S dt\Big)\Big\vert_{z=\eta} \circ \nabla_H \Phi.
\end{multline}
This also reads, equivalently,
\begin{multline}
    dF^\sigma(\eta, \Phi) := (r_\sigma \nabla_H \sigma^\Phi dW_t -  r_\sigma^2 \nabla_H \Phi^S dt) \circ \nabla_H \Phi + \Big(r_\sigma \nabla_H^\perp \sigma^{\chi} dW_t - r_\sigma^2 \nabla_H^\perp \chi^S dt\Big)\Big\vert_{z=\eta} \circ \nabla_H \Phi\\
    - \epsilon \mu (\nabla_H \eta \bcdot \nabla_H \Phi) \circ \mathcal{Z}(\eta)\Big(r_\sigma \sigma^\Phi dW_t - r_\sigma^2 \Phi^S dt\Big)\\
    + \epsilon^2 \mu^2 |\nabla_H \eta|^2 \mathcal{Z}(\eta)\Phi \circ \mathcal{Z}(\eta) \Big( r_\sigma \sigma^\Phi dW_t - r_\sigma^2 \Phi^S dt\Big).
\end{multline}
Thus, we infer the following expression for the LU Water Waves system, using that $r_\sigma = \epsilon_\sigma/\nu$,
\begin{gather}
    d_t \eta = \frac{1}{\mu \nu} G(\eta)\Phi + dF^\gamma(\eta) + dF^\omega(\eta), \label{eq-eta-final}\\
    d_t \Phi + \frac{\epsilon}{2 \nu} |\nabla \Phi|^2 + \eta - \frac{\epsilon \mu}{\nu} \frac{[\frac{1}{\mu} G(\eta)\Phi + \epsilon \nabla_H \Phi \bcdot \nabla_H \eta]^2}{2( 1 + \epsilon^2 \mu |\nabla \eta|^2)} \nonumber\\
    + \Big( \nabla_H (\frac{\epsilon_\sigma}{\nu} \sigma^\Phi dW_t - \frac{\epsilon_\sigma^2}{\nu^2} \Phi^S dt) + (\frac{\epsilon_\sigma}{\nu} \nabla_H^\perp \sigma^{\chi} dW_t - \frac{\epsilon_\sigma^2}{\nu^2} \nabla_H^\perp \chi^S dt)\vert_{z=\eta}\Big) \circ \nabla_H \Phi \nonumber \\
    - \epsilon \mu \Big(\nabla_H \eta \bcdot \nabla_H \Phi - \epsilon \mu |\nabla_H \eta|^2 \mathcal{Z}(\eta)\Phi\Big) \circ \mathcal{Z}(\eta)\Big( \frac{\epsilon_\sigma}{\nu} \sigma^\Phi dW_t - \frac{\epsilon_\sigma^2}{\nu^2} \Phi^S dt\Big) = 0, \label{eq-phi-final}
\end{gather}
reminding that
\begin{equation}
    \mathcal{Z}(\eta)(\bcdot) := \frac{1}{\mu} \frac{G(\eta)(\bcdot) + \epsilon \mu \nabla_H \eta \bcdot \nabla_H (\bcdot)}{1 + \epsilon^2 \mu |\nabla_H \eta|^2}.
\end{equation}
\begin{remark}
    At this point we can make a useful sanity check: assume that the noise is purely potential, that is $\chi_k^\sigma = 0$ for all $k \geq 1$, and that the It\={o}-Stokes drift $u_s = 0$. Thus, if we neglect the terms of order $O(\epsilon^2\mu^2)$, the system boils down to
    \begin{gather}
        d_t \eta = \frac{1}{\mu \nu} G(\eta)\Phi + dF^\gamma(\eta),\\
        d_t \Phi + \frac{\epsilon}{2 \nu} |\nabla \Phi|^2 + \eta - \frac{\epsilon \mu}{\nu} \frac{[\frac{1}{\mu} G(\eta)\Phi + \epsilon \nabla_H \Phi \bcdot \nabla_H \eta]^2}{2( 1 + \epsilon^2 \mu |\nabla \eta|^2)} + \frac{\epsilon_\sigma}{2\nu} \sum_k \nabla_H \Phi_k^\sigma d\beta^k \circ \nabla_H \Phi \nonumber\\
        - \epsilon \mu (\nabla_H \eta \bcdot \nabla_H \Phi) \circ \mathcal{Z}(\eta)\Big(\frac{\epsilon_\sigma}{\nu} \sigma^\Phi dW_t\Big) = O(\epsilon^2 \mu^2).
    \end{gather}
    This formulation is in turn equivalent to the system that is presented in \cite{Street_2023}. In particular, it shows that our model recovers known behaviours under suitable assumptions. Yet, the term of order $\epsilon^2 \mu^2$ is practically neglected by other arguments in the aforementioned article, as the vertical velocity $(\sigma^z dW_t)\vert_{z=\eta}$ cancels out in the surface potential equation. Importantly, our derivation is rooted in the LU formalism, while that of \cite{Street_2023} emerges from SALT. These two approaches yield in full generality the same wave equations whenever the noise is irrotational.
\end{remark}
A more rigorous and complete approach could be deriving a stochastic system in the fashion of \cite{CL2015}. There, the evolution of the rotational component is fully taken into account. However, the derivation would be much more challenging, for example since $\Phi$ cannot be defined the same way when accounting for rotational components of the flow. The interest of our approach is that it allows to derive a stochastic model with possibly rotational components, while remaining in the framework of the classical Water Waves theory.

\subsection{Another approach: ``direct'' stochastic Bernoulli equation}

\noindent
For potential flows and noises, it is possible to start from the Bernoulli equations and apply a reasoning that is similar to the derivation of LU. \cite{Mémin_2014} In the deterministic setting, the Bernoulli equations read, for an irrotational fluid, as follows,
\begin{gather}
    d_t \eta + (\nabla_H \phi)\vert_{z=\eta} \bcdot \nabla_H \eta - (\partial_z \phi) \vert_{z=\eta} = 0,\\
    d_t \phi + \frac{1}{2} |\nabla \phi |^2 + g z = - \frac{1}{\rho_0} (p dt + dp_t^\sigma),
\end{gather}
where, reusing the same notations, $\eta$ and $\phi$ denote the surface height and the velocity potential, respectively. Decomposing the latter as $\phi := \bar{\phi} + \phi'$, $\bar{\phi}$ and $\phi'$ being the time-correlated large-scale and time-decorrelated small-scale potentials, respectively, yields,
\begin{gather}
    d_t \eta + \nabla_H (\bar{\phi} + \phi')\vert_{z=\eta} \bcdot \nabla_H \eta - \partial_z (\bar{\phi} + \phi')\vert_{z=\eta} = 0,\\
    d_t (\bar{\phi} + \phi') + \frac{1}{2} |\nabla (\bar{\phi} + \phi') |^2 + g z = - \frac{1}{\rho_0} (p dt + dp_t^\sigma), \label{eq-potential-smallLargeScales}
\end{gather}
The only nonlinear term in the previous system is $\frac{1}{2} |\nabla (\bar{\phi} + \phi') |^2$, which we develop so that,
\begin{equation}
    \frac{1}{2} |\nabla (\bar{\phi} + \phi') |^2 = \frac{1}{2} |\nabla \bar{\phi}|^2 + \nabla \bar{\phi} \bcdot \nabla \phi' +  \frac{1}{2} |\nabla \phi' |^2.
\end{equation}
Hence, we may define two kinematic operators associated to potential, that is for all potential $f$,
\begin{equation}
    \mathbbm{D}_t^{WW} f = d_t f + \frac{1}{2} \nabla \bar{\phi} \bcdot \nabla f + \nabla \phi' \bcdot \nabla f, \quad D_t^{WW} f = d_t f + \frac{1}{2} |\nabla f|^2.
\end{equation}
Thus, equation \eqref{eq-potential-smallLargeScales} reads,
\begin{equation}
    \mathbbm{D}_t^{WW} \bar{\phi} + D_t^{WW} \phi' + g z = - \frac{1}{\rho_0} (p dt + dp_t^\sigma).
\end{equation}
Postulating a strong scale separation, we may neglect the term $D_t^{WW} \phi'$ in the previous equation, so that
\begin{equation}
    \mathbbm{D}_t^{WW} \bar{\phi} + g z = d_t \bar{\phi} + \frac{1}{2} | \nabla \bar{\phi} |^2 + \nabla \phi' \bcdot \nabla \bar{\phi} + g z = - \frac{1}{\rho_0} (p dt + dp_t^\sigma).
\end{equation}
In the following, we model $\phi' dt$ by $r_\sigma \sum_k \phi_k^\sigma d\beta^k$, and denote by $\sigma dW_t = \sum_k \nabla \phi_k^\sigma d\beta^k$. Thus, dropping the bar notation for readability, the system reads
\begin{gather*}
    d_t \eta + (\nabla_H \phi)\vert_{z=\eta} \circ \nabla_H \eta - (\partial_z \phi)\vert_{z=\eta} + r_\sigma \sum_k \Big[(\nabla_H \phi_k^\sigma)\vert_{z=\eta} \circ \nabla_H \eta - (\partial_z \phi_k^\sigma)\vert_{z=\eta} \Big] d\beta^k = 0,\\
    d_t \phi + \frac{1}{2} | \nabla \phi |^2 dt + r_\sigma \sum_k \nabla \phi_k^\sigma d\beta^k \circ \nabla \phi + g z dt = - \frac{1}{\rho_0} (p dt + dp_t^\sigma).
\end{gather*}
Hence, this system is contained in equation \eqref{eq:stochastic-Bernoulli}, and therefore in the previously derived LU Water-Waves model. Effectively, it is recovered by the assumptions of potential noise and zero It\={o}-Stokes drift -- that is $\sigma dW_t = \sum_k \nabla \phi_k^\sigma d\beta^k$ and $u_s = 0$.

\section{Derivation of stochastic interpretations of shallow-water wave models}

\noindent
Considering the LU Water Waves system in the form of equations \eqref{eq-eta-final} and \eqref{eq-phi-final}, we now express a few classical waves models in the LU framework as formal limits of the Water Waves one. The models studied are Saint-Venant, Serre-Green-Nagdhi and Boussinesq.

\subsection{The Saint-Venant system}

\noindent
The (deterministic) Saint-Venant system corresponds to a zeroth order approximation of the Water Waves system, in terms of the squared aspect ratio. As shown in \cite{lannes2013water}, the following identity holds,
\begin{equation}
    \frac{1}{\mu} G(\eta)\Phi = - \nabla_H \bcdot (h \bar{u}). \label{eq:relation-Dir-to-Neum}
\end{equation}
Moreover, in the same reference, it is proven that
\begin{equation}
    \bar{u} := \frac{1}{h} \int_{-1}^{\epsilon \eta} \nabla_H \phi \: dz = \nabla_H \Phi + O(\mu) = \nabla_H \phi + O(\mu). \label{eq:approx-uBar}
\end{equation}
Equation \eqref{eq:relation-Dir-to-Neum} arises directly from Leibniz formula, in both deterministic and stochastic settings. Moreover, since equation \eqref{eq:approx-uBar} holds in the deterministic regime, we assume that it also holds in the stochastic one. In turn, we find the following system,
\begin{align}
    d_t \eta &+ \frac{1}{\nu} \nabla_H \bcdot (h \bar{u}) dt \\
    &+ \frac{1}{\nu}\nabla_H \bcdot \Big(h \circ (\frac{\epsilon_\sigma}{\epsilon} \overline{\nabla_H \sigma^{\phi} dW_t} + \frac{\epsilon_\sigma}{\nu} \overline{\nabla_H^\perp \sigma^{\chi} dW_t} - \frac{\epsilon_\sigma^2}{\nu \epsilon} \overline{\nabla_H \phi^S} dt - \frac{\epsilon_\sigma^2}{\nu^2} \overline{\nabla_H ^\perp \chi^S} dt) \Big) = 0, \nonumber\\
    d_t \bar{u} &+ \frac{\epsilon}{2 \nu} \nabla_H |\bar{u}|^2 dt + \nabla_H \eta dt \\
    &+ \nabla_H \Big(( \nabla_H (\frac{\epsilon_\sigma}{\nu}\sigma^\Phi dW_t - \frac{\epsilon_\sigma^2}{\nu^2} \Phi^S dt) + (\frac{\epsilon_\sigma}{\nu} \nabla_H^\perp \sigma^{\chi} dW_t - \frac{\epsilon_\sigma^2}{\nu^2} \nabla_H^\perp \chi^S dt)\vert_{z=\eta}) \circ \bar{u}\Big) = O(\frac{\mu}{\nu}, \mu \epsilon_\sigma), \nonumber
\end{align}
where we use $h = 1 + \epsilon \eta$. Notice that this system is equivalent to the Saint-Venant system derived in \cite{DMM_2024} for divergence-free noises, except for the following term $T_{SV}$,
\begin{equation}
    T_{SV} := (\bar{u} \circ \nabla_H) \nabla_H (\frac{\epsilon_\sigma}{\nu}\sigma^\Phi dW_t - \frac{\epsilon_\sigma^2}{\nu^2} \Phi^S dt) + (\bar{u} \circ \nabla_H^\perp) \nabla_H (\frac{\epsilon_\sigma}{\nu} \sigma^{\chi} dW_t - \frac{\epsilon_\sigma^2}{\nu^2} \chi^S dt)\vert_{z=\eta}. \label{eq:possibly-negligible-term}
\end{equation}
This term arises from developing the last term in the vertically averaged velocity $\bar{u}$ equation above. At this point, we can remark the following: for all twice differentiable function $f$, we have that
\begin{equation}
    (\bar{u} \bcdot \nabla_H^\perp) \nabla_H f -  (\bar{u} \bcdot \nabla_H) \nabla_H^\perp f = \begin{pmatrix}
        \bar{u}_y \\ - \bar{u}_x  
    \end{pmatrix} (\partial_{xx} + \partial_{yy})f = - u^\perp \Delta_H f, \label{eq:commutator-grad-gradPerp}
\end{equation}
where $u^\perp = (-u_y\:u_x)^\tr$. The RHS of the equation \eqref{eq:commutator-grad-gradPerp} above has the structure of a Coriolis force contribution, namely $f_0 (\bar{u} dt + \overline{\sigma^H dW_t})^\perp$, with $f_0$ denoting the Coriolis parameter. Since the Coriolis force is implicitly neglected in our study, we can neglect $u^\perp \Delta_H \sigma^\chi dW_t$ as well whenever
\begin{equation}
    u^\perp \Delta_H \sigma^\chi dW_t \lesssim f_0 \: \overline{\sigma^H dW_t}^\perp.
\end{equation}
This last relation can be interpreted in terms of scaling parameters,
\begin{equation}
    \frac{L}{T} \frac{L_\sigma^\chi}{L} \lesssim \frac{L_\sigma^\phi}{T_0} \Rightarrow \frac{L_\sigma^\chi}{L_\sigma^\phi} \lesssim \frac{T}{T_0} = \frac{1}{Ro},
\end{equation}
where $L_\sigma^\chi$ and $L_\sigma^\phi$ are the characteristic lengths associated to the rotational and potential components of the noise, $T_0 = f_0^{-1}$, $L$ and $T$ denote respectively the characteristic length and time of the problem, and $Ro = (f_0 T)^{-1}$ is the Rossby number. Hence, as is written, this system corresponds to a grey zone where the Rossby number is large enough to neglect the Coriolis effects, yet not large enough to neglect the rotational effects of the noise.

\begin{remark}
    Notice that the commutator term \eqref{eq:commutator-grad-gradPerp} could be kept in the dynamics in the more general case, and would not yield any particular difficulty. Yet, it is omitted for simplicity. This term can also be seen as a way to incorporate some Coriolis-type contribution in the equations, without introducing a Coriolis term initially in the LU Euler equations. Indeed, adding a Coriolis force would dramatically alter the derivation of the Water Waves equations, since the Coriolis force cannot be expressed as a gradient for potential flows. For this reason, it is typically not present in the classical Water Waves literature \cite{lannes2013water}.
\end{remark}

\noindent
Hence, $T_{SV}$ can be rewritten as follows,
\begin{align*}
    T_{SV} &= (\bar{u} \circ \nabla_H) \nabla_H (\frac{\epsilon_\sigma}{\nu}\sigma^\Phi dW_t - \frac{\epsilon_\sigma^2}{\nu^2} \Phi^S dt) + (\bar{u} \circ \nabla_H) \nabla_H^\perp (\frac{\epsilon_\sigma}{\nu} \sigma^{\chi} dW_t - \frac{\epsilon_\sigma^2}{\nu^2} \chi^S dt)\vert_{z=\eta} \\
    &= (\bar{u} \circ \nabla_H) \Big[\nabla_H (\frac{\epsilon_\sigma}{\nu}\sigma^\Phi dW_t - \frac{\epsilon_\sigma^2}{\nu^2} \Phi^S dt) + \nabla_H^\perp (\frac{\epsilon_\sigma}{\nu} \sigma^{\chi} dW_t - \frac{\epsilon_\sigma^2}{\nu^2} \chi^S dt)\vert_{z=\eta}\Big]\\
    &= (\bar{u} \circ \nabla_H) \Big[ \frac{\epsilon_\sigma}{\nu} \overline{\sigma^H dW_t} - \frac{\epsilon_\sigma^2}{\nu^2} \bar{u}_S dt\Big].
\end{align*}
Thus, in this context, $T_{SV}$ can be interpreted as the vertical average of $(u \circ \nabla_H) \sigma^H dW_t$. To see this, we write
\begin{equation*}
    \overline{(u \circ \nabla_H) \sigma^H dW_t} = (\bar{u} \circ \nabla_H) \overline{\sigma^H dW_t} + \overline{(\tilde{u} \circ \nabla_H) \tilde{\sigma}^H dW_t},
\end{equation*}
where we use the notation $\tilde{f} = f - \bar{f}$, for a square integrable function $f$. We remark then that $\overline{(\tilde{u} \circ \nabla_H) \tilde{\sigma}^H dW_t} = O(\mu)$ whenever $\partial_z \nabla_H \sigma^H dW_t$ is bounded. In the LU formalism, $(u \circ \nabla_H) \sigma^H dW_t$ is initially neglected by the very application of the Stochastic Reynolds Transport Theorem. We neglect $T_{SV}$ as well, so that the previous stochastic Saint-Venant system becomes,
\begin{align}
    d_t \eta &+ \frac{1}{\nu} \nabla_H \bcdot (h \bar{u}) dt\\
    &+ \frac{1}{\nu}\nabla_H \bcdot \Big(h \circ (\frac{\epsilon_\sigma}{\epsilon} \overline{\nabla_H \sigma^{\phi} dW_t} + \frac{\epsilon_\sigma}{\nu} \overline{\nabla_H^\perp \sigma^{\chi} dW_t} - \frac{\epsilon_\sigma^2}{\nu \epsilon} \overline{\nabla_H \phi^S} dt - \frac{\epsilon_\sigma^2}{\nu^2} \overline{\nabla_H^\perp \chi^{S}} dt) \Big) = 0, \nonumber\\
    d_t \bar{u} &+ \frac{\epsilon}{\nu} (\bar{u} \bcdot \nabla_H) \bar{u} \: dt + \nabla_H \eta \: dt \\
    &+ \Big(( \nabla_H (\frac{\epsilon_\sigma}{\nu}\sigma^\Phi dW_t - \frac{\epsilon_\sigma^2}{\nu^2} \Phi^S dt) + \nabla_H^\perp (\frac{\epsilon_\sigma}{\nu} \sigma^{\chi} dW_t - \frac{\epsilon_\sigma^2}{\nu^2} \chi^S dt)\vert_{z=\eta})  \circ \nabla_H \Big) \bar{u} = O(\frac{\mu}{\nu}, \mu \epsilon_\sigma). \nonumber
\end{align}
Therefore we recover the LU Saint-Venant system of \cite{DMM_2024}. Notice that we have used $\frac{1}{2} \nabla_H |\bar{u}|^2 = (\bar{u} \bcdot \nabla_H) \bar{u} + O(\mu)$, which holds since $\bar{u} = \nabla_H \Phi + O(\mu)$. However, the argument above to neglect the term $T_{SV}$ is not possible for the models below, since they are associated to higher order expansion in $\mu$.

\subsection{The Serre-Green-Naghdi system}

\noindent
To derive the Serre-Green-Naghdi system, we make a first order approximation in $\mu$ of the Water Waves one. As shown in \cite{lannes2013water}, $\nabla_H \Phi = \bar{u} - \frac{\mu}{3h} \nabla (h^3 \nabla \bcdot \bar{u}) + O(\mu^2)$. Thus we derive the following system,
\begin{align}
    d_t \eta &+ \frac{1}{\nu} \nabla_H \bcdot (h \bar{u}) dt\\
    &+ \frac{1}{\nu}\nabla_H \bcdot \Big(h \circ (\frac{\epsilon_\sigma}{\epsilon} \overline{\nabla_H \sigma^{\phi} dW_t} + \frac{\epsilon_\sigma}{\nu} \overline{\nabla_H^\perp \sigma^{\chi} dW_t} - \frac{\epsilon_\sigma^2}{\nu \epsilon} \overline{\nabla_H \phi^S} dt - \frac{\epsilon_\sigma^2}{\nu^2} \overline{\nabla_H^\perp \chi^{S}} dt) \Big) = 0, \nonumber\\
    d_t \bar{v} &+ \frac{\epsilon}{ \nu} \nabla_H \Big(\bar{u} \bcdot \bar{v} - \frac{|\bar{u}|^2}{2} + \eta - \frac{\mu}{2} (h \nabla_H \bcdot \bar{u})^2\Big) \\
    &+ \nabla_H \Big(\big( \nabla_H (\frac{\epsilon_\sigma}{\nu}\sigma^\Phi dW_t - \frac{\epsilon_\sigma^2}{\nu^2} \Phi^S dt) + \nabla_H^\perp (\frac{\epsilon_\sigma}{\nu} \sigma^{\chi} dW_t - \frac{\epsilon_\sigma^2}{\nu^2} \chi^S dt)\vert_{z=\eta}\big) \circ \bar{v}\Big) \nonumber \\
    &+ \epsilon \mu \nabla_H \Big[h (\bar{u} \bcdot \nabla_H \eta) \circ \nabla_H \bcdot \Big(\frac{\epsilon_\sigma}{\nu} \overline{\nabla_H \sigma^{\phi} dW_t} - \frac{\epsilon_\sigma^2}{\nu^2} \overline{\nabla_H \phi^S} dt\Big)\Big] = O(\frac{\mu^2}{\nu}, \mu^2 \epsilon_\sigma), \nonumber\\
    &\qquad \qquad \qquad \bar{v} = \bar{u} - \frac{\mu}{3h} \nabla_H (h^3 \nabla_H \bcdot \bar{u}), \quad h = 1+\epsilon \eta.
\end{align}
This system is slightly different from that of \cite{DMM_2024}. After developing the stochastic contribution, we observe that the difference stems from the presence of commutation terms in the fashion of $T_{SV}$. As mentioned above, they cannot be neglected \emph{a priori} through a scaling argument, because some terms of order $\mu$ remain, as shown above for the Saint-Venant system.

\subsection{The Boussinesq system}

\noindent
The Boussinesq system is obtained by disregarding the terms of order $\epsilon \mu$ in the Serre-Green-Naghdi system. Doing so yields the following,
\begin{align}
    d_t \eta &+ \frac{1}{\nu} \nabla_H \bcdot (h \bar{u}) dt\\
    &+ \frac{1}{\nu}\nabla_H \bcdot \Big(h \circ (\frac{\epsilon_\sigma}{\epsilon} \overline{\nabla_H \sigma^{\phi} dW_t} + \frac{\epsilon_\sigma}{\nu} \overline{\nabla_H^\perp \sigma^{\chi} dW_t} - \frac{\epsilon_\sigma^2}{\nu \epsilon} \overline{\nabla_H \phi^S} dt - \frac{\epsilon_\sigma^2}{\nu^2} \overline{\nabla_H^\perp \chi^{S}} dt) \Big) = 0, \nonumber\\
    d_t \bar{v} &+ \frac{\epsilon}{ \nu} \nabla_H \Big(\frac{|\bar{u}|^2}{2} + \eta \Big) + \nabla_H \Big(\bar{v} \circ \nabla_H (\frac{\epsilon_\sigma}{\nu}\sigma^\Phi dW_t - \frac{\epsilon_\sigma^2}{\nu^2} \Phi^S dt)\Big)\\
    &+ \nabla_H \Big(\bar{v} \circ \nabla_H^\perp (\frac{\epsilon_\sigma}{\nu} \sigma^{\chi} dW_t - \frac{\epsilon_\sigma^2}{\nu^2} \chi^S dt)\vert_{z=\eta}\Big) = O(\frac{\mu^2}{\nu}, \mu^2 \epsilon_\sigma, \mu \epsilon), \nonumber\\
    &\qquad \qquad \qquad \qquad \bar{v} = \bar{u} - \frac{\mu}{3h} \nabla_H (h^3 \nabla_H \bcdot \bar{u}), \quad h = 1+\epsilon \eta.
\end{align}
Notice that we have not given any constraint on $\epsilon_\sigma$ a priori to derive this interpretation of the Boussinesq system. If we additionally assume that $\epsilon_\sigma \sim \epsilon$, then the term $T_{SV}$ can be neglected by the same arguments as for the Saint-Venant system, since it is of order $O(\epsilon \mu)$. Therefore, the stochastic terms appearing in the momentum equation above can be slightly simplified,
\begin{align}
    d_t \eta &+ \frac{1}{\nu} \nabla_H \bcdot (h \bar{u}) dt \\
    &+ \frac{1}{\nu}\nabla_H \bcdot \Big(h \circ (\overline{\nabla_H \sigma^{\phi} dW_t} + \frac{\epsilon}{\nu} \overline{\nabla_H^\perp \sigma^{\chi} dW_t} - \frac{\epsilon}{\nu} \overline{\nabla_H \phi^S} dt - \frac{\epsilon^2}{\nu^2} \overline{\nabla_H^\perp \chi^{S}} dt) \Big) = 0, \nonumber\\
    d_t \bar{v} &+ \frac{\epsilon}{ \nu} \nabla_H \Big(\frac{|\bar{u}|^2}{2} + \eta \Big) + \nabla_H (\frac{\epsilon}{\nu}\sigma^\Phi dW_t - \frac{\epsilon^2}{\nu^2} \Phi^S dt) \circ \nabla_H \bar{u}\\
    &+ \nabla_H^\perp (\frac{\epsilon}{\nu} \sigma^{\chi} dW_t - \frac{\epsilon^2}{\nu^2} \chi^S dt)\vert_{z=\eta} \circ \nabla_H \bar{u} = O(\frac{\mu^2}{\nu}, \mu \epsilon), \nonumber\\
    &\qquad \qquad \qquad \qquad \bar{v} = \bar{u} - \frac{\mu}{3h} \nabla_H (h^3 \nabla_H \bcdot \bar{u}), \quad h = 1+\epsilon \eta.
\end{align}
This system is in turn exactly the stochastic Boussinesq model that is derived in \cite{DMM_2024}. Notice that another simplification comes from the fact that, now, $\bar{v}$ appears only in the evolution term $d_t \bar{v}$. Hence, once $\bar{v}$ is computed, $\bar{u}$ is computed by solving the elliptic problem
\begin{equation}
     \bar{v} = \bar{u} - \frac{\mu}{3h} \nabla (h^3 \nabla \bcdot \bar{u}).
\end{equation}

\subsection{Asymptotic behaviour of the noise}

\noindent
In the following, we seek a suitable scaling for the noise, so that it is still present in the three LU equations derived from Water Waves. For this purpose, we denote by $\sigma_{WW} dW_t$ the noise in the stochastic Water Waves equations, and by $\sigma_{SV} d\tilde{W}_t$ that of the Saint-Venant equation. Then, we impose \emph{a priori} the following structure for the noise,
\begin{gather}
    \hat{\sigma}_{WW}^{H} (x,y,z,x',y',z') = \frac{1}{\mu^{3/4}} \hat{\sigma}_{SV}^{H} (x,y, \mu^{-1/2} z,x',y', \mu^{-1/2} z'), \\
    \hat{\sigma}_{WW}^{z} (x,y,z,x',y',z') = \frac{\alpha_{WW}^z}{\mu^{1/4}} \hat{\sigma}_{SV}^{z} (x,y, \mu^{-1/2} z,x',y', \mu^{-1/2} z'),
\end{gather}
where $\alpha_{WW}^z$ denotes a scaling parameter, which depends on $\mu$ \emph{a priori}, for the vertical component of the noise. We assume in the following that $\alpha_{WW}^z \rightarrow 0$ as $\mu \rightarrow 0$. 

Next, we define a scaled noise $\sigma_\mu dW_t$ for the scaled Water Waves equations, which is so that
\begin{align}
    \sigma_\mu^H dW_t(x,y,z,t) &= \sigma_{WW}^{H} dW_t(x,y,\mu^{1/2} z,t),\\
    \sigma_\mu^z dW_t(x,y,z,t) &= \frac{1}{\mu^{1/2}} \sigma_{WW}^{z} dW_t^z(x,y,\mu^{1/2} z,t).
\end{align}
As such, $\sigma_\mu dW_t$ is divergence-free whenever $\sigma_{WW} dW_t$ is. Furthermore,
\begin{align*}
    \sigma_\mu^H (x,y,z) dW_t &= \int_{\mathcal{S_\mu}} \hat{\sigma}_{SV}^{H}(x,y,\mu^{1/2} z,x',y',z') dW_t(x',y',z') dx' dy'dz'\\
    &=\mu^{1/2} \int_{\mathcal{S}} \hat{\sigma}_{SV}^{H} (x,y,\mu^{1/2} z,x',y',\mu^{1/2} z'') dW_t(x',y',\mu^{1/2} z'') dx' dy'dz''\\
    &= \mu^{3/4} \int_{\mathcal{S}} \hat{\sigma}_{SV}^{H} (x,y,\mu^{1/2} z,x',y',\mu^{1/2} z'') d\tilde{W}_t(x',y', z'') dx' dy'dz''\\
    &= \sigma_{SV}^H(x,y, z) d\tilde{W}_t.
\end{align*}
Similarly,
\begin{align*}
    \sigma_\mu^z (x,y,z) dW_t &= \frac{\alpha_{WW}^z}{\mu^{1/2}} \int_{\mathcal{S_\mu}} \hat{\sigma}_{SV}^{z}(x,y,\mu^{1/2} z,x',y',z') dW_t(x',y',z') dx' dy'dz'\\
    &= \alpha_{WW}^z \int_{\mathcal{S}} \hat{\sigma}_{SV}^{z} (x,y,\mu^{1/2} z,x',y',\mu^{1/2} z'') dW_t(x',y',\mu^{1/2} z'') dx' dy'dz''\\
    &= \alpha_{WW}^z \mu^{1/4} \int_{\mathcal{S}} \hat{\sigma}_{SV}^{z} (x,y,\mu^{1/2} z,x',y',\mu^{1/2} z'') d\tilde{W}_t(x',y', z'') dx' dy'dz''\\
    &= \alpha_{WW}^z \sigma_{SV}^z(x,y, z) d\tilde{W}_t.
\end{align*}
Since $\alpha_{WW}^z \rightarrow 0$ and $\sigma_{SV} dW_t$ is divergence-free, it follows that $\sigma_{SV}^H dW_t$ is divergence-free as a 2D vector.

\section{2D Numerical simulations}

\noindent
In this section, we present some numerical simulations and results on the previously derived models. The domain of interest is typically the 2D torus, where the variables $\eta$, $\Phi$ and $\bar{u}$ are defined and computed. The Python algorithms have been adapted from \cite{website_Haefner_2026} and are available at the following link: \hyperlink{https://github.com/Lord-Fufu/stochastic-wave-equations}{https://github.com/Lord-Fufu/stochastic-wave-equations}.

For simplicity, we disregard the It\={o}-Stokes drift terms to focus on the Stratonovitch contributions. A more precise focus on the effect of these contributions would be the subject of another work. Importantly, the Water Waves equation are numerically solved using a ``direct'' method that is presented below, while the other models are solved using pseudo-spectral methods.

\subsection{Presentation of the numerical schemes}

\noindent
We introduce hereafter the numerical schemes we used for solving the four models presented before, which are the LU interpretations of Water Waves, Saint-Venant, Boussinesq and Serre-Green-Naghdi.

\subsubsection{The Water Waves system}

\noindent
Let us present first the method we use for solving numerically the Water Waves equations. As mentioned earlier, the elevation $\eta$ and the surface potential $\Phi$ are solved on a C-grid, both at the centre of the cells, while their gradient are computed on the edges.

\begin{figure}[h!]
\centering
\begin{tikzpicture}[scale=2]

    \tikzstyle{fleche}=[->,>=stealth,thick];
    \tikzstyle{doubleFleche}=[<->,>=stealth,thick];

    \tikzstyle{stylebloc}=[draw,rectangle,thick];
    
    \draw[thick] (-4cm, 0cm) -- (0cm, 0cm);
    \draw[thick] (-4cm,-1cm) -- (0cm,-1cm);
    \draw[thick] (-4cm,-2cm) -- (0cm,-2cm);

    \draw[thick] (-3cm,0.5cm) -- (-3cm,-2.5cm) ;
    \draw[thick] (-2cm,0.5cm) -- (-2cm,-2.5cm) ;
    \draw[thick] (-1cm,0.5cm) -- (-1cm,-2.5cm) ;
    
    \fill[opaque, black]  (-2.6,-1.4) rectangle (-2.4,-1.6);
    \fill[opaque, black]  (-1.6,-1.4) rectangle (-1.4,-1.6); 
    \fill[opaque, black]  (-1.6,-0.4) rectangle (-1.4,-0.6); 
    \fill[opaque, black]  (-2.6,-0.4) rectangle (-2.4,-0.6);


    \node[mark=*, mark size=3pt] at (-2.5,0) {\pgfuseplotmark{*}};
    \node[mark=*, mark size=3pt] at (-2.5,-1) {\pgfuseplotmark{*}};
    \node[mark=*, mark size=3pt] at (-2.5,-2) {\pgfuseplotmark{*}};
    \node[mark=*, mark size=3pt] at (-1.5,0) {\pgfuseplotmark{*}};
    \node[mark=*, mark size=3pt] at (-1.5,-1) {\pgfuseplotmark{*}};
    \node[mark=*, mark size=3pt] at (-1.5,-2) {\pgfuseplotmark{*}};
    
    \node[mark=*, mark size=3pt] at (-3,-0.5) {\pgfuseplotmark{triangle*}};
    \node[mark=*, mark size=3pt] at (-2,-0.5) {\pgfuseplotmark{triangle*}};
    \node[mark=*, mark size=3pt] at (-1,-0.5) {\pgfuseplotmark{triangle*}};
    \node[mark=*, mark size=3pt] at (-3,-1.5) {\pgfuseplotmark{triangle*}};
    \node[mark=*, mark size=3pt] at (-2,-1.5) {\pgfuseplotmark{triangle*}};
    \node[mark=*, mark size=3pt] at (-1,-1.5) {\pgfuseplotmark{triangle*}};

    \node[black] at (-5, -0.5) {$\partial_y \Phi, \partial_y \eta$};
    \draw[dotted, line width=0.5mm] (-2.5cm,-1cm) -- (-4.6cm,-0.5cm) ;
    
    \node[black] at (1, -1) {$\partial_x \Phi, \partial_x \eta$};
    \draw[dotted, line width=0.5mm] (-1cm,-0.5cm) -- (0.6cm,-1cm) ;
    
    \node[black] at (1, -2) {$\Phi, \eta$};
    \draw[dotted, line width=0.5mm] (-1.5cm,-1.5cm) -- (0.7cm,-2cm) ;

\end{tikzpicture}

\caption[Representation of the C-grid for for solving the Water Waves system.]{Representation of the C-grid for solving the Water Waves system. The variables $\Phi, \eta$ are computed at the center of cells (squares), and their gradients on the edges (circles and triangles).}
\end{figure}
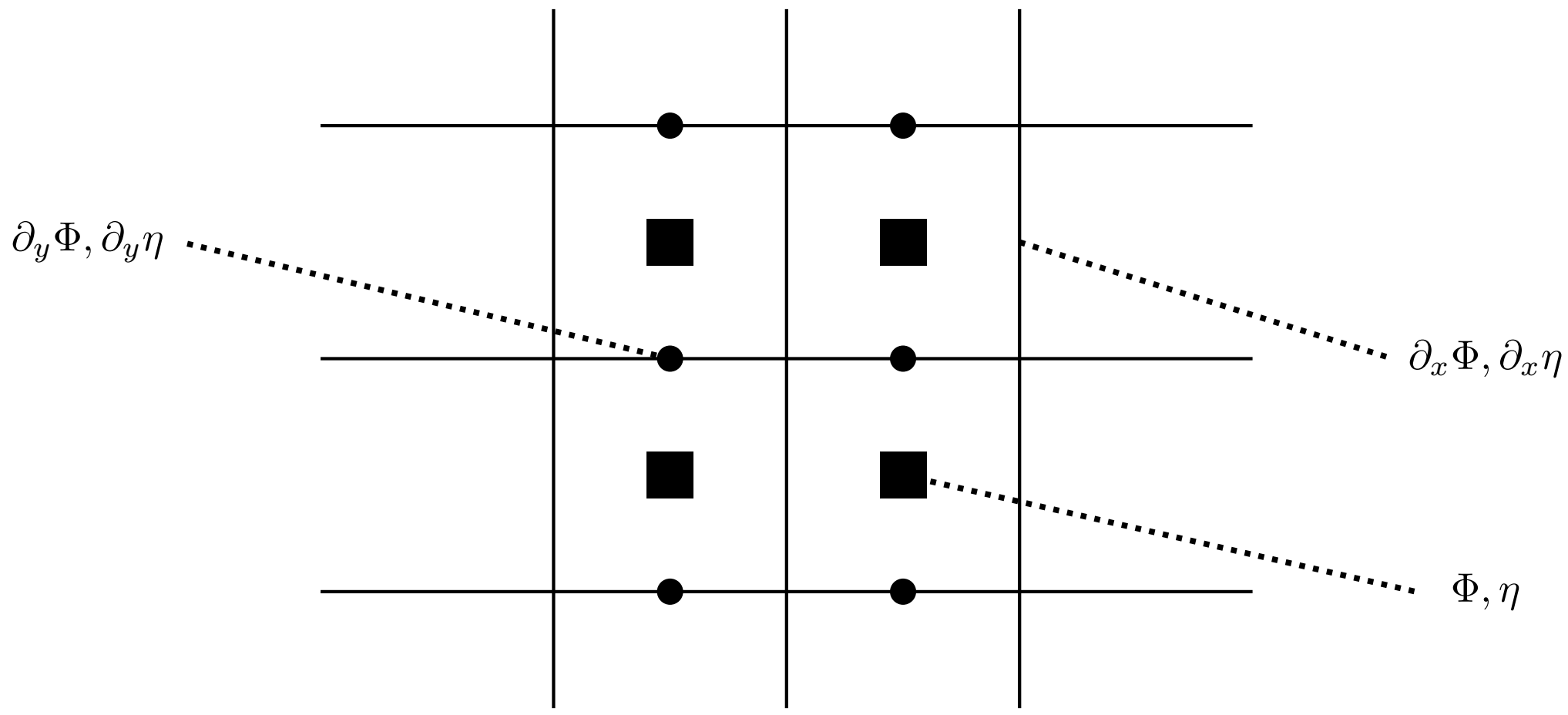

\noindent
Before presenting the numerical simulations, we briefly explain the method used to solve the Water Waves equations. As previously mentioned, the method is ``direct'', as opposed to the pseudo-spectral methods used for the three other derived models.

\paragraph{Computation of the Dirichlet-to-Neumann operator.}
In order to ensure the conservation of mass, we compute the Dirichlet-to-Neumann operator as follows,
\begin{equation}
    G(\eta)\Phi = - \nabla_H \bcdot \Big(\int_{-h_0}^\eta \nabla_H \phi dz'\Big).
\end{equation}
As such, $G(\eta)\Phi$ reads as the divergence of a flux, and thus the evolution equation of $\eta$ is in conservative form.

To compute $G(\eta)\Phi$, we must solve the following elliptic problem at each timestep, the surface potential $\Phi$ being given,
\begin{equation}
    (\mu \Delta_H + \partial_{zz}) \phi = 0, \quad \partial_z \phi = 0 \text{ at $z=-h_0$}, \quad \phi = \Phi \text{ at $z=\eta$.}
\end{equation}
To solve this equation, we rewrite the problem in terms of scaled space variables. The reasoning is similar as in \cite{lannes2013water}. Let $\tilde{\phi}(\tilde{x},\tilde{y},\tilde{z}) = \phi(l_x \tilde{x}, l_y \tilde{y}, (h_0+\eta) \tilde{z})$, so that
\begin{equation}
    \phi(x,y,z) = \tilde{\phi}\Big(\frac{x}{l_x}, \frac{y}{l_y}, \frac{z}{h_0+\eta}\Big).
\end{equation}
Then,
\begin{multline}
    0 = \Delta_3 \phi = \Big(\frac{1}{l_x^2} \partial_{\tilde{x}\tilde{x}} + \frac{1}{l_y^2} \partial_{\tilde{y}\tilde{y}} + \frac{1}{(h_0+\eta)^2} \partial_{\tilde{z}\tilde{z}}\Big) \tilde{\phi} - \frac{2\tilde{z}}{(h_0+\eta)} (\frac{1}{l_x} \partial_{\tilde{z}\tilde{x}} + \frac{1}{l_y} \partial_{\tilde{z}\tilde{y}}) \tilde{\phi}\\
    + \frac{\tilde{z}^2 |\nabla_H \eta|^2}{(h_0+\eta)^2} \partial_{\tilde{z}\tilde{z}} \tilde{\phi} - \frac{\tilde{z} \Delta_H \eta}{(h_0+\eta)} \partial_{\tilde{z}} \tilde{\phi} + \frac{2 \tilde{z} |\nabla_H \eta|^2}{(h_0+\eta)^2} \partial_{\tilde{z}} \tilde{\phi}.
\end{multline}
Hence, the Dirichlet-to-Neumann operator is computed by solving first the equation above in the 3D domain, subject to the horizontal periodic boundary conditions, together with $\tilde{\phi}\vert_{z=\eta} = \Phi$ and $\tilde{\phi}\vert_{z=-h_0}$. Then, the gradient of the 3D solution $\tilde{\phi}$ is integrated along $\tilde{z}$ to infer $G(\eta)\Phi$.

In addition, the time scheme that we use is intermediate between a deterministic Runge-Kutta method of order 4, and a ``double advection'' method, which we could also refer to as ``stochastic Euler-Heun'' method. It is essentially in the same spirit as what is presented in \cite{DMM_2024}. To be more precise:
\begin{itemize}
    \item[$\bullet$] on the one hand, the non-stochastic terms -- i.e. not explicitly dependent of $\sigma dW_t$ -- are computed as in the standard fourth order Runge-Kutta approach,
    \item[$\bullet$] on the other hand, the stochastic terms are computed via the ``stochastic Euler-Heun'' method, as presented in \cite{XT_2016}. Notice that we have neglected the It\={o}-Stokes drift contributions, therefore stochasticity only appears through the Stratonovitch products. 
\end{itemize}
Such type of schemes have been studied in more details in \cite{GF_2019, PS_2008, FBLM_2023}.

\subsubsection{Saint-Venant, Boussinesq and Serre-Green-Naghdi systems}

\noindent
Remind that we derived three classical wave models from the Water Waves one, which are Saint-Venant, Boussinesq and Serre-Green-Naghdi. To solve them numerically, we use a C-grid again, but this time in terms of the elevation $\eta$ and the momentum $h\bar{u}$. While $\eta$ is computed at the centre of the cells again, the variable $h\bar{u} = (h\bar{u}_x, h\bar{u}_y)^\tr$ is computed on the edges. We represent hereafter the points where the gradients are computed. 

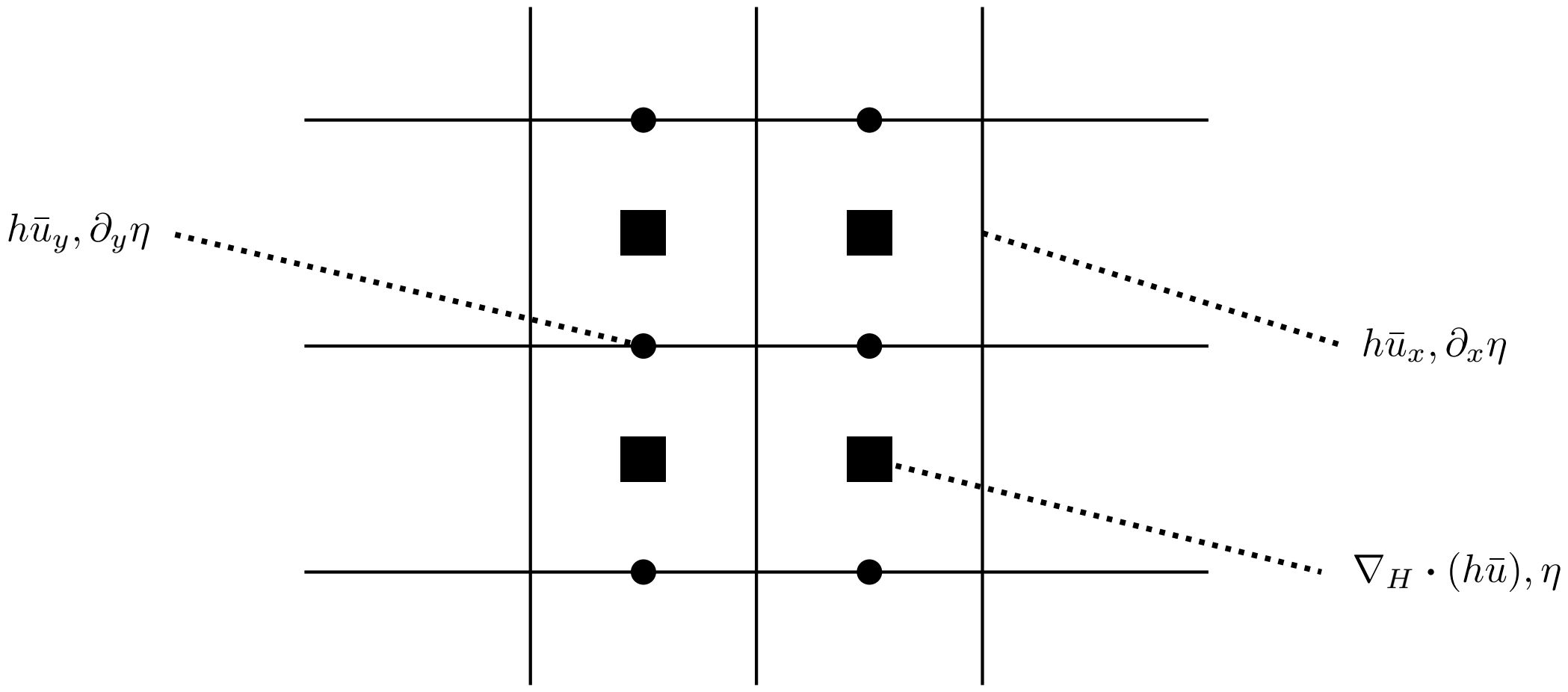
\begin{figure}[h!]
\centering
\begin{tikzpicture}[scale=2]

    \tikzstyle{fleche}=[->,>=stealth,thick];
    \tikzstyle{doubleFleche}=[<->,>=stealth,thick];

    \tikzstyle{stylebloc}=[draw,rectangle,thick];
    
    \draw[thick] (-4cm, 0cm) -- (0cm, 0cm);
    \draw[thick] (-4cm,-1cm) -- (0cm,-1cm);
    \draw[thick] (-4cm,-2cm) -- (0cm,-2cm);

    \draw[thick] (-3cm,0.5cm) -- (-3cm,-2.5cm) ;
    \draw[thick] (-2cm,0.5cm) -- (-2cm,-2.5cm) ;
    \draw[thick] (-1cm,0.5cm) -- (-1cm,-2.5cm) ;
    
    \fill[opaque, black]  (-2.6,-1.4) rectangle (-2.4,-1.6);
    \fill[opaque, black]  (-1.6,-1.4) rectangle (-1.4,-1.6); 
    \fill[opaque, black]  (-1.6,-0.4) rectangle (-1.4,-0.6); 
    \fill[opaque, black]  (-2.6,-0.4) rectangle (-2.4,-0.6);


    \node[mark=*, mark size=3pt] at (-2.5,0) {\pgfuseplotmark{*}};
    \node[mark=*, mark size=3pt] at (-2.5,-1) {\pgfuseplotmark{*}};
    \node[mark=*, mark size=3pt] at (-2.5,-2) {\pgfuseplotmark{*}};
    \node[mark=*, mark size=3pt] at (-1.5,0) {\pgfuseplotmark{*}};
    \node[mark=*, mark size=3pt] at (-1.5,-1) {\pgfuseplotmark{*}};
    \node[mark=*, mark size=3pt] at (-1.5,-2) {\pgfuseplotmark{*}};
    
    \node[mark=*, mark size=3pt] at (-3,-0.5) {\pgfuseplotmark{triangle*}};
    \node[mark=*, mark size=3pt] at (-2,-0.5) {\pgfuseplotmark{triangle*}};
    \node[mark=*, mark size=3pt] at (-1,-0.5) {\pgfuseplotmark{triangle*}};
    \node[mark=*, mark size=3pt] at (-3,-1.5) {\pgfuseplotmark{triangle*}};
    \node[mark=*, mark size=3pt] at (-2,-1.5) {\pgfuseplotmark{triangle*}};
    \node[mark=*, mark size=3pt] at (-1,-1.5) {\pgfuseplotmark{triangle*}};

    \node[black] at (-5, -0.5) {$h\bar{u}_y, \partial_y \eta$};
    \draw[dotted, line width=0.5mm] (-2.5cm,-1cm) -- (-4.6cm,-0.5cm) ;
    
    \node[black] at (1, -1) {$h\bar{u}_x, \partial_x \eta$};
    \draw[dotted, line width=0.5mm] (-1cm,-0.5cm) -- (0.6cm,-1cm) ;
    
    \node[black] at (1.1, -2) {$\nabla_H \bcdot (h\bar{u}), \eta$};
    \draw[dotted, line width=0.5mm] (-1.5cm,-1.5cm) -- (0.5cm,-2cm) ;

\end{tikzpicture}

\caption[Representation of the C-grid for for solving the Saint-Venant, Boussinesq and Serre-Green-Naghdi systems.]{Representation of the C-grid for solving the Saint-Venant, Boussinesq and Serre-Green-Naghdi systems. The variable $\eta$ is computed at the center of cells (squares), and $h\bar{u} = (h\bar{u_x}, h\bar{u}_y)^\tr$ on the edges (circles and triangles).}
\end{figure}

The main difference compared to the Water Waves simulation is that we use pseudo-spectral methods to compute the aforementioned models. Nevertheless, the terms involving Stratonovitch products are again computed with a double advection method, a fourth order Runge-Kutta method being used for the remaining terms  -- see \cite{XT_2016}.

\subsection{Heap-of-water simulations}

\noindent
In the following, we present and compare simulations of the different models with a \emph{heap-of-water} type initial condition.

\noindent
The setting for the Water Waves system is the following,
\begin{itemize}
    \item[$\bullet$] the horizontal domain is a square of 3.2 km side length, which is discretised using a square mesh of format $32 \times 32$, with spatial steps $dx = dy = 0.1 km$, 
    \item[$\bullet$] the depth of the basin is assumed to be constant and equal to $100 m$,
    \item[$\bullet$] the height of the initial condition is $10 m$ in our study.
\end{itemize}
More precisely, the shape of the initial condition is given by the following formula,
\begin{equation}
    \eta_0 = A \exp\Bigg[-10^4 \times \Big[(\frac{x}{l_x} - 1/2)^2 + (\frac{y}{l_y} - 1/2)^2\Big]^2\Bigg], \quad \Phi_0 = 0 \text{ or } \bar{u}_0 = 0,
\end{equation}
where $A=10m$, $l_x = l_y = 3.2 km$.
\begin{figure}[h!]
    \centering
    \begin{tabular}{c|c}
        \includegraphics[trim={1cm 1cm 3cm 2.5cm}, clip, width=6.05cm]{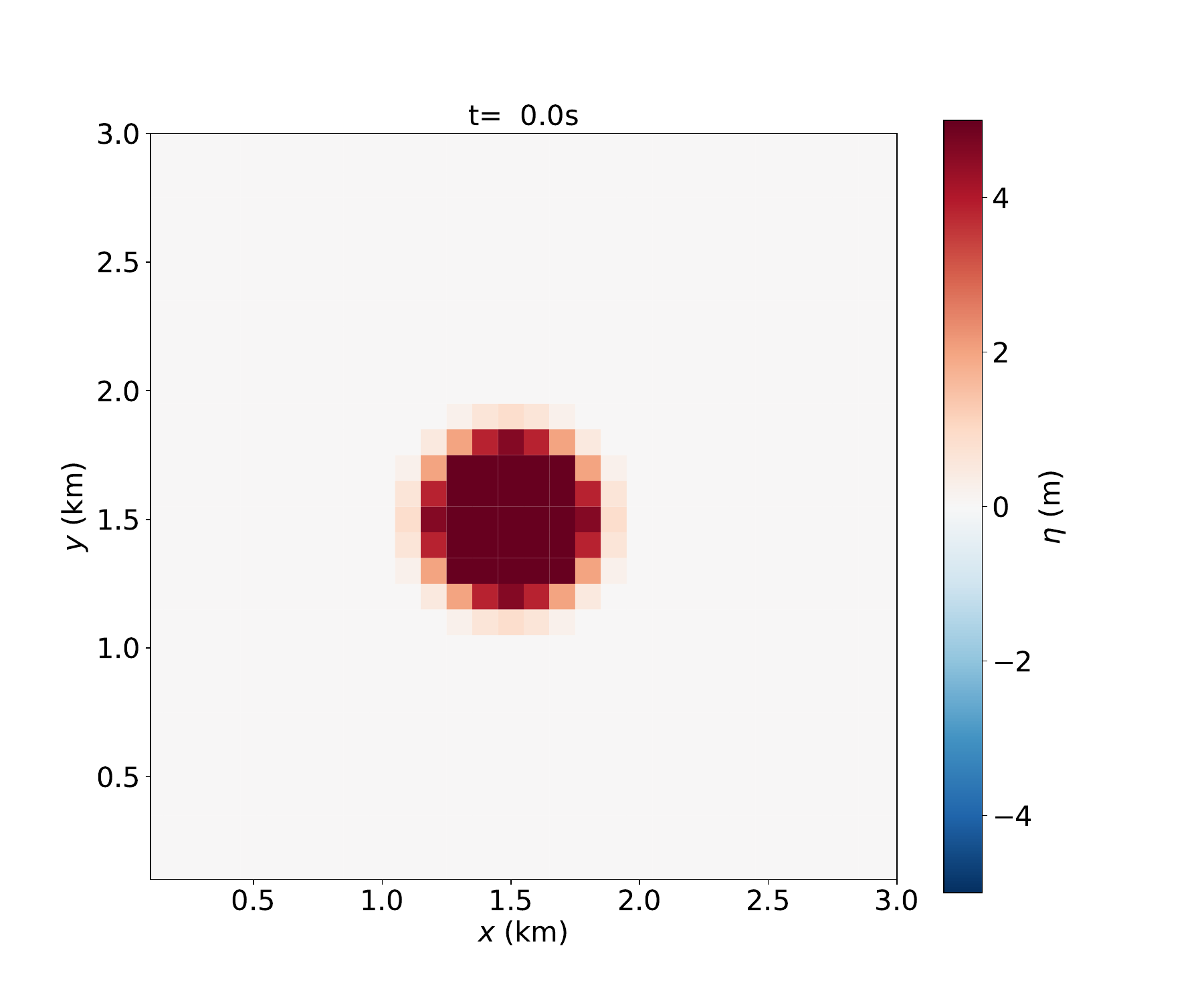} & \includegraphics[trim={17cm 0.3cm 3cm 0}, clip, width=8cm]{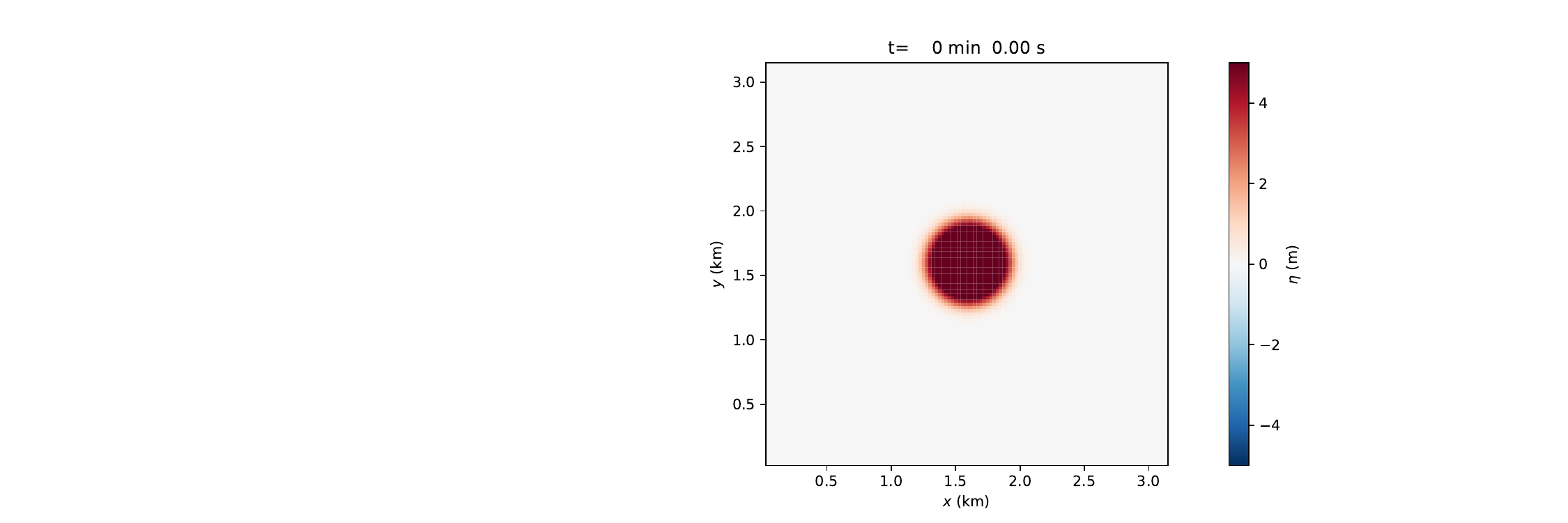}
    \end{tabular}
    \caption[Initial elevation projected onto two horizontal grids.]{Initial elevation $\eta_0$, projected over the $32 \times 32$ horizontal grid (left) and over the $128 \times 128$ one (right).}
\end{figure}

\noindent
In this setting, the values of the nonlinearity and shallowness parameters are $\epsilon = 0.1$ and $\mu = (1/32)^2 \approx 10^{-3}$, respectively. Since our study focuses on the shallow water limit of the Water Waves equations, we consider $\nu = 1$ in the following. Furthermore, as presented above, the noise is decomposed into potential and rotational contributions, i.e.
\begin{equation}
    \sigma dW_t = \sqrt{\Upsilon_\phi} \nabla \sigma^\phi dW_t + \sqrt{\Upsilon_\chi} \nabla_H^\perp \sigma^\chi dW_t,
\end{equation}
where $\Upsilon_\phi$ and $\Upsilon_\chi$ are scaling parameters for the noise. In our study, we choose $\Upsilon_\phi = 10^4 m^2/s$ and $\Upsilon_\chi = 10^8 m^2/s$, and we define
\begin{equation}
    \sigma^\phi dW_t = \sum_{k=1}^4 \phi_k d\beta^k, \quad \sigma^\chi dW_t = \sum_{k=5}^8 \chi_k d\beta^k,
\end{equation}
where
\begin{gather*}
    \phi_1(x,y) = \sin(k(x/l_x+y/l_y)), \quad \phi_3(x,y) = \cos(k(x/l_x+y/l_y)), \\
    \phi_2(x,y) = \sin(k(x/l_x-y/l_y)), \quad \phi_4(x,y) = \cos(k(x/l_x-y/l_y)),
\end{gather*}
and
\begin{gather*}
    \chi_5(x,y) = \sin(k(x/l_x+y/l_y)), \quad \chi_7(x,y) = \cos(k(x/l_x+y/l_y)), \\
    \chi_6(x,y) = \sin(k(x/l_x-y/l_y)), \quad \chi_8(x,y) = \cos(k(x/l_x-y/l_y)),
\end{gather*}
with $k = 2 \pi$.

\noindent
For the other three systems derived from Water Waves, the setting is similar, except for the fact that we investigate their behaviour on a thinner mesh over the same horizontal domain. It is of format $128 \times 128$, thus with spatial steps $dx = dy = 0.025 km$. Moreover, the LU noise is chosen as follows,
\begin{align*}
    \sigma dW_t := \sqrt{\Upsilon_{LU}} f_{\alpha}^k dW_t,
\end{align*}
where
\begin{equation}
    f_{\alpha}^k dW_t := \frac{1}{C} \mathcal{F}_{2D}^{-1} \Bigg[\exp\Big(- \alpha (|\xi|-k)^2\Big) \begin{pmatrix}
        - \xi_y \\ \xi_x
    \end{pmatrix} dW_t\Bigg].
\end{equation}
Here, $\mathcal{F}_{2D}^{-1}$ denotes the classical inverse 2D Fourier transform. In addition, $\Upsilon_{LU}$ is a scaling parameter for the noise in the Saint-Venant system, taken such that $\Upsilon_{LU} = 5.0 \: 10^1$. The parameter $k$ is interpreted as the peak wavenumber of the noise, and is chosen to be $2\pi$, and $C$ is a normalisation constant. To make the noise close to a space monochromatic wave, we set $\alpha = 10^7$.
 
\subsubsection{One trajectory comparison}

\noindent
We compare hereafter the deterministic Water Waves solution to a single trajectory of its LU counterpart.
\begin{figure}[h!]
    \centering
    \begin{tabular}{c|c}
        \includegraphics[trim={10cm 0 4cm 0}, clip, width=6.5cm]{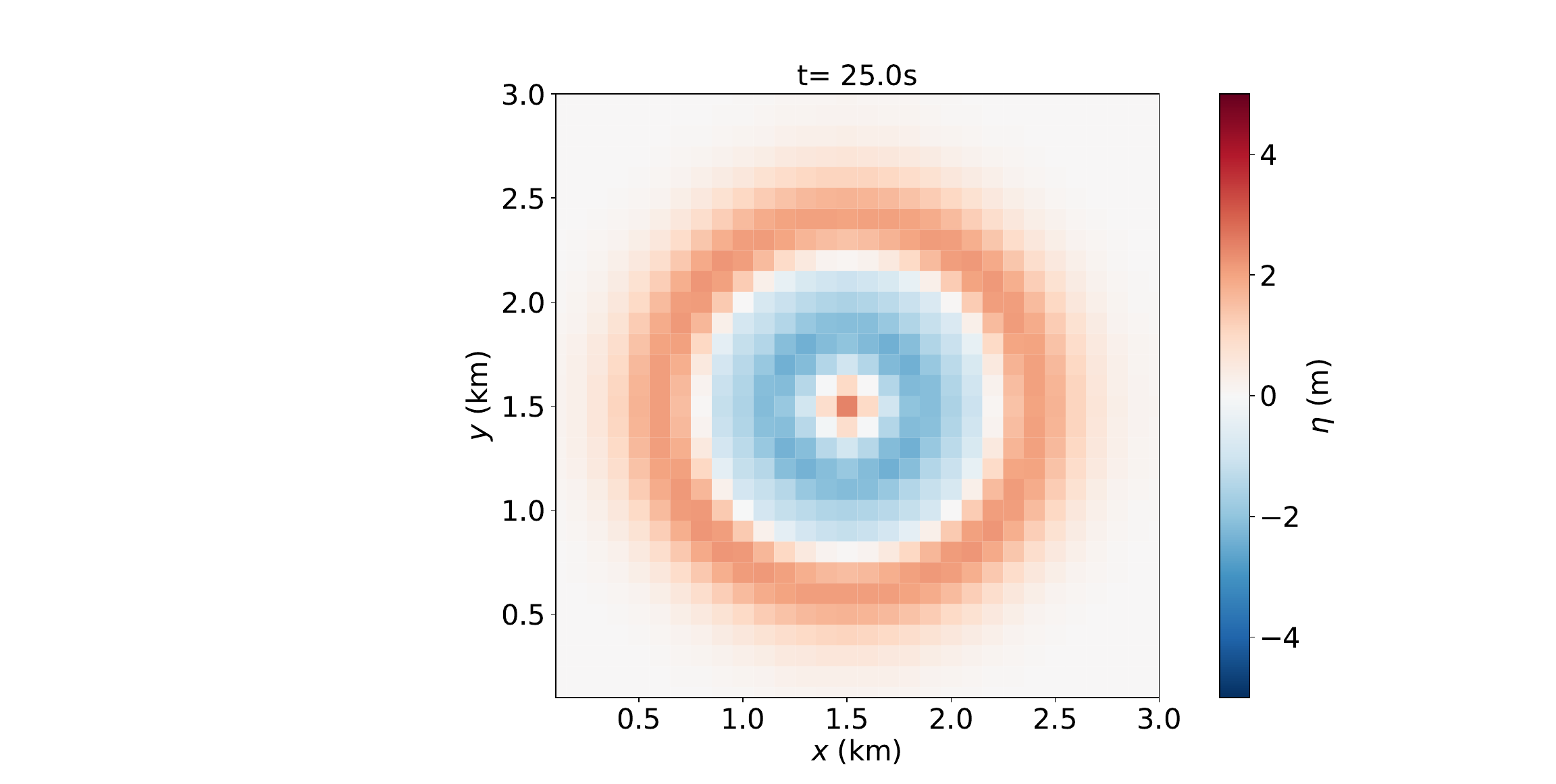} & \includegraphics[trim={10cm 0 4cm 0}, clip, width=6.5cm]{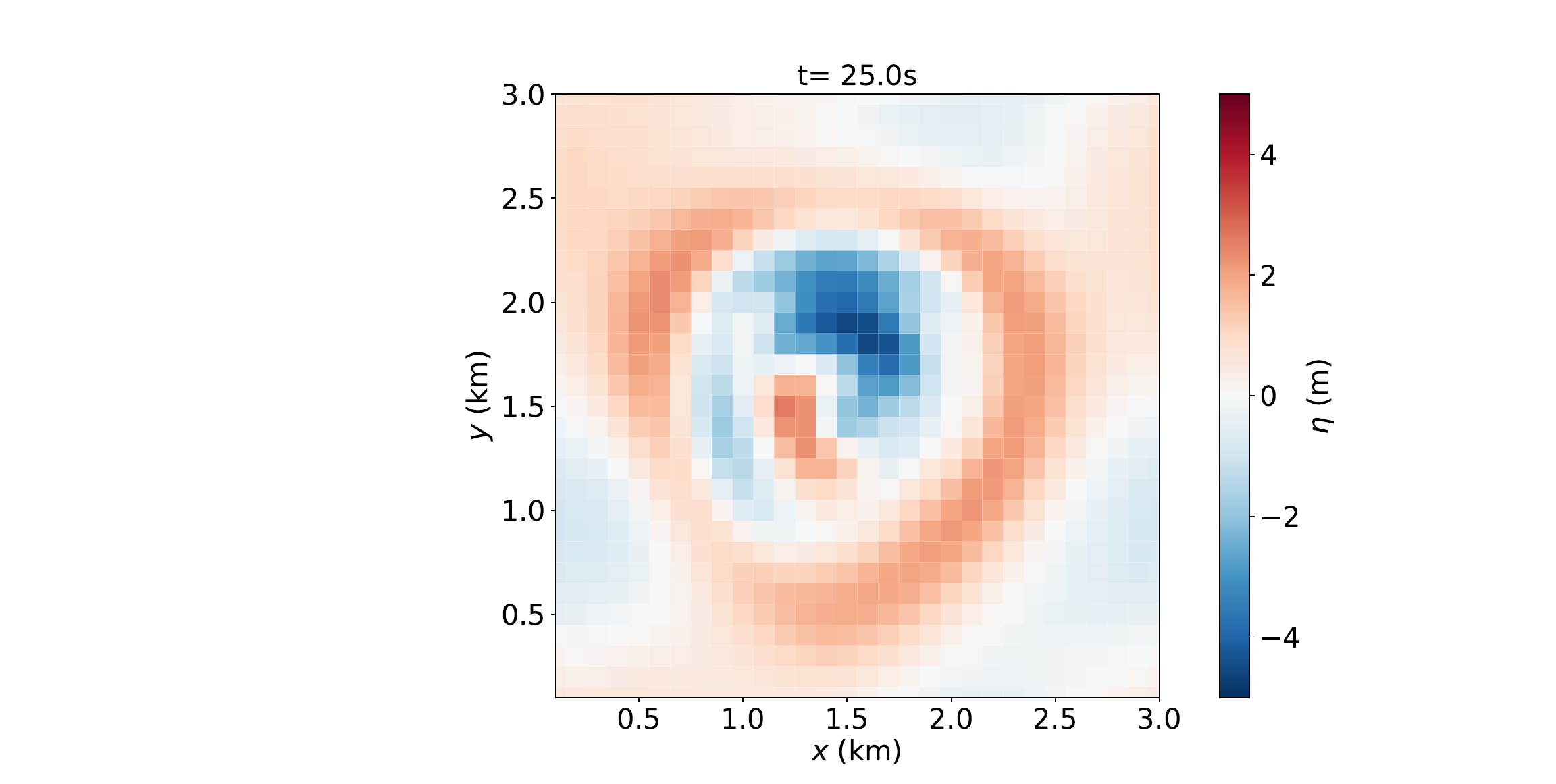}\\
        \hline
        \includegraphics[trim={10cm 0 4cm 0}, clip, width=6.5cm]{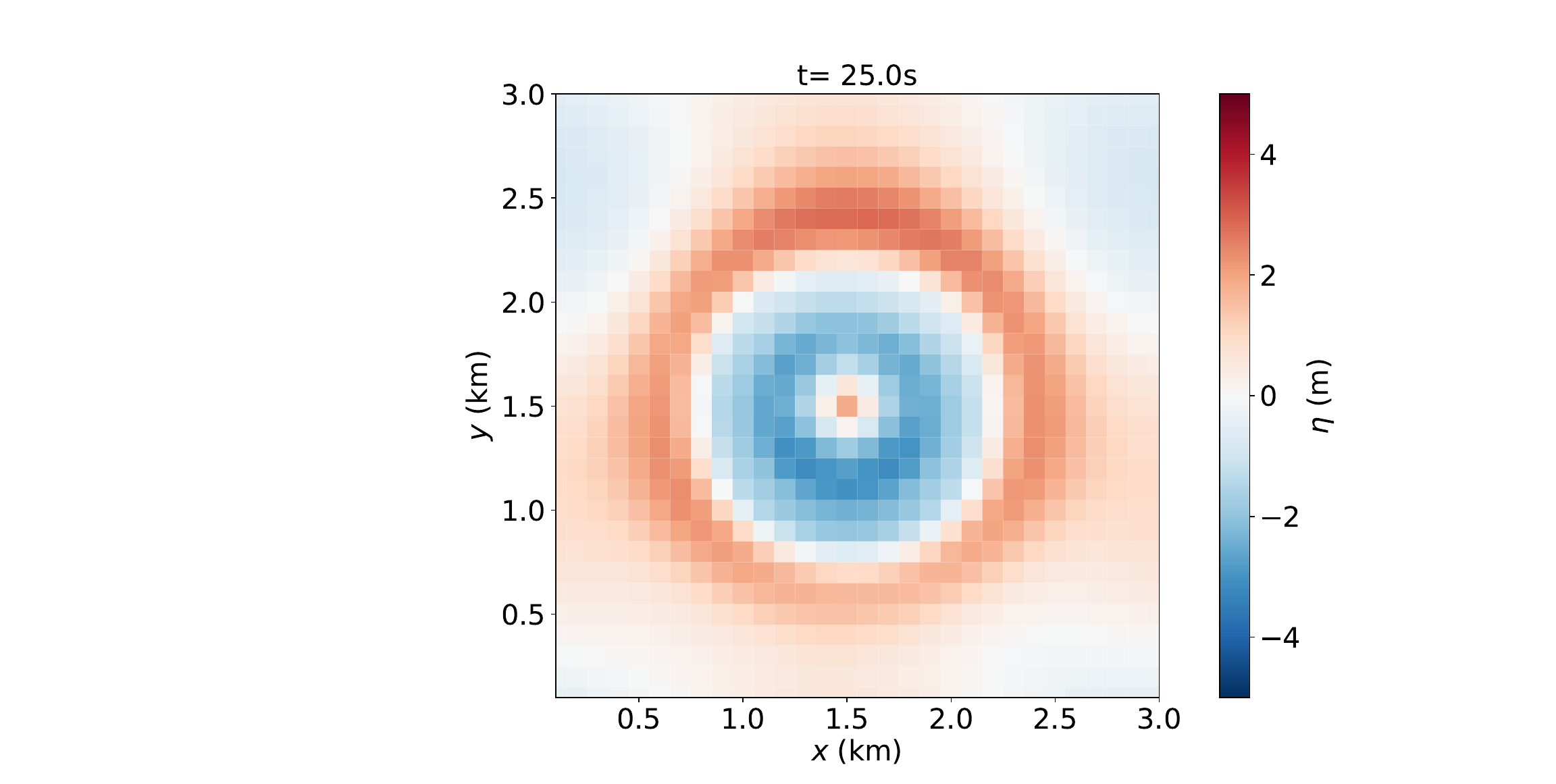} & \includegraphics[trim={10cm 0 4cm 0}, clip, width=6.5cm]{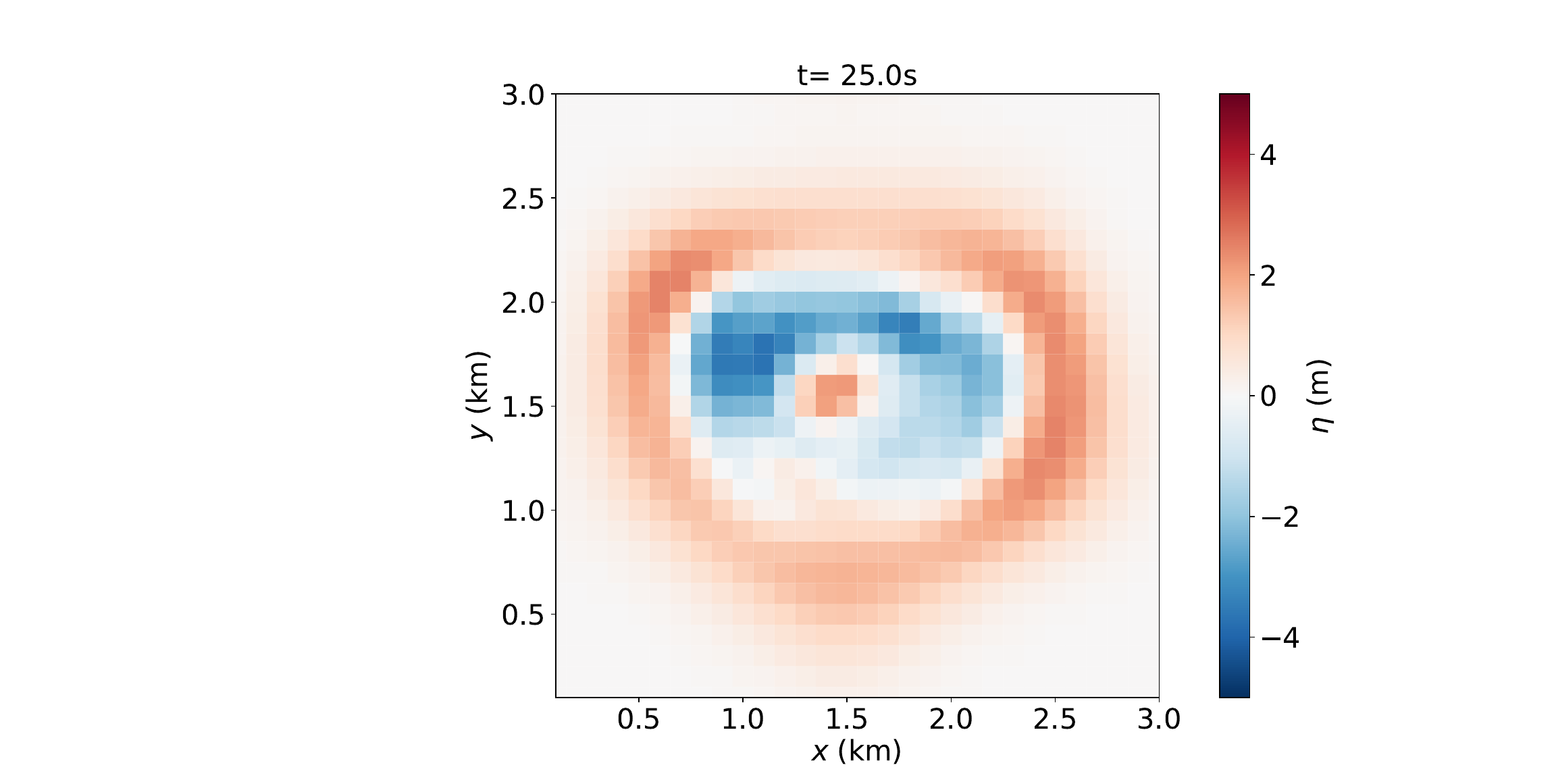}
    \end{tabular}
    \caption[Elevation of the deterministic Water Waves solution and its stochastic counterpart for different types of noises, at fixed time.]{Elevation of the deterministic Water Waves solution (upper left), and of three trajectories of the LU Water Waves system for different types of noises, at $t=25s$. The LU interpretations are proposed with potential noise (lower left), rotational noise (lower right) and both noises (upper right).} \label{fig:plots-WW} 
\end{figure}
First, we observe on Figure \ref{fig:plots-WW} that the introduction of a noise breaks the symmetry that is present in the deterministic context -- see \cite{DMM_2024}. However, adding a potential or a rotational noise does not impact the system in the same way. We observe that potential noise affects the water elevation but conserves the ``circular'' shape of the wave, while the rotational component modifies dramatically the elevation profile. More specifically, the ``bounce'' of the initial condition is visible in the deterministic setting. This ``bounce'' is shifted due to rotational noises in the stochastic interpretation, yet, in comparison, it is left almost unchanged by potential noises. Moreover, potential noises affect the downstream of the wave while rotational noise do not. This is due to the vertical velocity term $(\partial_z \phi_k^\sigma)\vert_{z=\eta} d\beta^k$ involved in the expression of $G(\eta) \Phi_k^\sigma$, as it is not a transport noise. In fact, this term can be considered as an additive noise term in the approximation $\eta = 0$, which is valid in particular at the initial time $t=0$ ``away'' from the heap-of-water. For this reason, potential noises tend to have a ``blurring'' effect on the elevation profile. However, this term cannot be considered additive for the whole simulation, since it depends on $\eta$ explicitly.

\begin{figure}[H]
    \centering
    \begin{tabular}{c|c}
        \includegraphics[trim={10cm 0 4cm 0}, clip, width=6.5cm]{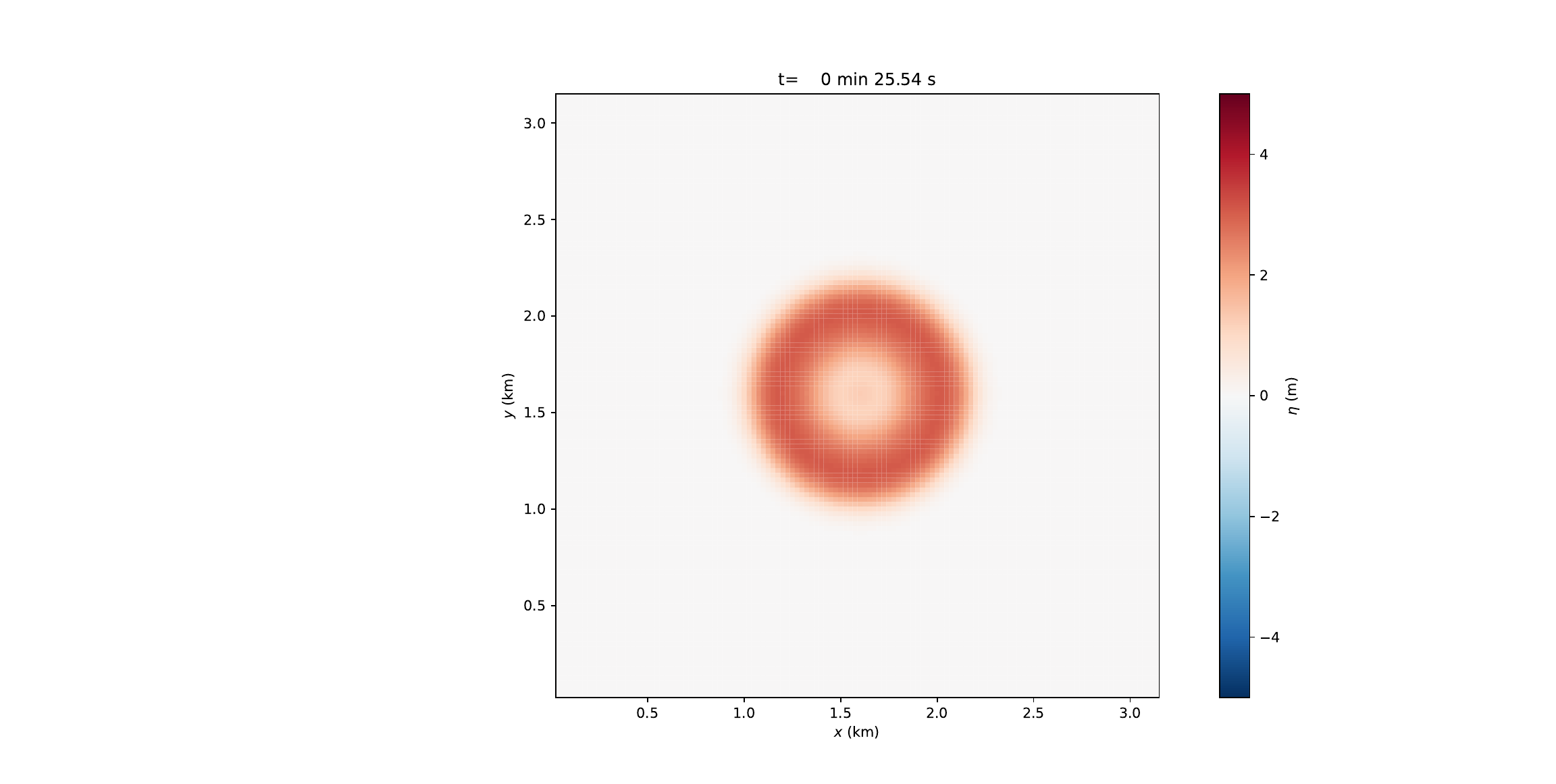} & \includegraphics[trim={10cm 0 4cm 0}, clip, width=6.5cm]{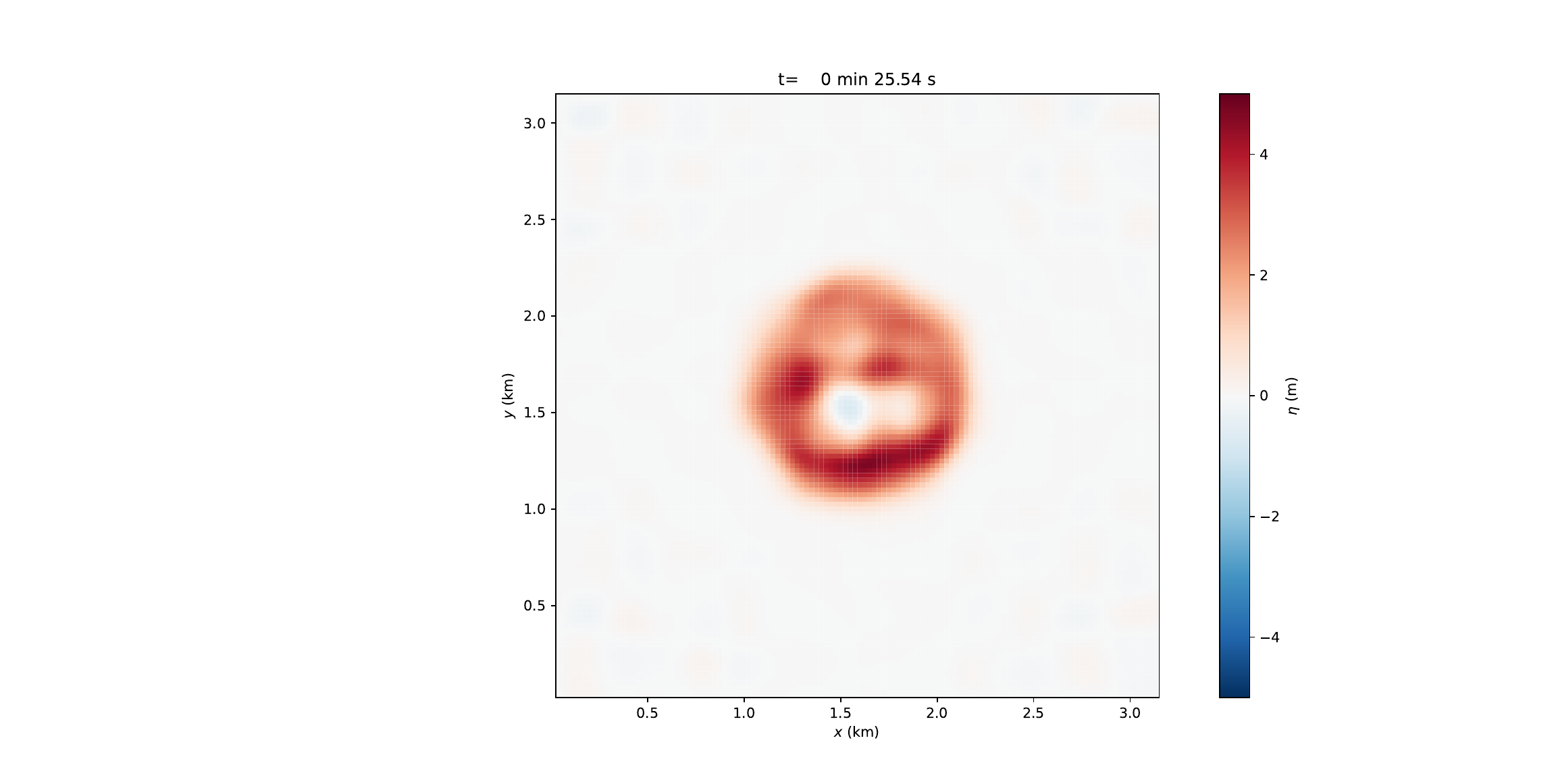}\\
        \hline
        \includegraphics[trim={10cm 0 4cm 0}, clip, width=6.5cm]{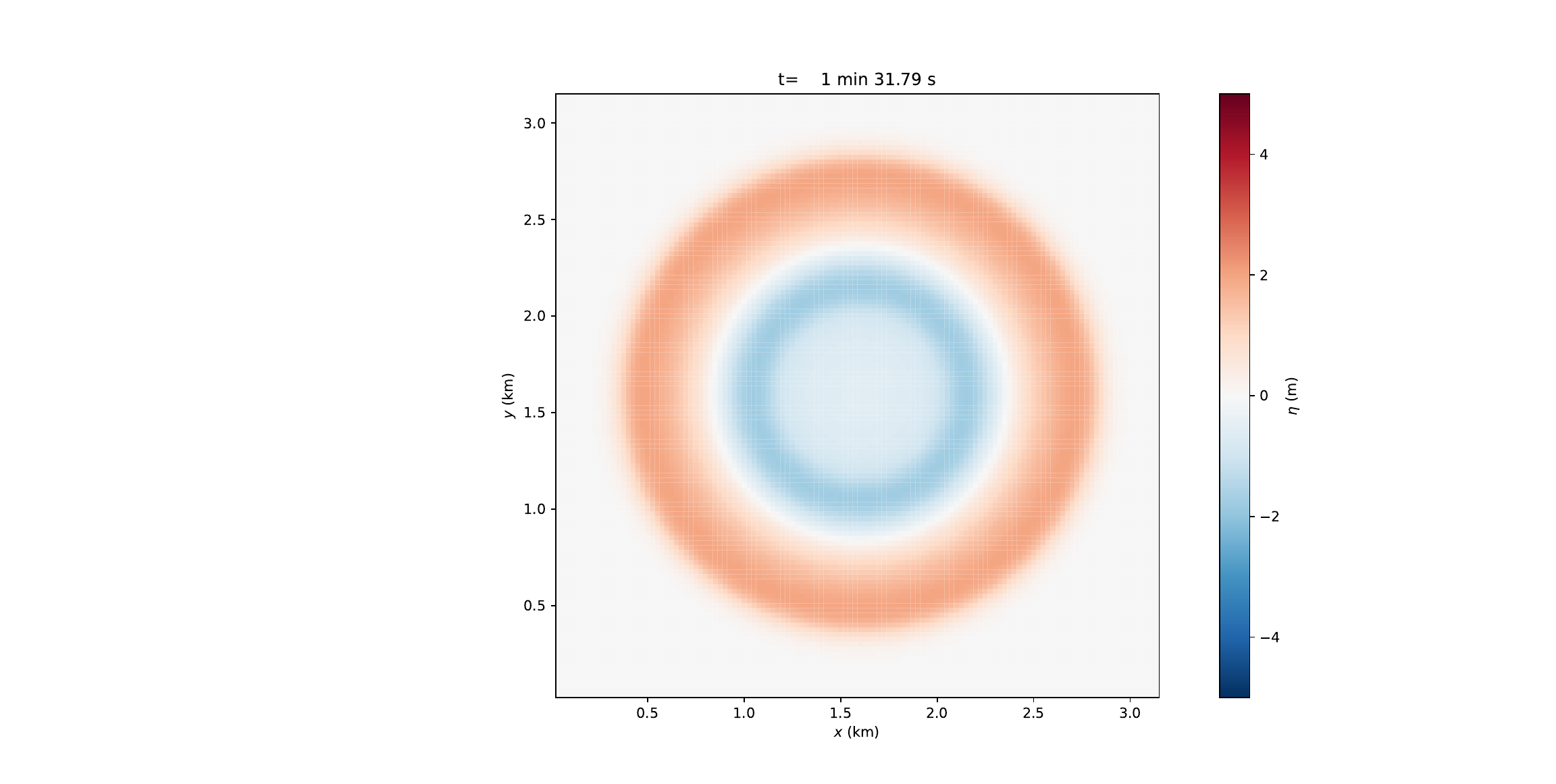} & \includegraphics[trim={10cm 0 4cm 0}, clip, width=6.5cm]{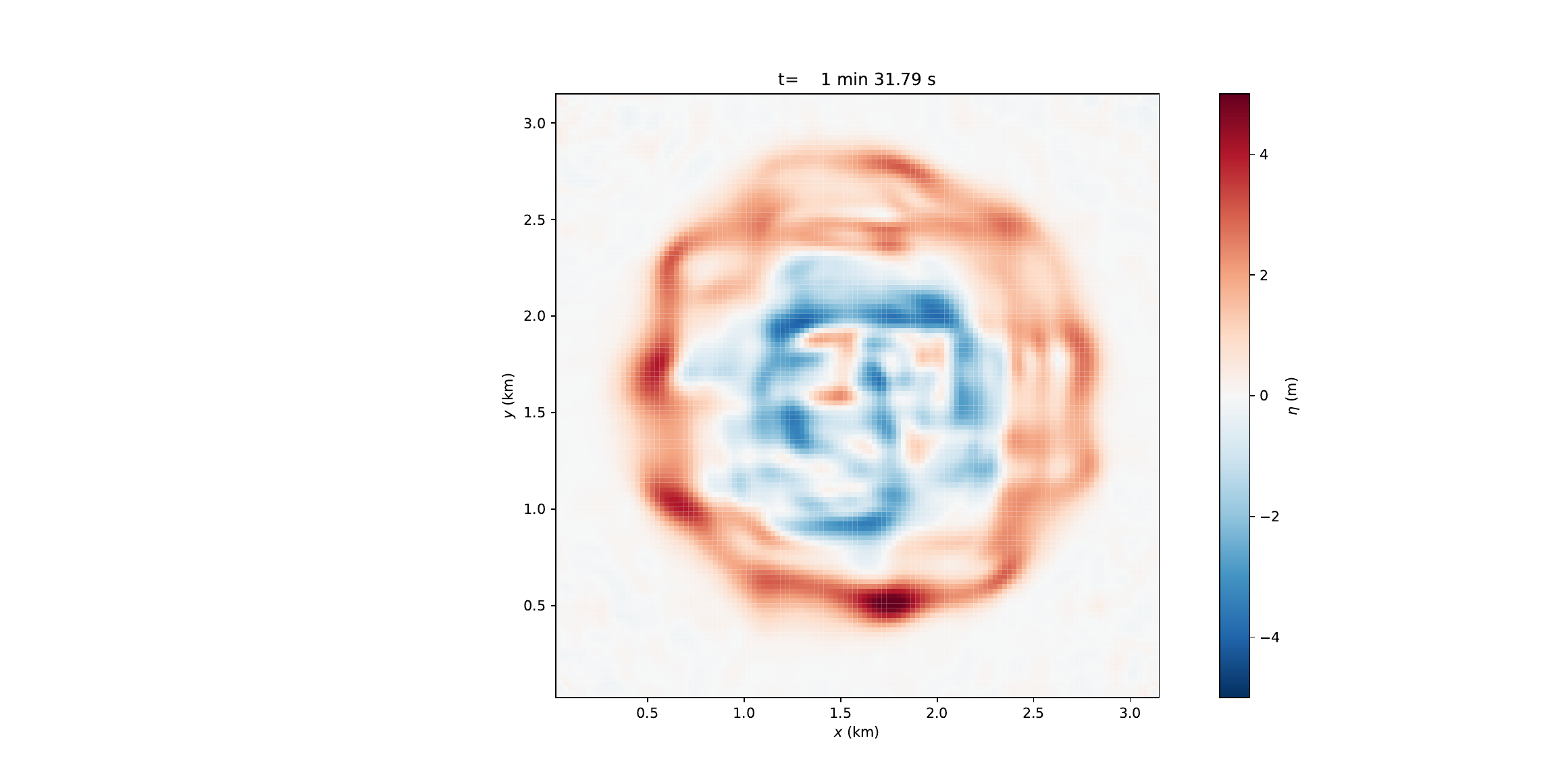}
    \end{tabular}
    \caption[Elevation of the deterministic and LU Saint-Venant systems, at fixed time.]{Elevation of the deterministic and LU Saint-Venant systems, at $t = 25.54s$ (upper left: deterministic system, upper right: LU system) and $t= 1 min \: 31.79s$ (lower left: deterministic system, lower right: LU system). The upper panels display the elevation at approximately the same time as in figure \ref{fig:plots-WW}, while the lower panels show the elevation at a ``similar advance'' of the wave.}
    \label{fig:plots-SV}
\end{figure}

In addition, we compare the deterministic Saint-Venant system to its LU counterpart. First, we observe on Figure \ref{fig:plots-SV} that the advance of the wave is much different than that of the Water Waves simulation for a similar simulated time. This can be imputed to the coarser resolution of the Water Waves system, which affects dramatically the characteristic time of evolution. This justifies the presence of the lower panels, as they show the elevation for a similar advance of the wave. In both cases, the presence of the noise tends to break the symmetry of the wave and to induce a slight deformation, which is coherent with the noise being rotational. 

\subsubsection{Statistical comparison}

\noindent
As mentioned above, the presence of the noise tends to break the symmetry of the system. However, we show that the system remains symmetric statistically, which is in the line of what is presented in \cite{DMM_2024}. As a matter of fact, we computed the empirical mean and variance over 100 realisations of the LU Saint-Venant system to compare it with the deterministic solution.
\begin{figure}[H]
    \centering
    \includegraphics[trim={4cm 0 4cm 0}, clip, width=12cm]{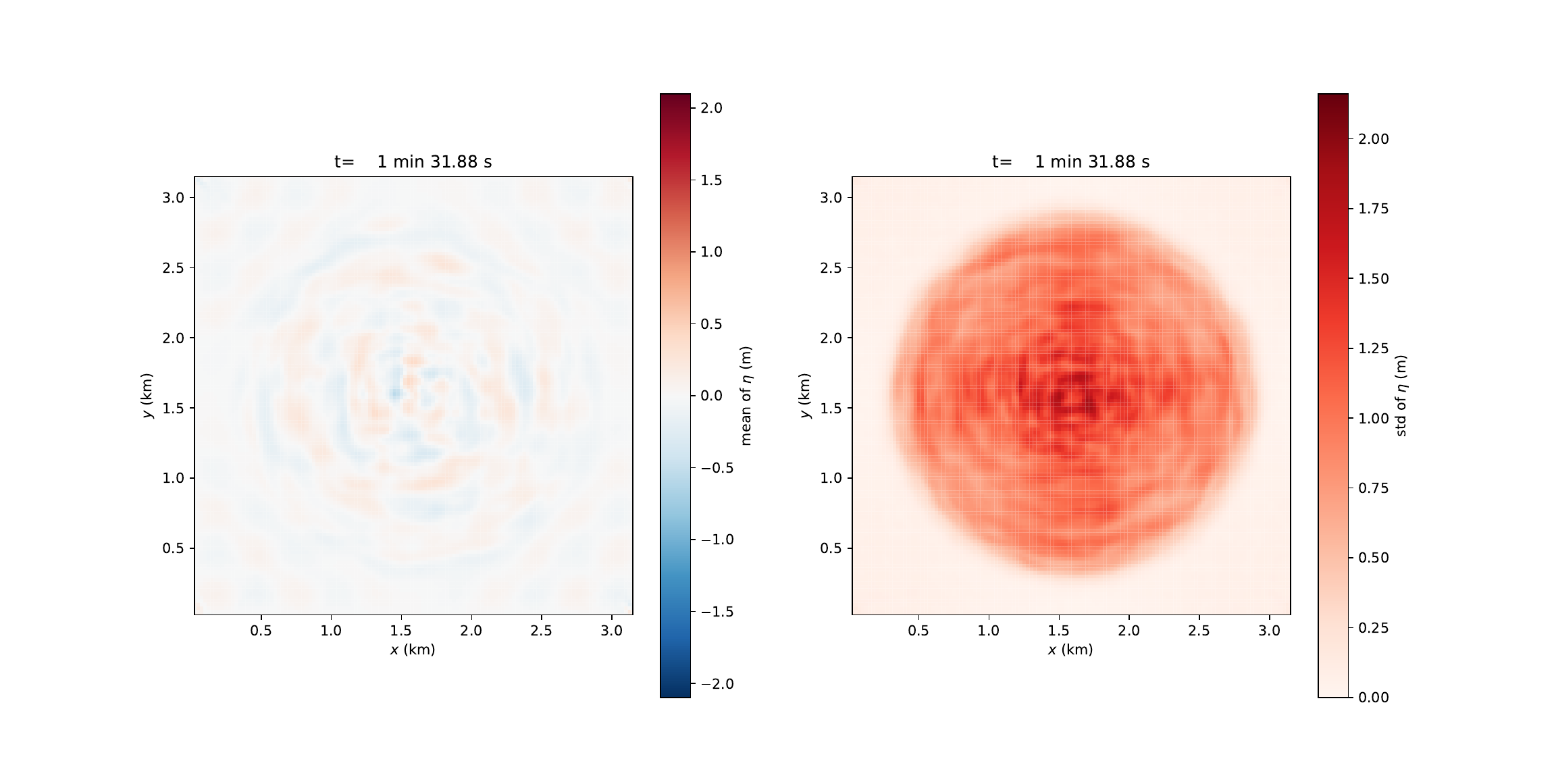}
    \caption[Statistics of the LU Saint-Venant system, at fixed time.]{Difference between the mean and the deterministic solution (left), and variance (right) for the LU Saint-Venant system, at $t = 1min 31.88s$. These empirical statistical parameters are computed over 100 realisations. We remind that it initial condition is set to be a heap of water of 10m. As the left profile is roughly contained in [-0.5,0.5], it can be considered negligible compared to the initial condition.}
\end{figure}

\subsection{Influence of an additive forcing}

\noindent
In this subsection, we analyse the Saint-Venant equation in a slightly different context. The initial condition is now set to zero, but the equation on average velocity $\bar{u}$ is forced by an additive noise. Hence, the deterministic Saint-Venant equation with additive noise forcing reads as follows,
\begin{gather}
    d_t \eta + \frac{1}{\nu} \nabla_H \bcdot (h \bar{u}) dt = 0,\\
    d_t \bar{u} + \frac{\epsilon}{\nu} (\bar{u} \bcdot \nabla_H) \bar{u} \: dt + \nabla_H \eta \: dt = \sigma_{Add} dW_t.
\end{gather}
We consider the following LU counterpart of the aforementioned problem (i.e. LU Saint-Venant system with additive noise forcing). For simplicity, we disregard the contribution of the It\={o}-Stokes drift terms, so that the system reads as follows,
\begin{equation}
\begin{cases}
    d_t \eta &+ \frac{1}{\nu} \nabla_H \bcdot (h \bar{u}) dt + \frac{1}{\nu}\nabla_H \bcdot \Big(h \circ (\frac{\epsilon_\sigma}{\epsilon} \overline{\nabla_H \sigma^{\phi} dW_t^1} + \frac{\epsilon_\sigma}{\nu} \overline{\nabla_H^\perp \sigma^{\chi} dW_t^1}) \Big) = 0,\\
    d_t \bar{u} &+ \frac{\epsilon}{\nu} (\bar{u} \bcdot \nabla_H) \bar{u} \: dt + \nabla_H \eta \: dt \\
    &+ \frac{\epsilon_\sigma}{\nu} \Big(\nabla_H \sigma^\Phi dW_t^1 \circ \nabla_H + \nabla_H^\perp (\sigma^{\chi} dW_t^1)\vert_{z=\eta} \circ \nabla_H \Big) \bar{u} = \sigma_{Add} dW_t^2, 
\end{cases}
\label{eq:forced-LUSV}
\end{equation}
where $W^i$, $i \in \{1,2\}$, are independent cylindrical Wiener processes. The LU noise has the same shape as in the noise is chosen as follows,
\begin{align*}
    \sigma dW_t^1 = \sqrt{\Upsilon_{LU}} f_{\alpha}^{k} dW_t^1, \quad \sigma_{Add} dW_t^2 = \sqrt{\Upsilon_{Add}} f_{\alpha}^{k} dW_t^2.
\end{align*}
with $\Upsilon_{LU} = 5.0 \: 10^1$ and $\Upsilon_{Add} = 1.0 \: 10^5$. Again, the wavenumber and the coefficient $\alpha$ are chosen to be $\omega = 2\pi$ and $\alpha = 10^7$.

\begin{figure}[h!]
    \centering
    \begin{tabular}{c|c}
        \includegraphics[trim={4cm 1cm 24cm 1cm}, clip, width=6.5cm]{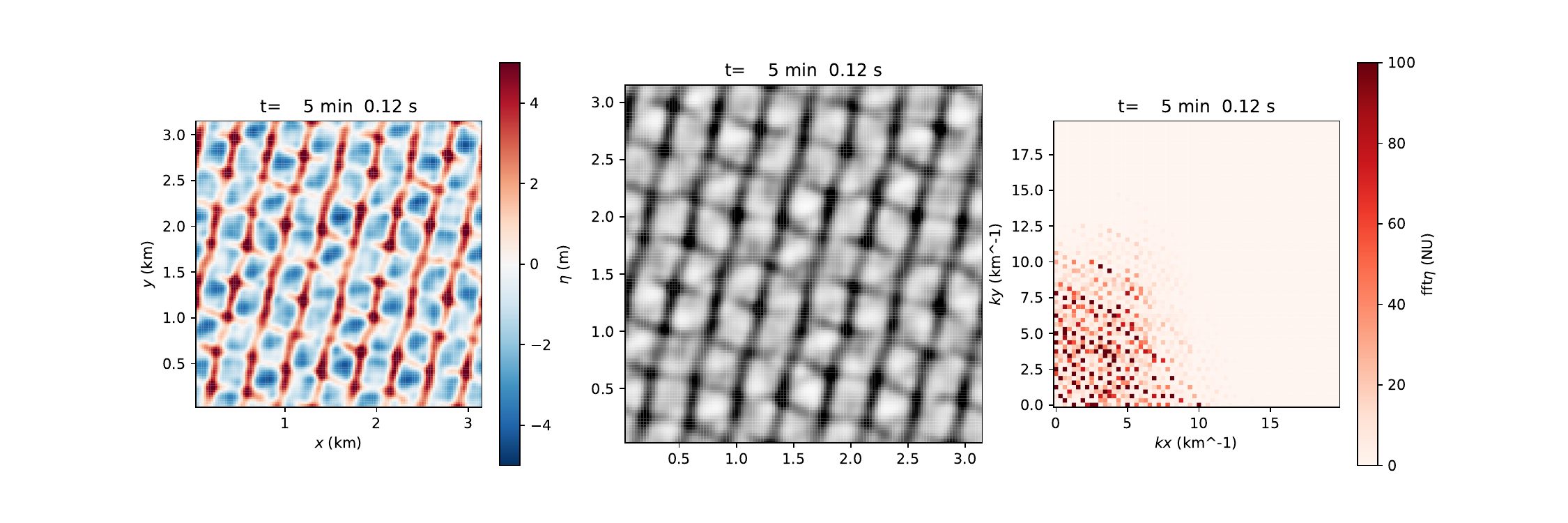} & \includegraphics[trim={4cm 1cm 24cm 1cm}, clip, width=6.5cm]{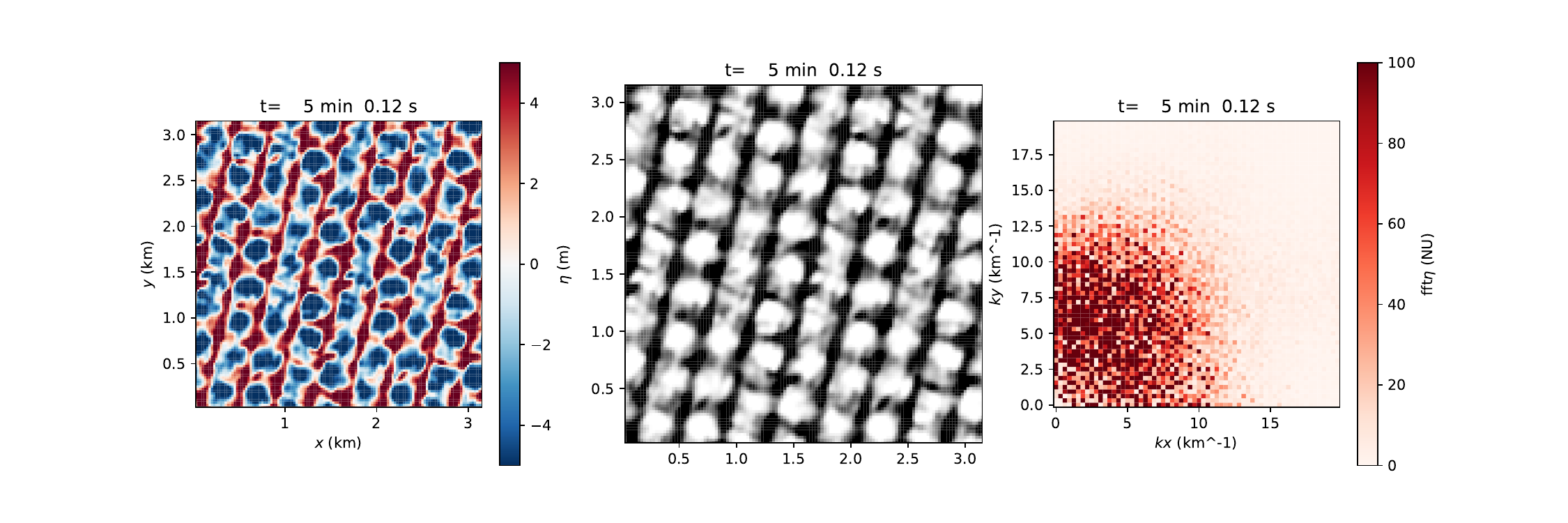}
    \end{tabular}
    \caption[Comparison between Saint-Venant system and its LU counterpart under additive noise forcing.]{Trajectories of the Saint-Venant (left) and LU Saint-Venant (right) systems, both forced with additive noise, at $t= 5 min \: 0.12s$.}
    \label{fig:comparison-detLUSV-addNoise}
\end{figure}

\begin{figure}[h!]
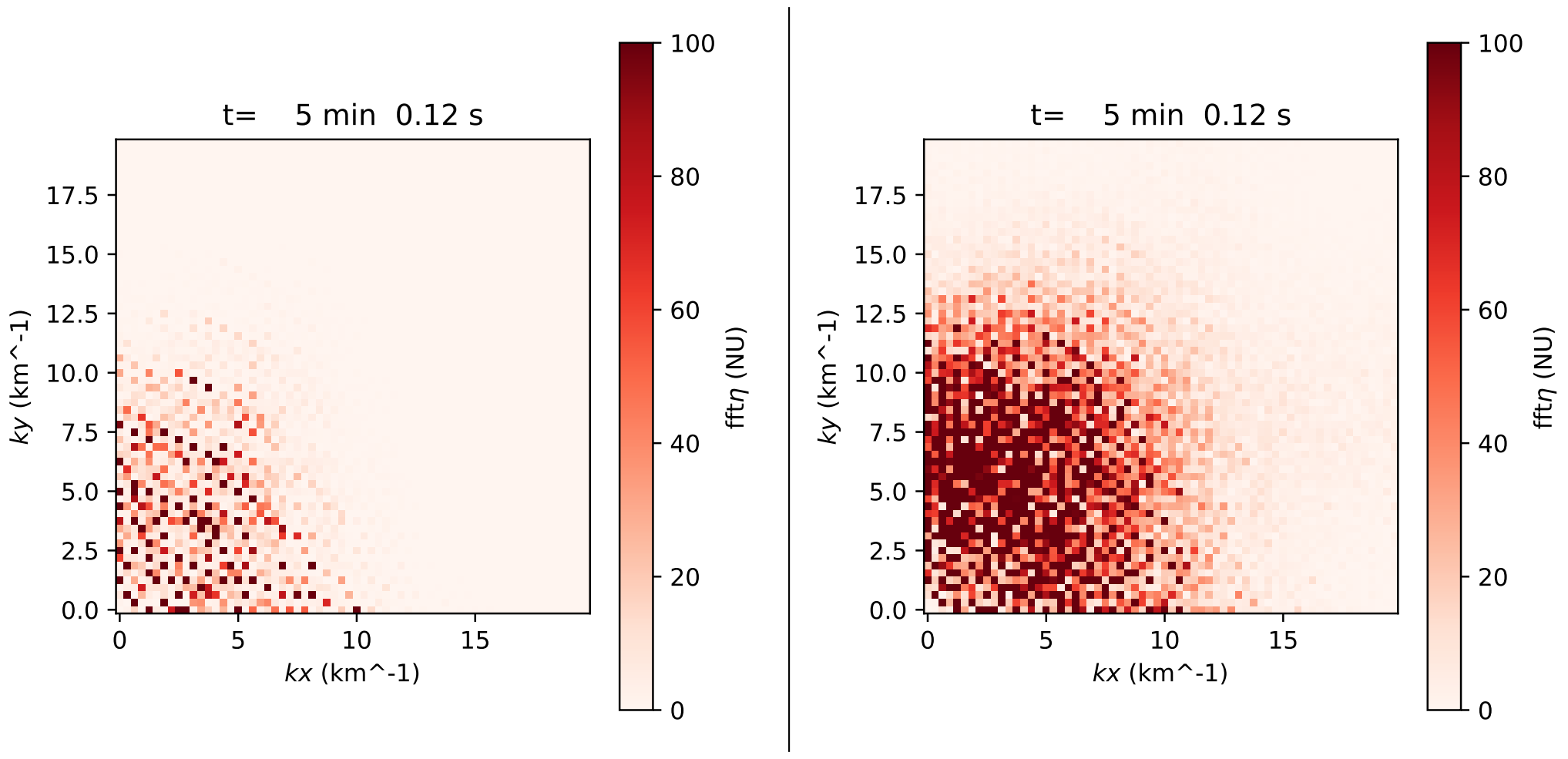

    \centering
    \begin{tabular}{c|c}
        \includegraphics[trim={24cm 1cm 3cm 1cm}, clip, width=6.5cm]{Chapitre-waves-WW/figures-additive/376_SV_noLU.pdf} & \includegraphics[trim={24cm 1cm 3cm 1cm}, clip, width=6.5cm]{Chapitre-waves-WW/figures-additive/376_SV-LU.pdf}
    \end{tabular}
    \caption[Comparison of the spectral behaviour of the deterministic and LU Saint-Venant systems]{Comparison of the spectral behaviour of the deterministic and LU Saint-Venant systems, both forced with additive noise. These 2D spectra are those associated to the plots of Figure \ref{fig:comparison-detLUSV-addNoise}.} \label{fig:spectral-behaviour}
\end{figure}

Regarding figure \ref{fig:comparison-detLUSV-addNoise}, which compares the two aforementioned systems, we observe that the elevation of the forced Saint-Venant displays a wave behaviour, and that the amplitude of these waves is increased in the presence of the transport noise arising from the LU formalism. In addition, we observe that small scale waves -- i.e. waves of small wavelength -- are induced by this transport noise. This is in the line of Figure \ref{fig:spectral-behaviour}, which compares the Fourier transforms of the realisations above -- see Figure \ref{fig:comparison-detLUSV-addNoise}. It shows that the LU Saint-Venant system displays a spectrum that is more spread than the forced Saint-Venant system. Hence, more small scale phenomena, such as small scale waves, are present in the former than in the later.

\begin{figure}[h!]
    \centering
    \begin{tabular}{c|c|c}
        \includegraphics[trim={12cm 0 6cm 0}, clip, width=4cm]{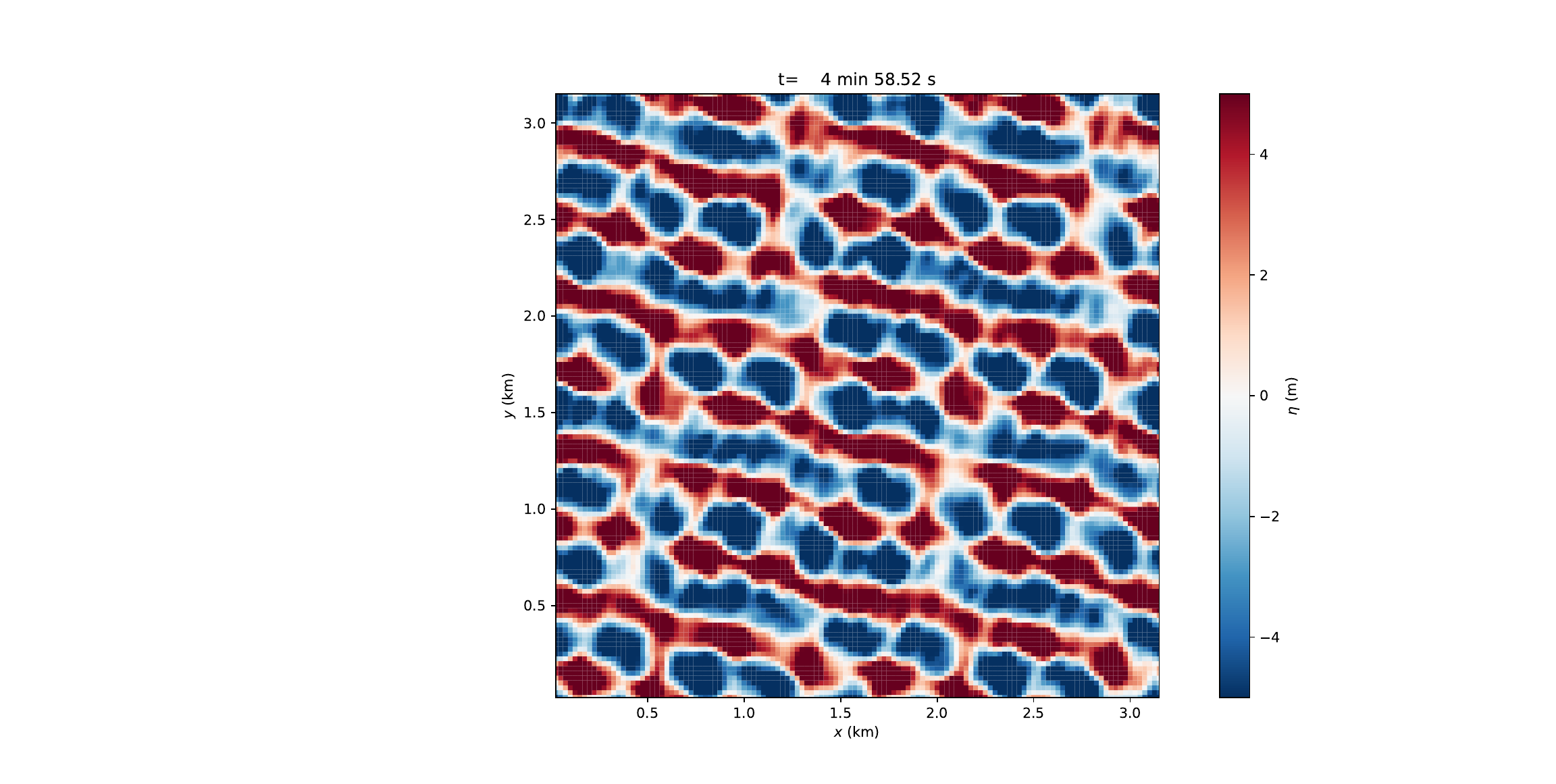} & \includegraphics[trim={12cm 0 6cm 0}, clip, width=4cm]{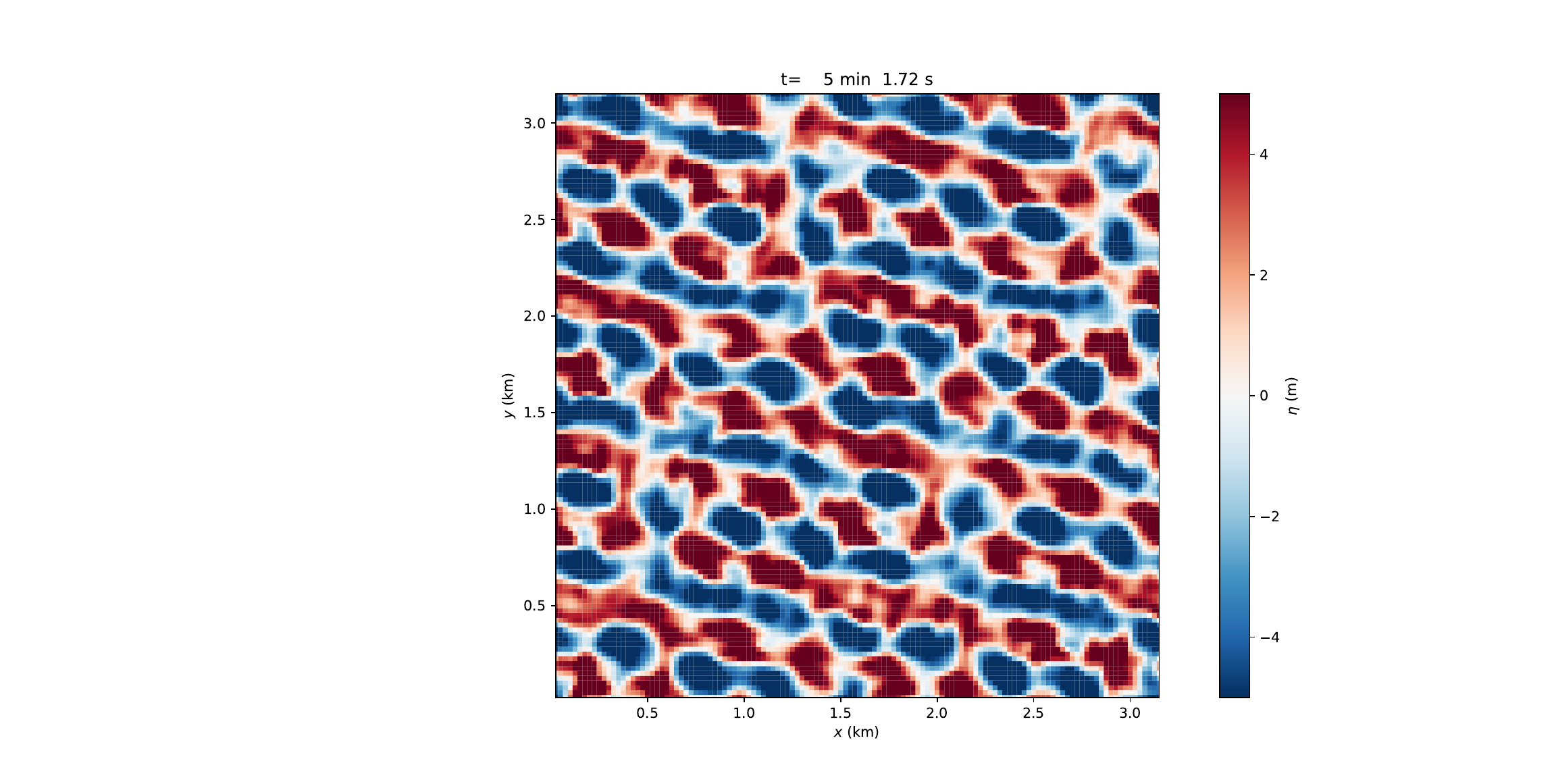} &
        \includegraphics[trim={12cm 0 6cm 0}, clip, width=4cm]{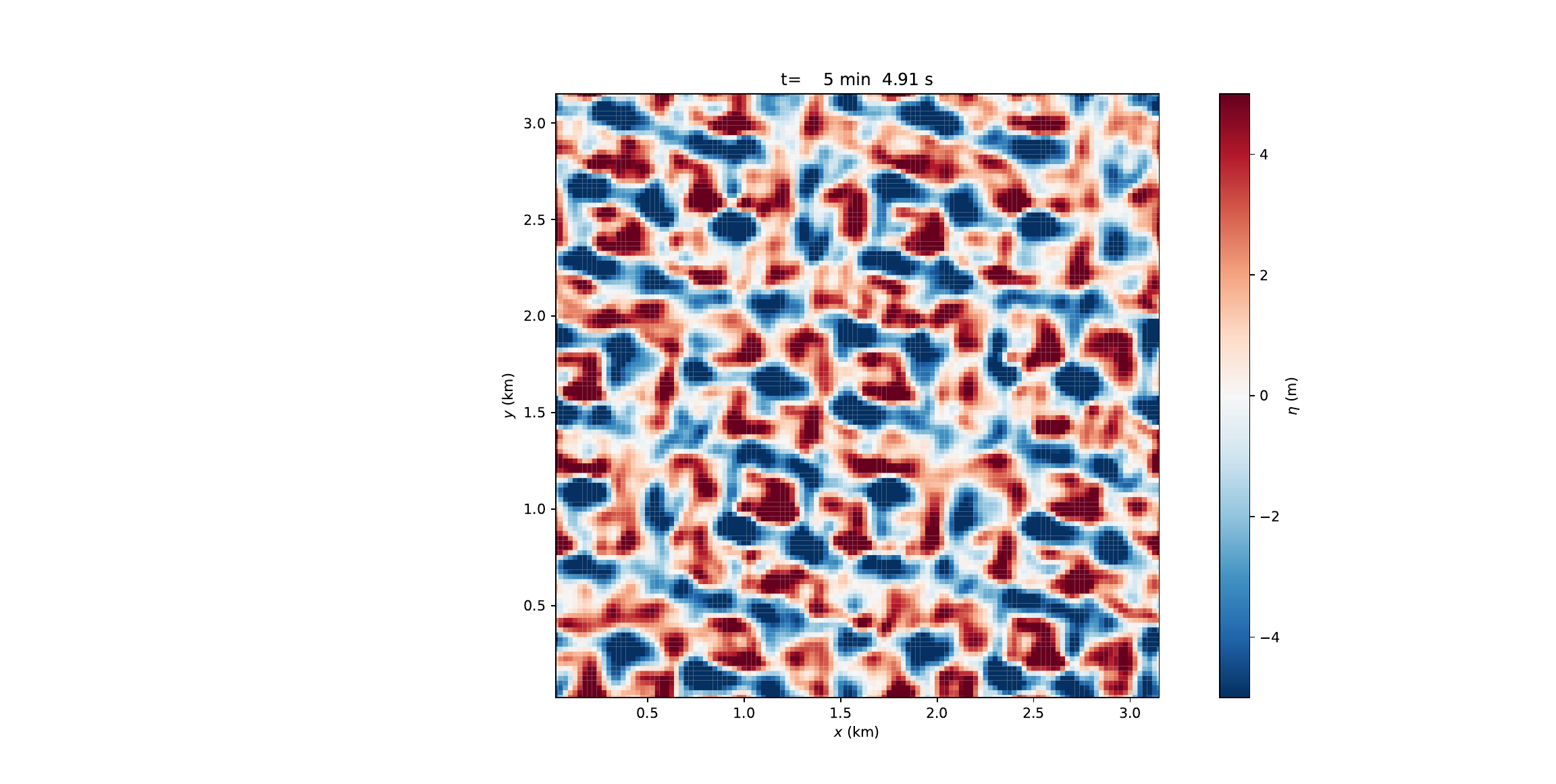} \\
        \hline
        \includegraphics[trim={12cm 0 6cm 0}, clip, width=4cm]{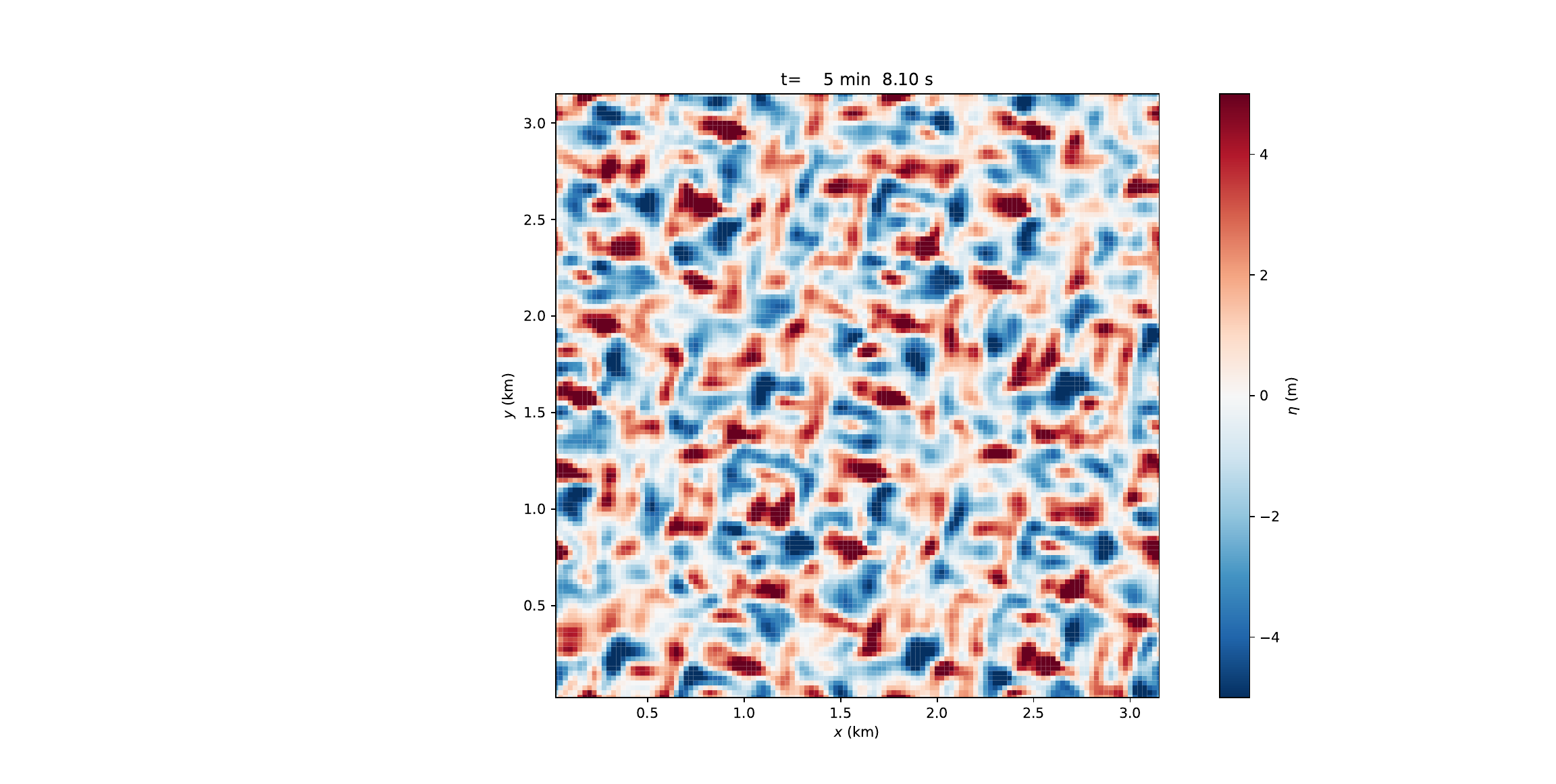} & \includegraphics[trim={12cm 0 6cm 0}, clip, width=4cm]{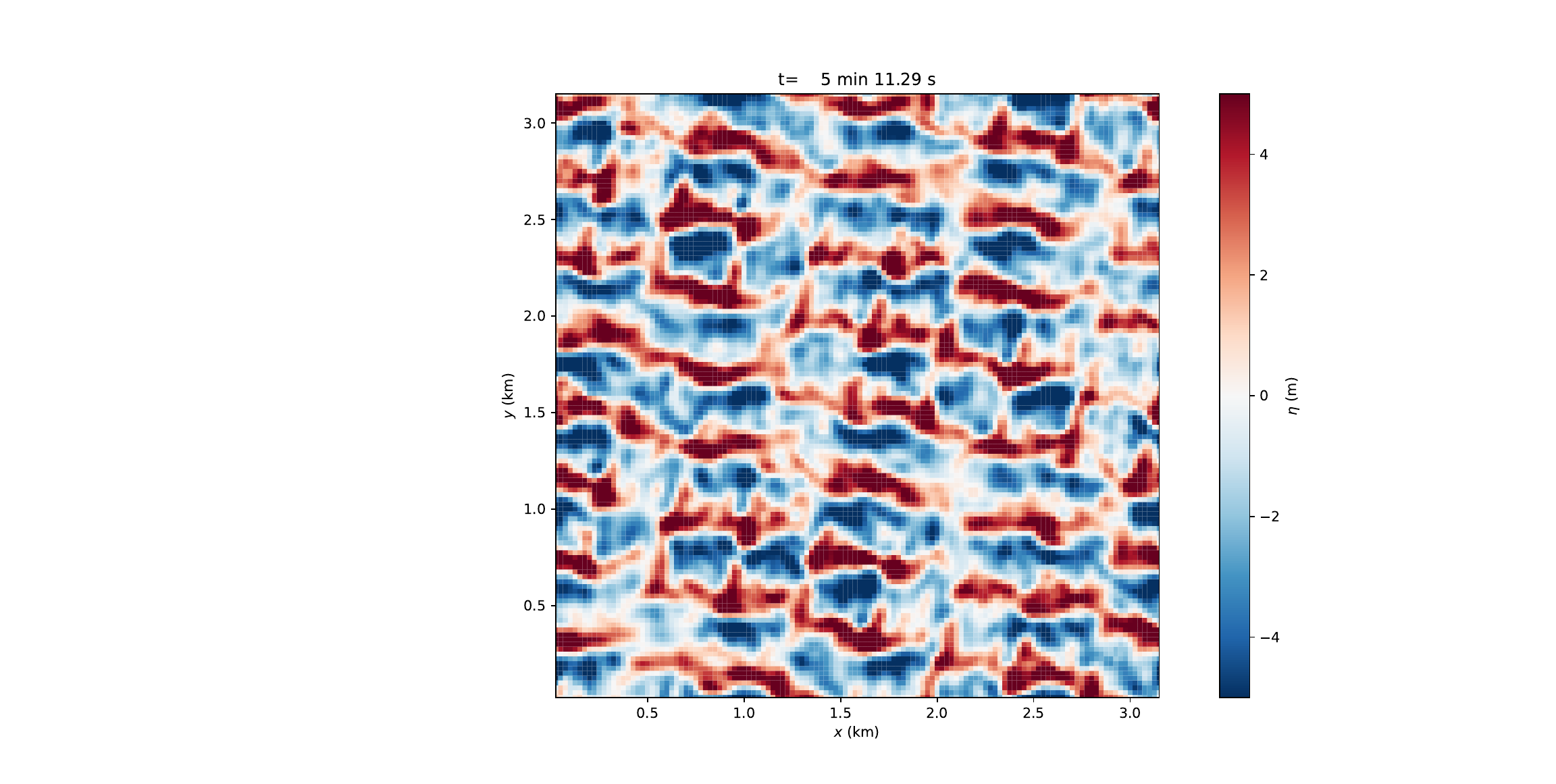} &
        \includegraphics[trim={12cm 0 6cm 0}, clip, width=4cm]{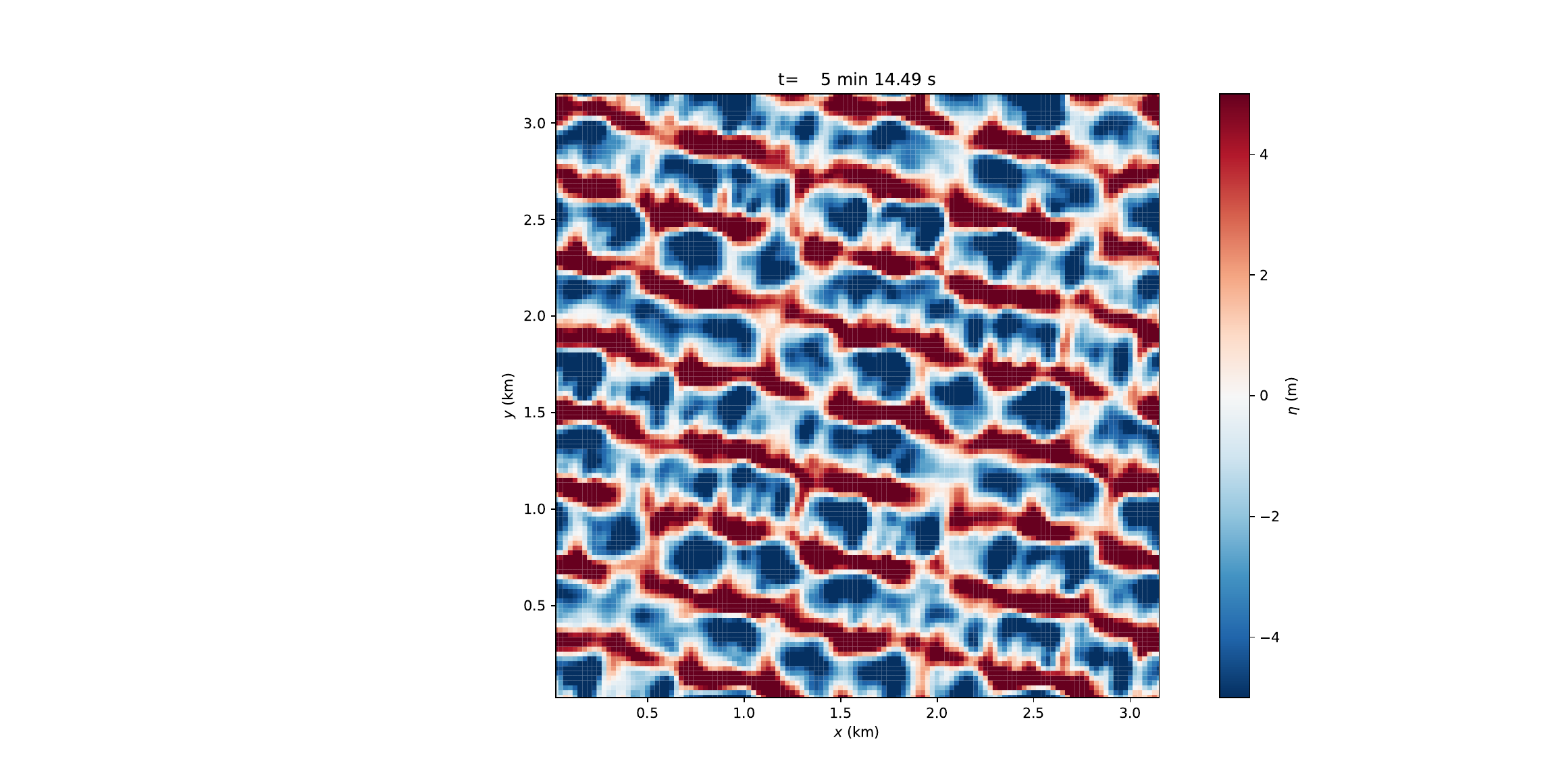}
    \end{tabular}
    \caption[Snapshots of a realisation of the LU Saint-Venant system with additive noise forcing]{Six snapshots of a realisation of the LU Saint-Venant system with additive noise forcing, taken between $t=4 min \: 12.23s$ and $t= 4 min \: 28.19s$. The time step between each snapshot is equal to $3.19s$. }
    \label{fig:fast-small-scale-waves}
\end{figure}

Furthermore, these arising small scale phenomena have a rather short time scale, compared to the rest of the system. Indeed, Figure \ref{fig:fast-small-scale-waves} displays six snapshots spanning approximately $16s$ of simulation time. These plots highlight an oscillatory behaviour that characterises the LU Saint-Venant system, which we comment below. For that purpose, we denote the beginning of the time range by $t_0 = 4 min \: 12.23s$, and by $t_k = t_0 + k \Delta t$ the observation times, where $\Delta t = 3.19s$ is the time step between observation times. The oscillatory phenomenon displayed on Figure \ref{fig:fast-small-scale-waves} is the following,
\begin{enumerate}
    \item at $t_0$, the system is ``coherent'' as an organised successions of waves, which look quasi-periodic,
    \item at $t_1$, then $t_2$, the system enters a decoherent state, where the larger scale phenomena vanish and the smaller scale ones predominate,
    \item at $t_3$, the system is fully decoherent, and consists almost only in short range waves,
    \item at $t_4$, then $t_5$, the system evolves back to an ``coherent'' state, which is similar to that of $t_0$.
\end{enumerate}
This suggests that stochasticity, as introduced in the LU formalism, can yield fast small scale behaviour phenomena such as waves, by interacting with additive random forcing terms. We interpret this behaviour as follows: since the additive and LU noises both consist in waves, the so-called ``coherent'' state of the system presented above stems from constructive interference between these noises. Then, whenever the interference becomes destructive, the system relaxes and yields small scale phenomena. Eventually, the ``coherent'' state is recovered whenever the interference is constructive again. This finding may open the door to a complementary understanding of rapid phenomena such as rogue waves, which are known to appear, move, and disappear very fast. \cite{Gemmrich2026} It is worth noticing that the aforementioned noise interaction is not directly related to the It\={o}-Stratonovitch cross-correction term between the LU noise and the additive noise, since they are chosen to be independent in our simulations.

\subsection{Capturing non-linear phenomena}

\noindent
A simple generalisation of the models presented above consists in replacing the noise $\sigma dW_t$ involved above by a correlated-in-time stochastic process. This can be achieved using Ornstein-Uhlenbeck processes. In practise, we express the noise $\sigma dW_t$ as follows from the Karhunen-Loève decomposition, where $\beta^k$ are Brownian motions,
\begin{equation}
    \sigma dW_t = \sum_k \phi_k d\beta_t^k.
\end{equation}
From this relation, we may define a new noise,
\begin{equation}
    \tilde{\sigma} dU_t = \sum_k \tilde{\phi}_k d\gamma_t^k,
\end{equation}
where $\gamma_t^k$ are Ornstein-Uhlenbeck processes, namely,
\begin{equation}
    \gamma_t^k := c_\gamma \int_0^t e^{(s-t)/\tau} d\beta_s^k,
\end{equation}
where $\tau$ is the characteristic correlation time, and $c_\gamma$ a coefficient.

The LU approach can be interpreted as way to account non-linear phenomena by representing them through the noise. To assess the capability of the Saint-Venant system to represent them, we base our work on the real data analysis presented in \cite{Gemmrich2026}. In this article, the study rests on the statistical analysis of a dataset of various real waves. Notably, the ``non-Gaussian'' character of the distribution is investigated through its skewness and kurtosis.

\subsubsection{Monochromatic Ornstein-Uhlenbeck LU noise with decorrelated additive noise}

\begin{figure}
    \centering
    \includegraphics[trim={0 0 0 0}, clip, width=11cm]{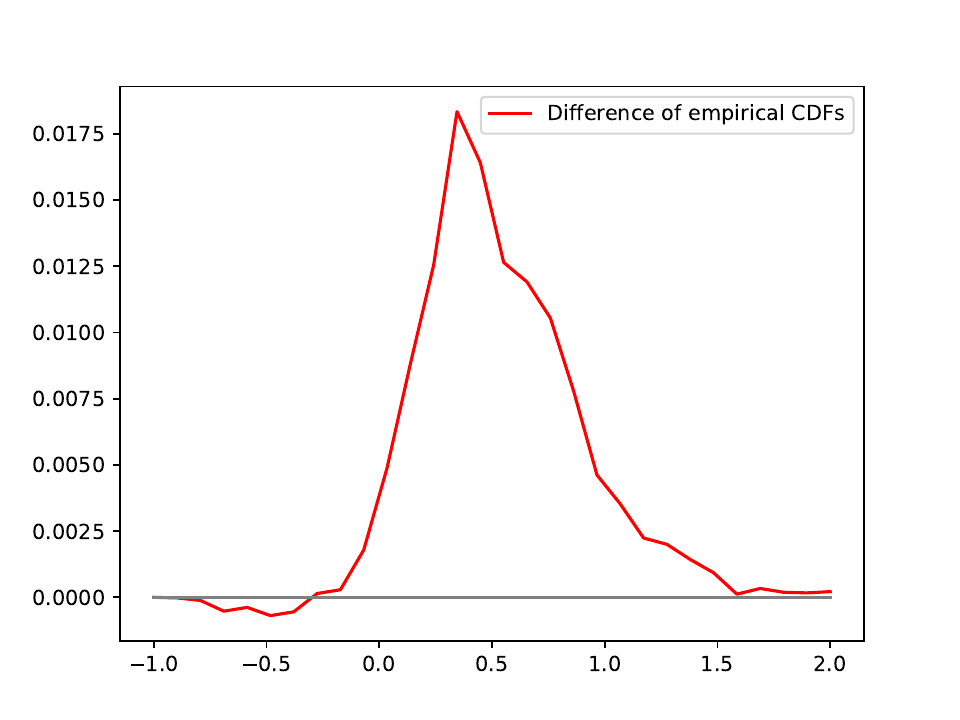}
    \caption[Representation of the empirically computed difference of repartition functions.]{Representation of the empirically computed function $F$ for Datasets 1.1 and 1.2, as defined in \eqref{eq:diffCDF}.}
    \label{fig:comparison-CDF-monochromatic}
\end{figure}

\begin{figure}
    \centering
    \includegraphics[trim={0 2cm 0 2cm}, clip, width=12cm]{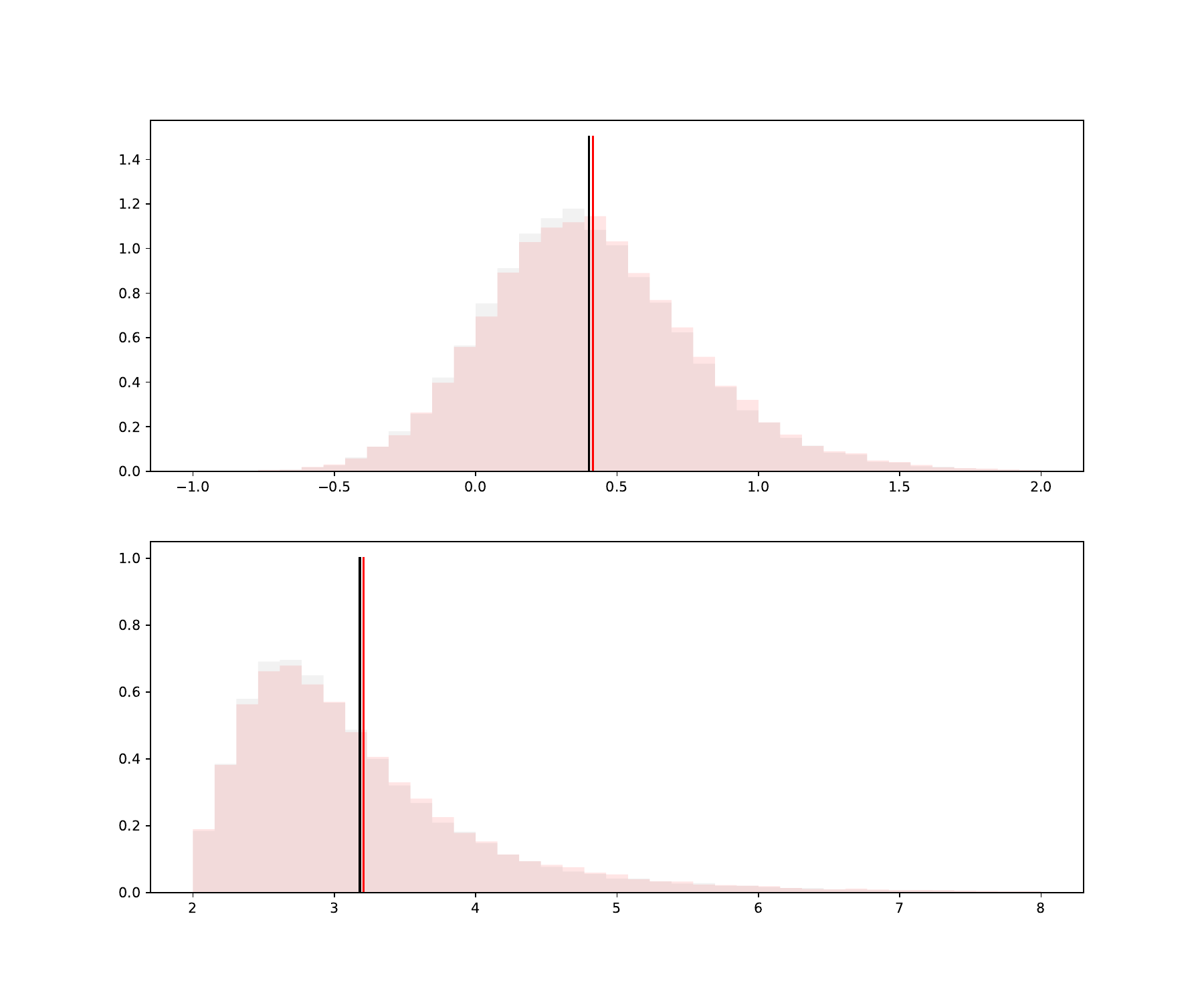}
    \caption[Distribution in skewness and kurtosis axes for two datasets.]{Distribution of skewness (upper) and kurtosis (lower) for Dataset 1.1 (red) and 1.2 (grey). The means are represented by vertical lines.}
    \label{fig:hist-skew-kurt}
\end{figure}

\noindent
To adapt the approach of \cite{Gemmrich2026}, we create two datasets of waves from the trajectories of the stochastic Saint-Venant system \eqref{eq:forced-LUSV}. For this purpose, we consider only monochromatic noises first. In addition, we use the Ornstein-Uhlenbeck type described above as a LU type noise, the additive noise being unchanged. The coefficient $\tau$ is chosen to be $-1/\log(0.9)$ and $c_\gamma = -\log(0.9)$, so that they are independent of $k$. Doing so allows to investigate the behaviour when the smaller scale dynamics is slightly correlated in time, while the forcing remains decorrelated.
\begin{itemize}
    \item[$\bullet$] Dataset 1.1: 50 trajectories are generated in the setting presented above when $t \leq 20 min$, with the following parameters: $\Upsilon_{LU} = 1.0 \: 10^7$, $\Upsilon_{Add} = 1.0 \: 10^4$, $k=2 \pi$ and $\alpha = 10^7$. Whenever $t \geq 20 min$, the additive noise $\sigma_{Add} dW_t^2$ is set to zero. Then, the dataset is formed by space-time sampling, over a (space) grid of $32 \times 32$ positions, with $41$ sample times over the interval $[20 min, 40 min]$. For each space-time points, the skewness and kurtosis are computed over the 50 realisations mentioned above. 
    \item[$\bullet$] Dataset 1.2: This dataset is generated similarly, except for the fact that $\sigma_{LU} dW_t^2$ is set to zero. It corresponds to a (deterministic) Saint-Venant system that is forced by a decorrelated additive noise. 
\end{itemize}
As derived in \cite{MK1998}, an ideal law relating the skewness $\mu_3$ to the kurtosis $\mu_4$ for wave height distribution is given by the following formula,
\begin{equation}
    \mu_4 = 3 + (\frac{4}{3} \mu_3)^2.
\end{equation}
The law above can be generalised as,
\begin{equation}
    \mu_4 = B + (A \mu_3)^2, \label{eq:ideal-law}
\end{equation}
where the coefficients $A$ and $B$ are estimated through linear regression.
\begin{figure}
    \centering
    \includegraphics[trim={0 0 0 0}, clip, width=6cm]{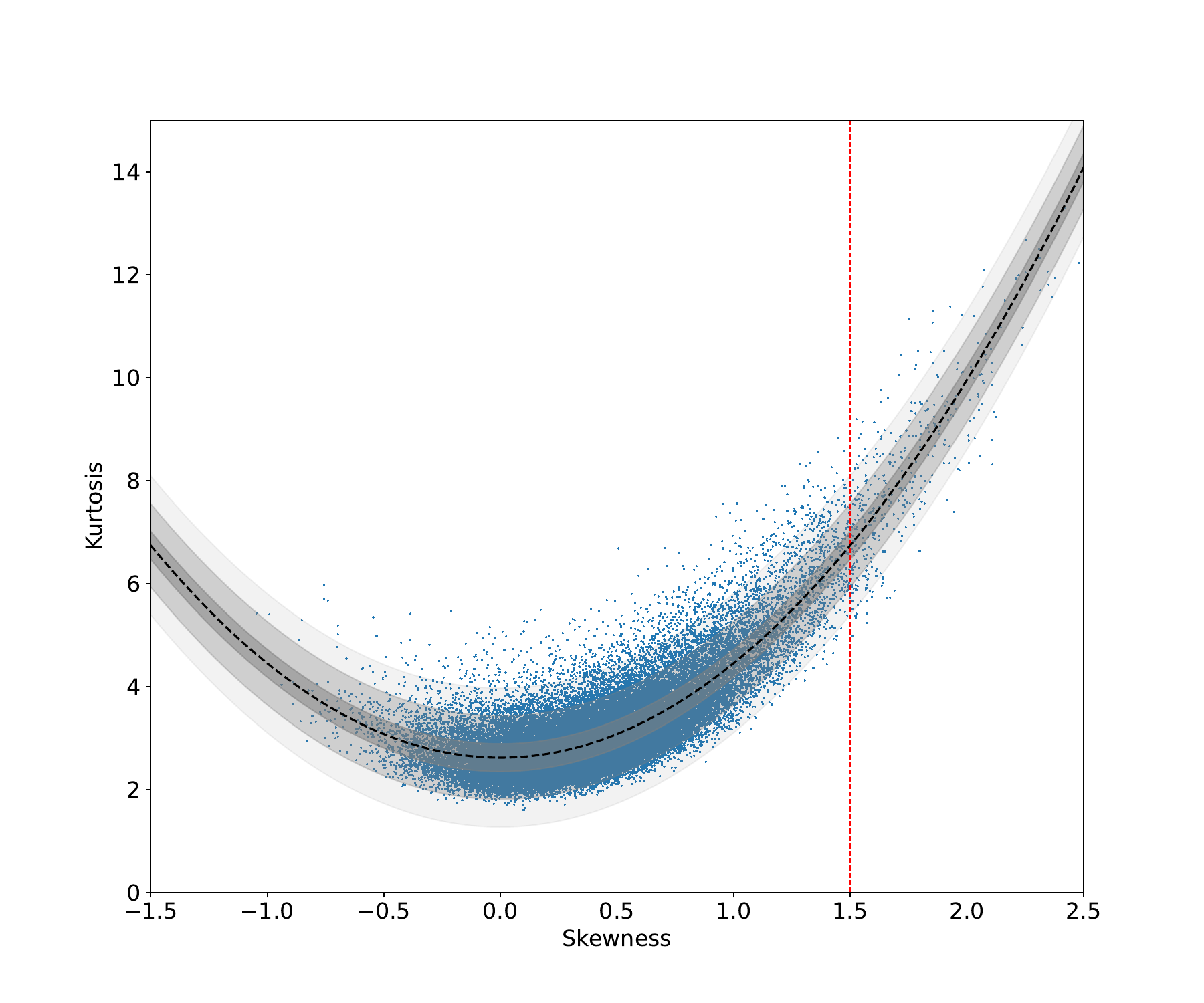}
    \includegraphics[trim={0 0 0 0}, clip, width=6cm]{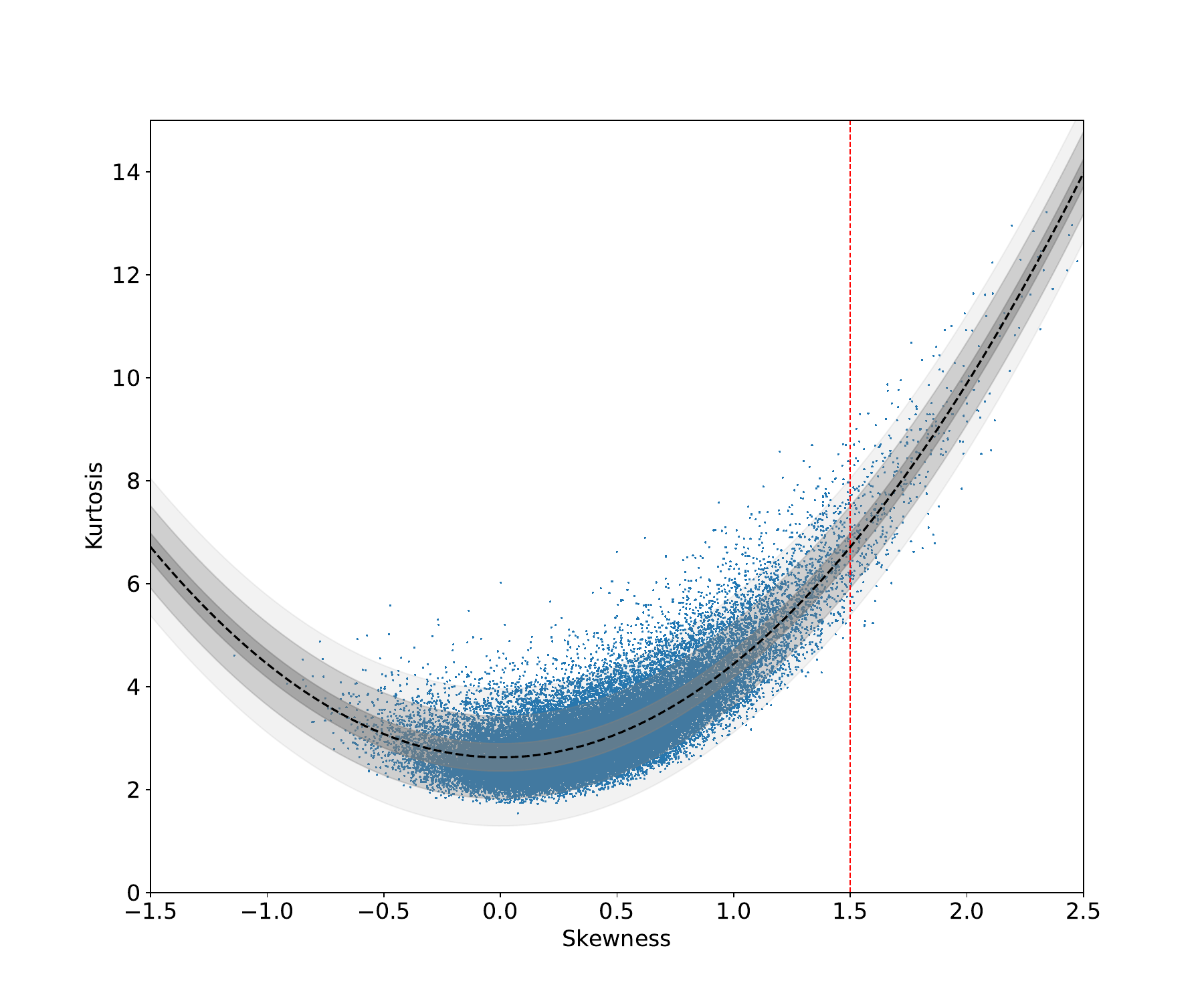}
    \caption[Scatter plots of the LU Saint-Venant system with Ornstein-Ulhenbeck noise for two datasets.]{Scatter plot of Dataset 1.1 (left) and Dataset 1.2 (right) in skewness-kurtosis axes. The dotted black line is given by the ideal law \ref{eq:ideal-law}, where the coefficients are estimated according to each dataset by linear regression. The grey zones, from darker to lighter, correspond respectively to a spread of $k \: S$ around the ideal law, $k \in \{1,3,5\}$, where $S$ is the RMSE between the dataset and this ideal law. The vertical dotted red line represents an arbitrary limit of ``non-Gaussian'' distribution with respect to the skewness $\mu_3$, which is set here to $\mu_3 = 1.5$. Total number of realisations: 41984. Number of realisations with $(\mu_3 \geq 1.5)$: 409 (Dataset 1.1) and 370 (Dataset 1.2). Parameters for Dataset 1.1: $A=1.35$, $B=2.62$. Parameters for Dataset 1.2: $A=1.35$, $B=2.63$.}
    \label{fig:scatter-plots-monochromatic}
\end{figure}
To compare the distributions with and without the LU noise, we define the cumulative distribution function of a random variable $X$ (abbreviated as CDF) as follows, for $x \in \R$,
\begin{equation}
    F_{X}(x) = \mathbbm{P}(X \leq x).
\end{equation}
Next, we denote by $F_{\mu_3}^{SV}$ and $F_{\mu_3}^{LUSV}$ the CDF of skewness associated respectively to the classical and stochastic Saint-Venant models. With this notation, the following quantity can be interpreted as a ``distance'' between these two distributions, in terms of their capability to capture events of high skewness,
\begin{equation}
    F(x) := F_{\mu_3}^{SV}(x) - F_{\mu_3}^{LUSV}(x) = \mathbbm{P}(\mu_3^{LUSV} > x) - \mathbbm{P}(\mu_3^{SV} > x). \label{eq:diffCDF}
\end{equation}
Figure \ref{fig:comparison-CDF-monochromatic} shows that $F$ is consistently positive for positive values of skewness. This indicates that the LU Saint-Venant system provides overall more realisations of higher skewness than the classical one. Consequently, LU captures more non-Gaussian -- hence non-linear -- phenomena than its deterministic counterpart. The histograms of the skewness and kurtosis support this observation, as shown in Figure \ref{fig:hist-skew-kurt}. As a matter of fact, the mean skewness and kurtosis are (slightly) greater for Dataset 1.1 than for Dataset 1.2. Both histograms also show that the former yields consistently greater values of skewness and kurtosis than the latter.

Yet, as mentioned in the caption of Figure \ref{fig:scatter-plots-monochromatic}, the values of $A$ and $B$ are similar for both datasets. These values seem relatively coherent to those provided in \cite{Gemmrich2026}, since $A \sim 4/3$ and $B$ belongs to the interval $[2.5,3.0]$. Although the differences between these two datasets are small overall, the results are encouraging, especially since the noise is chosen ``arbitrarily''. We expect that further experiments with suitably calibrated noise yield even more non-linear phenomena captured by the LU system. 

\subsubsection{Monochromatic decorrelated LU noise with Ornstein-Uhlenbeck additive noise}

\begin{figure}
    \centering
    \includegraphics[trim={0 0 0 0}, clip, width=6cm]{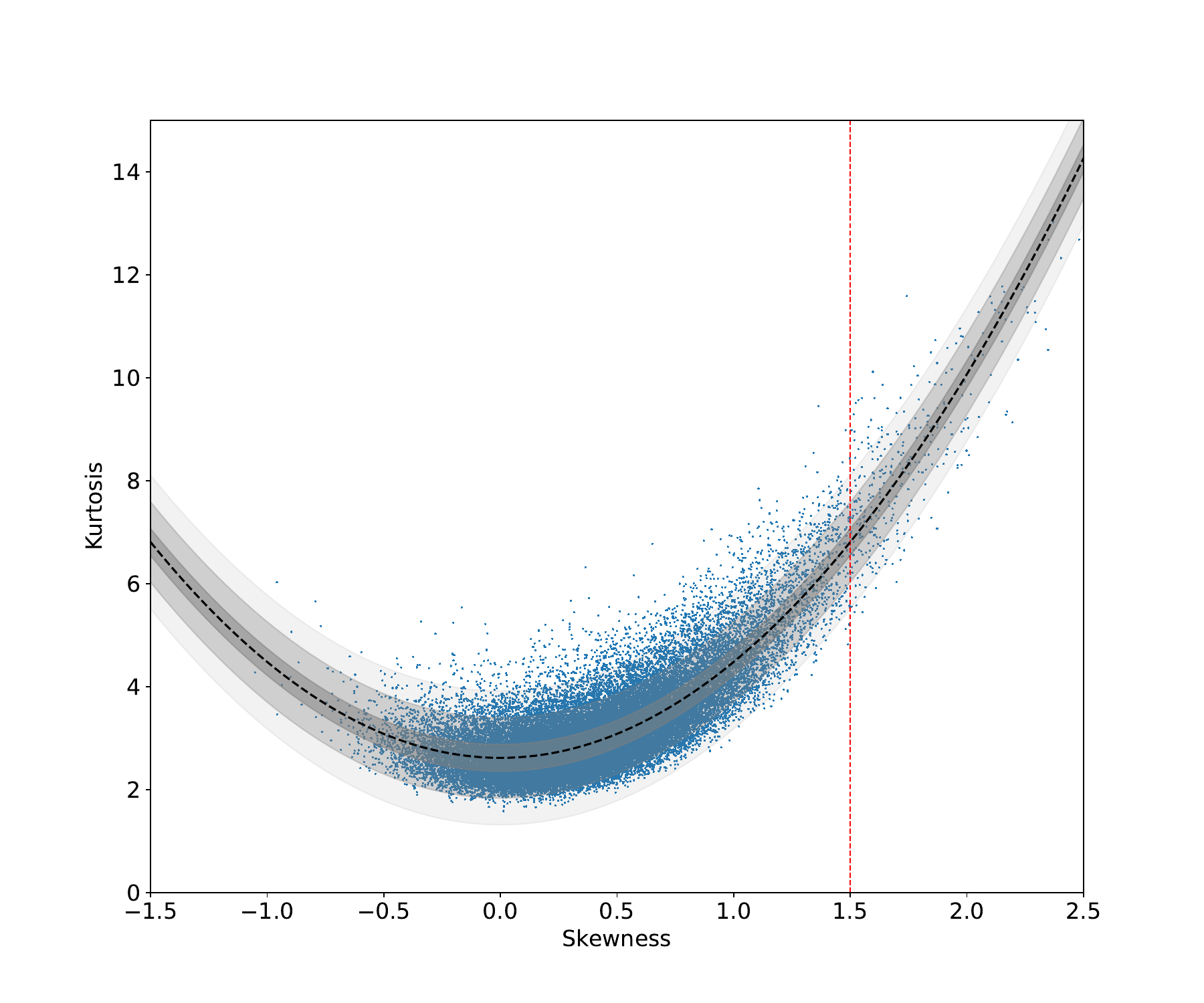}
    \includegraphics[trim={0 0 0 0}, clip, width=6cm]{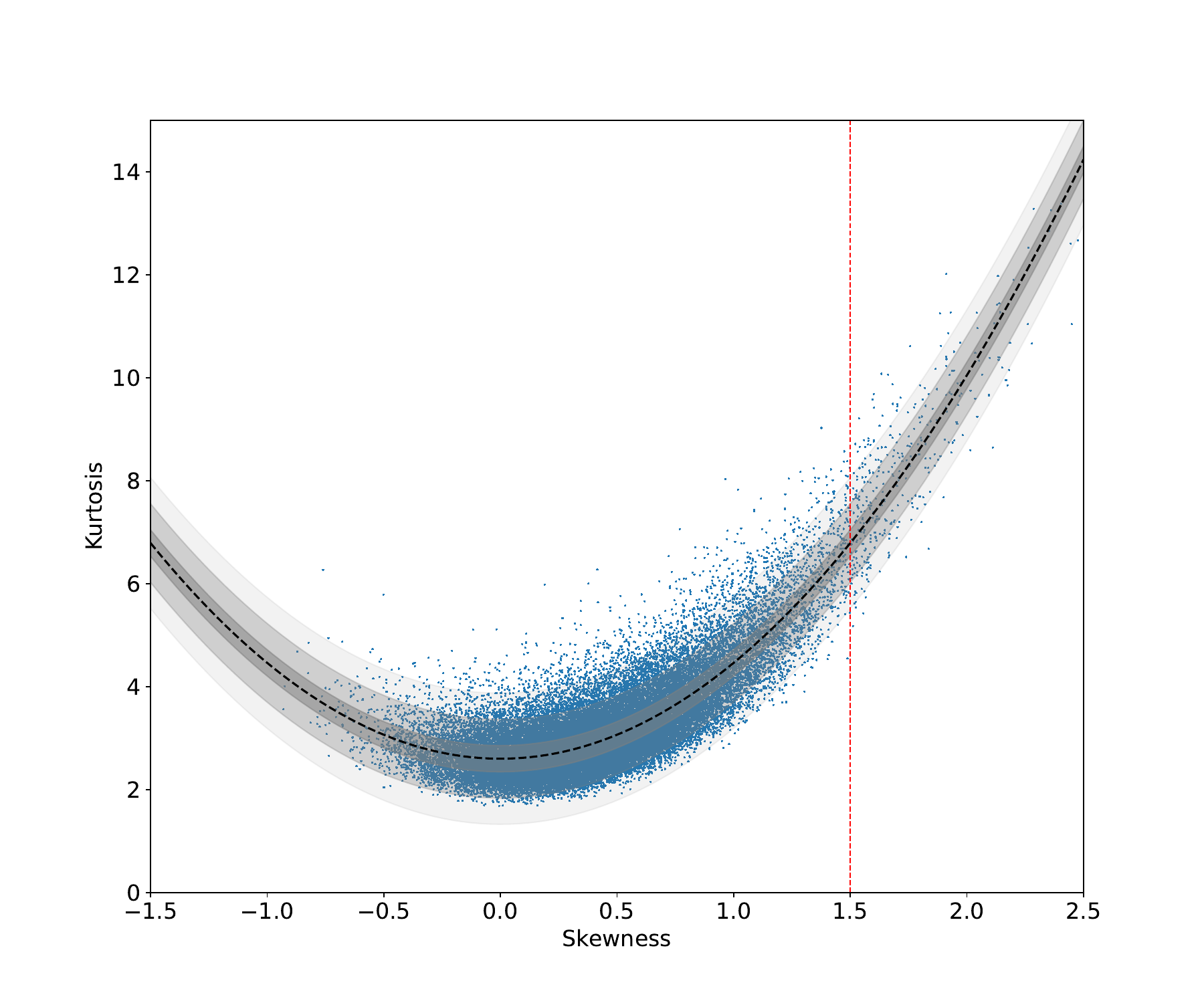}
    \caption[Scatter plots of the LU Saint-Venant system with Ornstein-Ulhenbeck noise for two datasets.]{Scatter plot of Dataset 2.1 (left) and Dataset 2.2 (right) in skewness-kurtosis axes, similar to Figure \ref{fig:scatter-plots-monochromatic}. Total number of realisations: 41984. Number of realisations with $(\mu_3 \geq 1.5)$: 349 (Dataset 2.1) and 357 (Dataset 2.2). Parameters for Dataset 2.1: $A=1.37$, $B=2.62$. Parameters for Dataset 2.2: $A=1.36$, $B=2.60$.}
    \label{fig:scatter-plots-monochromatic-corrAdd}
\end{figure}

\begin{figure}
    \centering
    \includegraphics[trim={0 0 0 0}, clip, width=11cm]{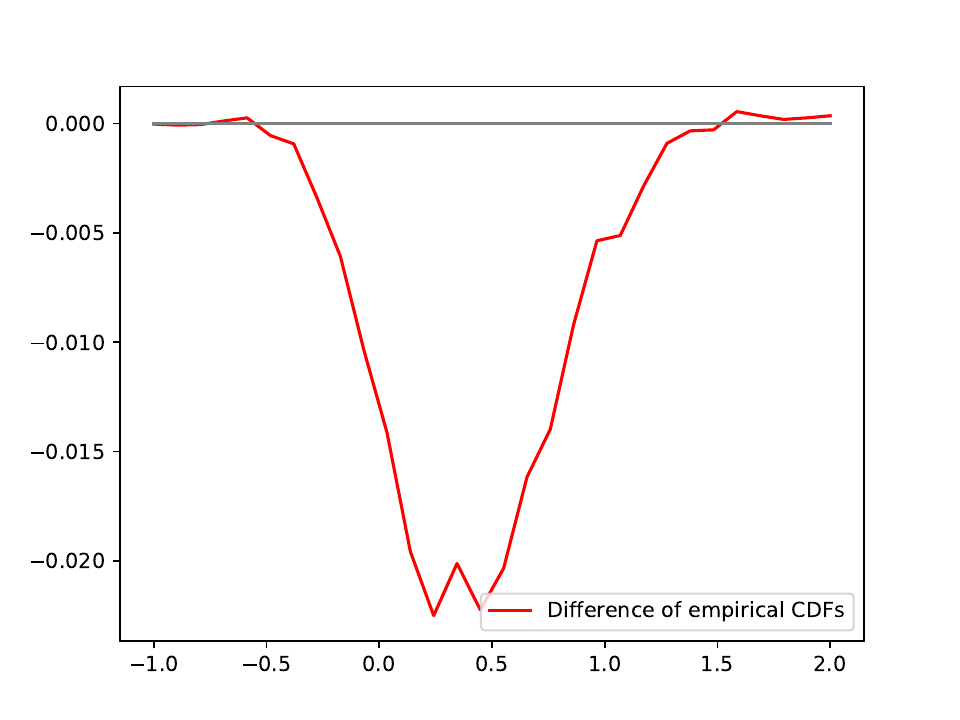}
    \caption[Representation of the empirically computed difference of repartition functions.]{Representation of the empirically computed function $F$ for Datasets 2.1 and 2.2, as defined in \eqref{eq:diffCDF}.}
    \label{fig:comparison-CDF-monochromatic-corrAdd}
\end{figure}

\begin{figure}
    \centering
    \includegraphics[trim={0 2cm 0 2cm}, clip, width=12cm]{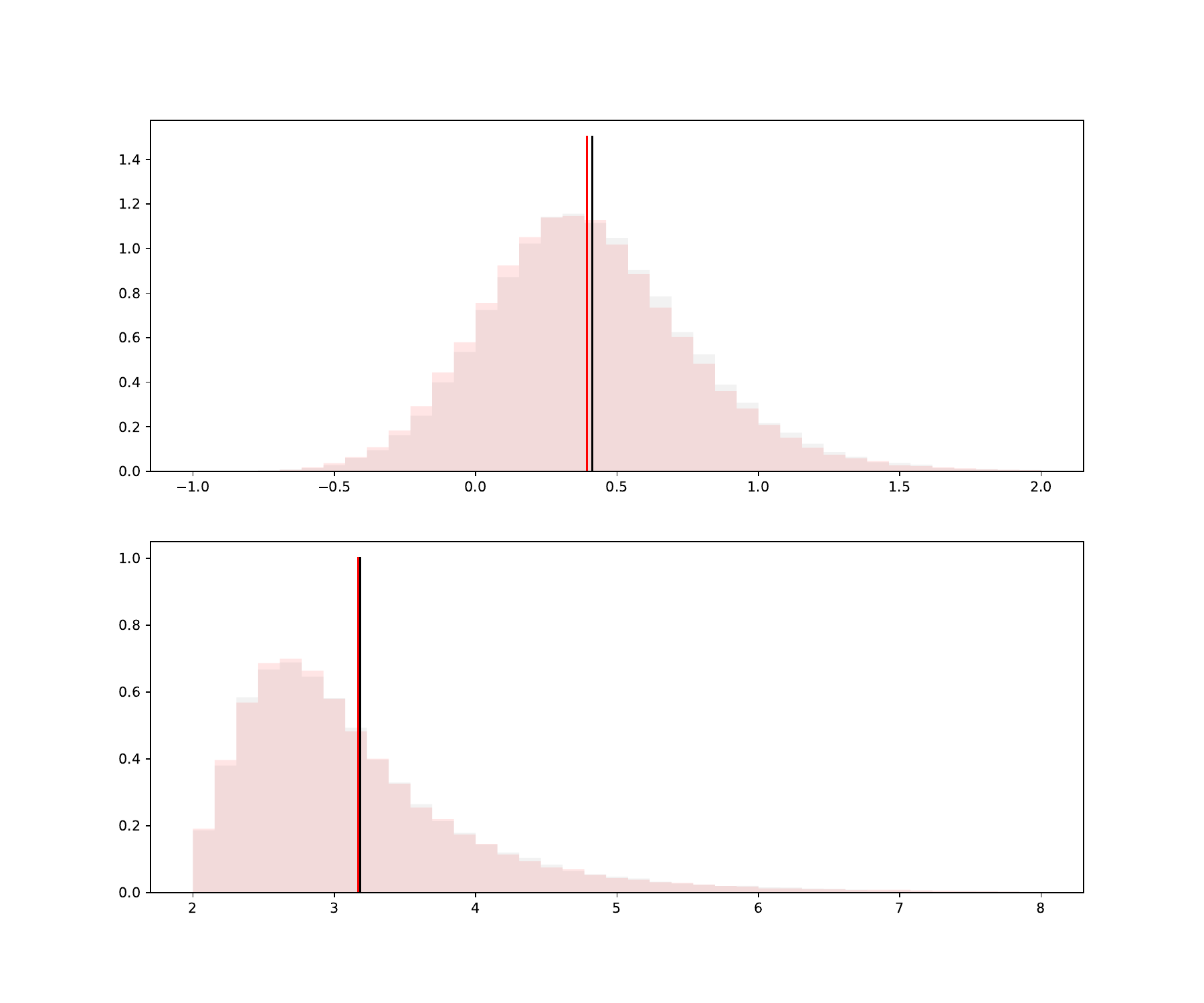}
    \caption[Distribution in skewness and kurtosis axes for two datasets.]{Distribution of skewness (upper) and kurtosis (lower) for Dataset 2.1 (red) and 2.2 (grey). The means are represented by vertical lines.}
    \label{fig:hist-skew-kurt-corrAdd}
\end{figure}

\noindent
We are now interested in the case of LU noise with Ornstein-Uhlenbeck additive noise, both being monochromatic. The coefficient $\tau$ associated to the additive noise is chosen to be $-1/\log(0.9)$. Similarly, we generate two datasets according to our new setting.
\begin{itemize}
    \item[$\bullet$] Dataset 2.1: 50 trajectories are generated in the setting presented above, with the following parameters: $\Upsilon_{LU} = 1.0 \: 10^7$, $\Upsilon_{Add} = 5.0 \: 10^5$, $k=2 \pi$ and $\alpha = 10^7$.
    \item[$\bullet$] Dataset 2.2: This dataset is generated similarly, except for the fact that $\sigma_{LU} dW_t^2$ is set to zero. It corresponds to a (deterministic) Saint-Venant system that is forced by an Ornstein-Uhlenbeck additive noise.
\end{itemize}
These datasets give similar kurtosis-skewness profiles as the two previous ones, as shown on Figure \ref{fig:scatter-plots-monochromatic-corrAdd}. However, the influence on the skewness is now opposite to what was observed in the case of correlated LU noise and decorrelated additive noise. As a matter of fact, Figure \ref{fig:comparison-CDF-monochromatic-corrAdd} shows that the empirically computed function $F$ of equations \eqref{eq:diffCDF} is now consistently negative for positive values of $x$. This is coherent with Figure \ref{fig:hist-skew-kurt-corrAdd}, where the values of skewness and kurtosis are lower overall for the dataset with LU decorrelated noise than without.

The comparison between correlated and decorrelated noise cases suggest that the LU formalism captures better non-linear dynamics when the subgrid phenomena it models is more correlated in time than the forcing, typically by the wind in our study. Conversely, a highly decorrelated subgrid noise would likely induce small-scale mixing phenomena, which in turn yields dissipation and less non-linearity. Hence, LU can be seen a new degree of freedom, which allows to tune ``how well'' the system captures non-linearity. In this tuning, the time-correlation of the noise appears to be a very impactful parameter.


\subsubsection{An example of (decorrelated) polychromatic noises}

\begin{figure}
    \centering
    \includegraphics[trim={0 0 0 0}, clip, width=6cm]{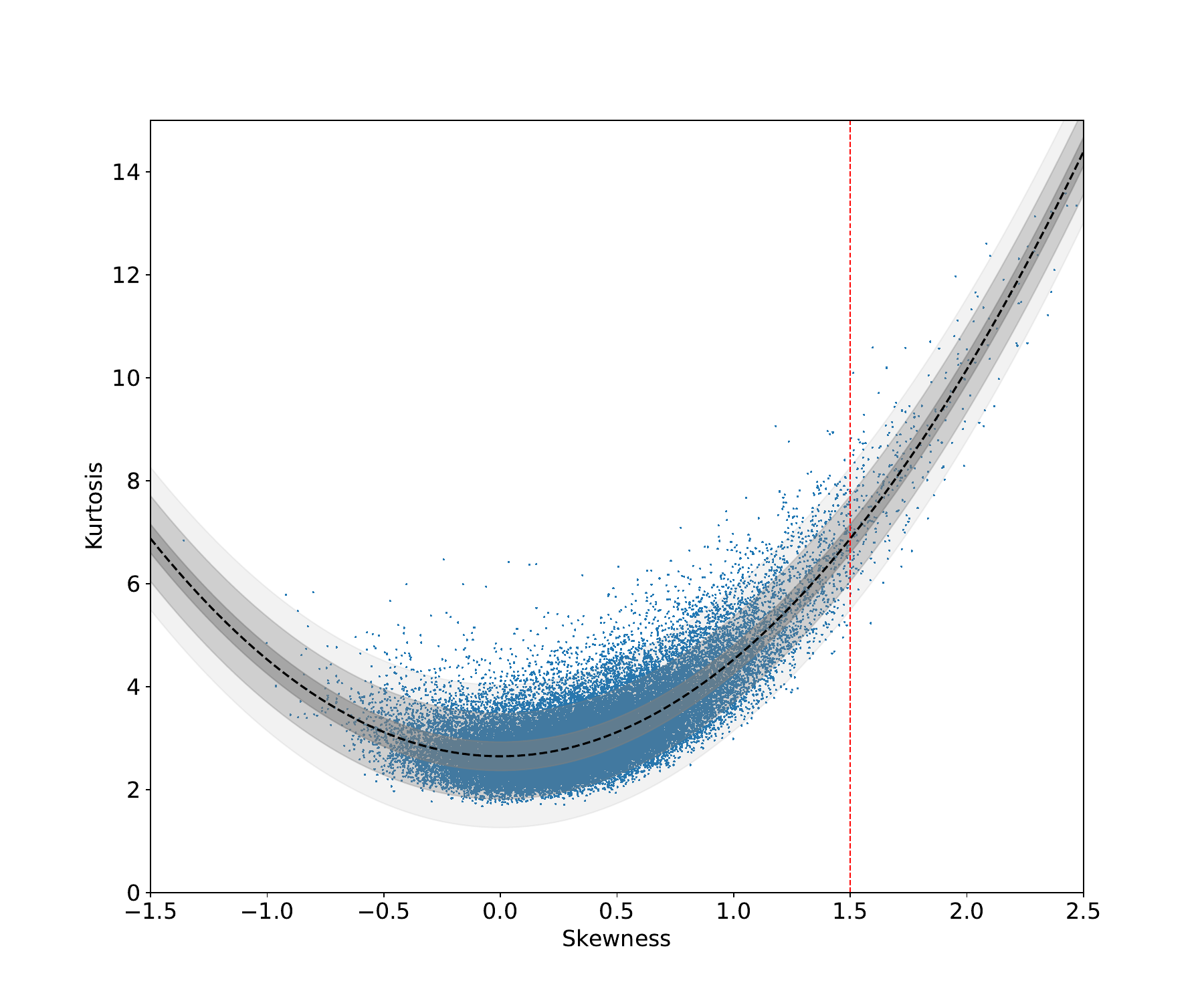}
    \includegraphics[trim={0 0 0 0}, clip, width=6cm]{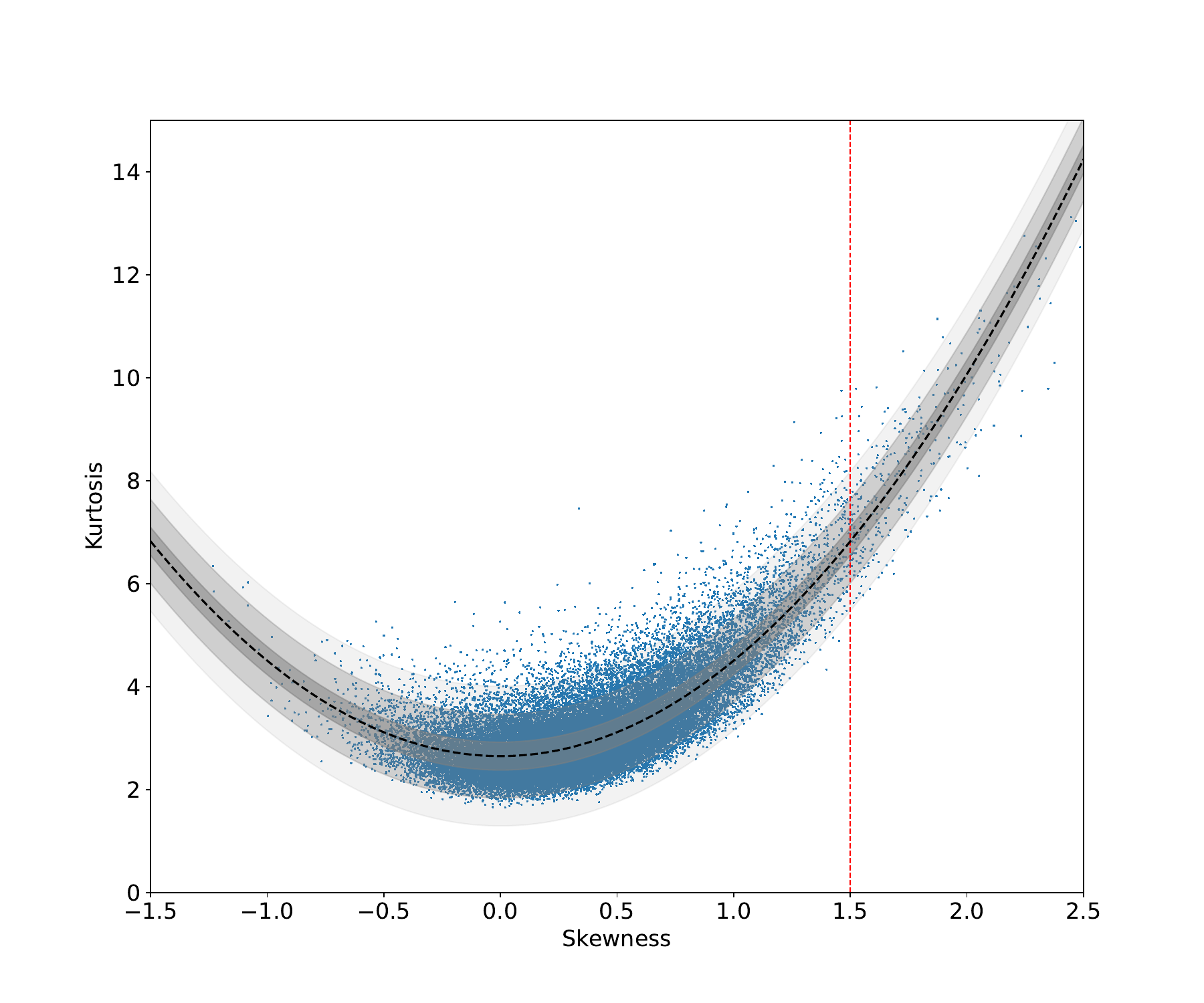}
    \caption[Scatter plots of the LU Saint-Venant system with Ornstein-Ulhenbeck noise for two datasets.]{Scatter plot of Dataset 3.1 (left) and Dataset 3.2 (right) in skewness-kurtosis axes, similar to Figure \ref{fig:scatter-plots-monochromatic}. Total number of realisations: 41984. Number of realisations with $(\mu_3 \geq 1.5)$: 298 (Dataset 3.1) and 304 (Dataset 3.2). Parameters for Dataset 3.1: $A=1.37$, $B=2.65$. Parameters for Dataset 3.2: $A=1.36$, $B=2.65$.}
    \label{fig:scatter-plots-monochromatic-polyChrom}
\end{figure}

\begin{figure}
    \centering
    \includegraphics[trim={0 2cm 0 2cm}, clip, width=12cm]{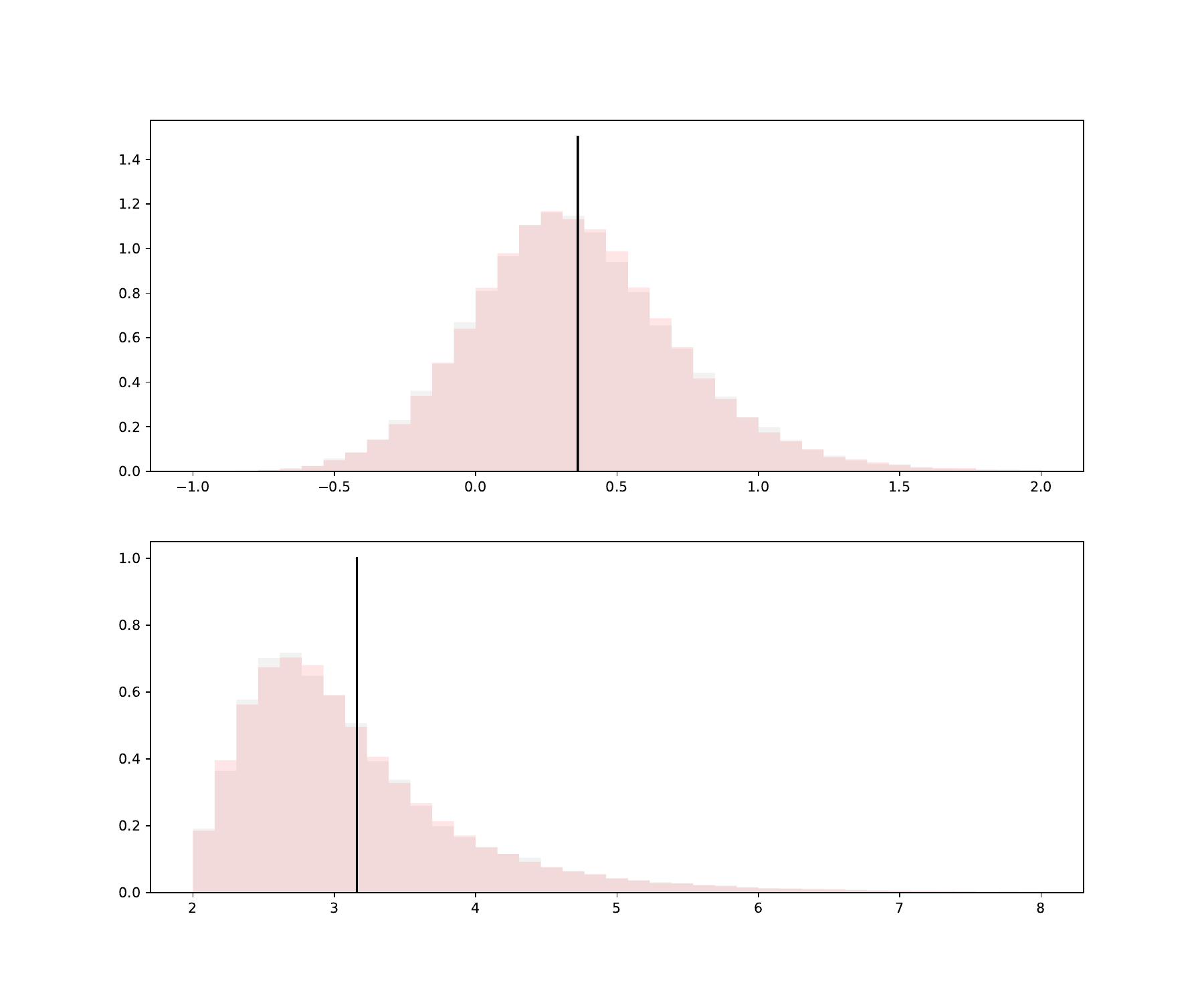}
    \caption[Distribution in skewness and kurtosis axes for two datasets.]{Distribution of skewness (upper) and kurtosis (lower) for Dataset 3.1 (red) and 3.2 (grey). The means are represented by vertical lines.}
    \label{fig:hist-skew-kurt-polyChrom}
\end{figure}

\begin{figure}
    \centering
    \includegraphics[trim={0 0 0 0}, clip, width=11cm]{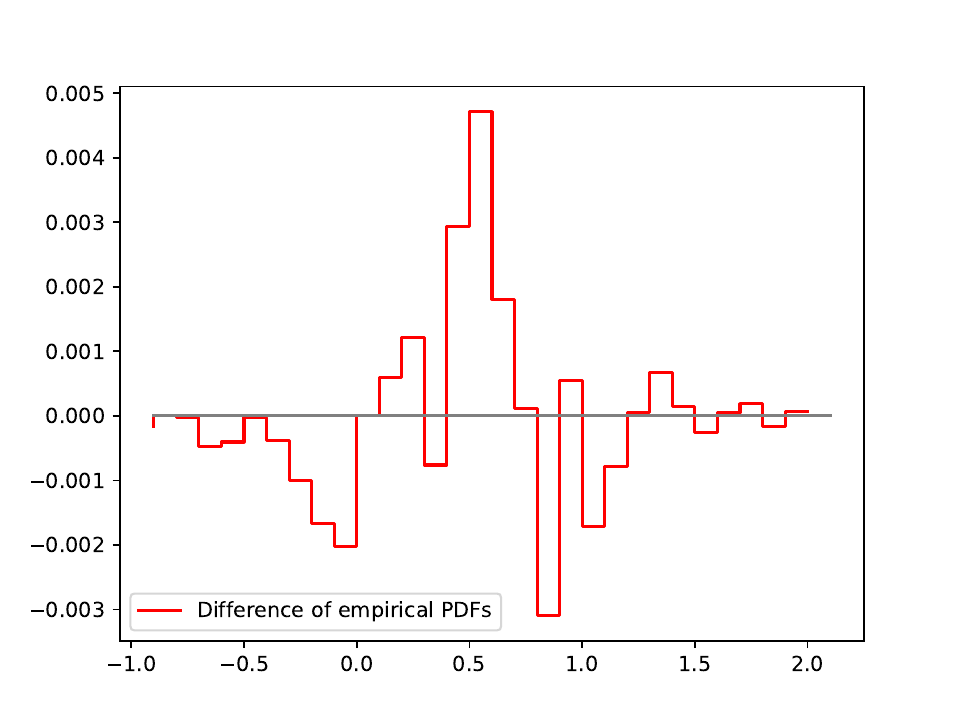}
    \caption[Representation of the empirically computed difference of probability density functions.]{Representation of the empirically computed derivative of the function $F$ as defined in \eqref{eq:diffCDF}. This corresponds to the difference of empirical PDFs (probability density functions).}
    \label{fig:comparison-PDF-monochromatic-polyChrom}
\end{figure}

\noindent
In this subsection, we commute back to decorrelated additive and LU noise, and we consider that they have a polychromatic spatial profile, which is expected to lead to wave-wave interaction and non-linear phenomena. Our setting is the same as \eqref{eq:forced-LUSV} -- with decorrelated noise -- but this time with the following additive and LU noises,
\begin{align*}
    \sigma dW_t^1 = \sqrt{\Upsilon_{LU}} f dW_t^1, \quad \sigma_{Add} dW_t^2 = \sqrt{\Upsilon_{Add}} f dW_t^2,
\end{align*}
where
\begin{equation}
    f = \frac{1}{1+(\frac{k_1}{k_2})^{5/2} + (\frac{k_1}{k_3})^{5/2}} \Big(f_{\alpha}^{k_1} + \big(\frac{k_1}{k_2}\big)^{5/2} f_{\alpha}^{k_2} + \big(\frac{k_1}{k_3}\big)^{5/2} f_{\alpha}^{k_3}\Big).
\end{equation}
The noise modes are weighted as in \cite{MLLTC_2023}, that is with a decreasing power law of exponent $5/2$. Additionally, we set $\Upsilon_{LU} = 1.0 \: 10^7$ and $\Upsilon_{Add} = 5.0 \: 10^4$. The wavenumbers are chosen to be $k_1 = 2\pi$, $k_2 = 5\pi$ and $k_3 = 10\pi$. As before, we set the coefficient $\alpha = 10^1$, and generate two datasets in this setting.
\begin{itemize}
    \item[$\bullet$] Dataset 3.1: 50 trajectories are generated in the setting presented above.
    \item[$\bullet$] Dataset 3.2: 50 trajectories are generated in the same setting, except that $\Upsilon_{LU} = 0$. It corresponds to a (deterministic) Saint-Venant system that is forced by a (decorrelated and polychromatic) additive noise.
\end{itemize}
Figure \ref{fig:scatter-plots-monochromatic-polyChrom} shows that considering this setting does not affect much the parameters $A$ and $B$. Moreover, it is now unclear from the comparison of the skewness and kurtosis distributions whether the presence of the LU noise favours or penalises non-Gaussian events -- see Figure \ref{fig:hist-skew-kurt-polyChrom}. This is supported by Figure \ref{fig:comparison-PDF-monochromatic-polyChrom}, which shows that the difference of empirical PDFs (probability density functions) does not exhibit a a remarkable behaviour. In turn, the polychromatic noise setting appears not to have a strong influence on the dynamics overall. Hence, it suggests that the space frequency profile has a less impact than the time correlation of the noise which was previously presented. Yet, considering more general noises stands as a natural approach to generate other wave models, where many parameters of interest -- such as wavenumber profiles, noise amplitudes, correlation lengths etc. -- can be tuned to represent the dynamics faithfully. This can be the subject of further studies.

\section{Conclusion}

\noindent
Hence, a stochastic interpretation of the deterministic Water Waves can be inferred from Euler equations under location uncertainty. New interpretations of the wave models studied in \cite{DMM_2024} can then be derived in the limit of small aspect ratios. These are very similar to what was previously established, but involve extra ``commutation'' terms, which can be neglected in suitable limits. Such interpretations open the door to deriving new stochastic versions of other classical models, as the Water Waves equations are central in wave modelling. Additionally, the stochastic Water Waves system thus inferred allows incorporate rotationality into the system through the noise. This can be seen, to a limited extent, as a way to bridge the study of currents and waves.

We have shown numerically that the LU noise alters significantly the shape of the waves, either in the Water Waves or the Saint-Venant systems. Moreover, it can induce higher and steeper waves when associated to an additive noise forcing. The time correlation of the noise also appears to be an important parameter, as it strongly influences the skewness and kurtosis of the distribution, while the polychromaticity of the noise seems to have little impact on the dynamics. On the one hand, high time correlation of the LU noise compared to the additive one seems to favour high skewness events, which are interpreted as non-linear phenomena. On the other hand, decorrelated LU noise with correlated additive noise is associated to fewer of these events in our simulations. Yet, more simulations would be needed to fully assess the properties of these wave models. We could also choose different types of noise, by taking a smaller parameter $\alpha$ to cover a wider range of frequencies or consider various polychromatic profiles for instance. Performing a statistical study of the Water Waves system would most likely provide insightful results as well, since it is the most general model studied here. Furthermore, investigating the Boussinesq system should be simultaneously staightforward and enlightening, as it is a relatively simple generalisation of the Saint-Venant system, and  yields natural regularisation properties. Although it does not ensure the conservation of energy, the Boussinesq system often remains a good approximation of the much more complex Serre-Green-Naghdi system, and thus must induce more complex dynamics than the ones we have studied here.


\bibliography{sn-bibliography}

\end{document}